\documentclass[11pt,a4paper, twoside]{article}

\usepackage{amsmath,amssymb,amsthm,amsfonts,mathrsfs,amscd,environ}
\usepackage{dsfont}
\usepackage{latexsym,enumerate,color,geometry,extarrows}

\def\beg{\begin}
\def\bequ{\begin{equation}}
\def\enqu{\end{equation}}
\def\bes{\begin{split}}
\def\ens{\end{split}}
\def\bews{\begin{ews}}
\def\beqn{\begin{eqnarray}}
\def\enqn{\end{eqnarray}}
\def\beq*{\begin{equation*}}
\def\enq*{\end{equation*}}
\def\bqn*{\begin{eqnarray*}}
\def\eqn*{\end{eqnarray*}}
\def\bary{\begin{array}}
\def\eary{\end{array}}
\def\bpma{\begin{pmatrix}}
\def\epma{\end{pmatrix}}
\def\bvma{\begin{Vmatrix}}
\def\evma{\end{Vmatrix}}

 \numberwithin{equation}{section}

\def\al{\alpha}
\def\be{\beta}
\def\ga{\gamma}
\def\de{\delta}
\def\ep{\epsilon}

\def\th{\theta}

\def\si{\sigma}

\def\ph{\phi}

\def\ps{\psi}

\def\Ga{\Gamma}

\def\Th{\Theta}

\def\R{\mathbb R}
\def\P{\mathbb P}
\def\E{\mathbb E}

\def\sF{\mathscr F}

\def\sB{\mathscr B}

\def\sK{\mathscr K}

\def\sG{\mathscr G}

\def\d{\mathrm{d}}

\def\ff{\frac}
\def\ra{\rightarrow}
\def\nn{\nabla}

\def\<{\langle}
\def\>{\rangle}
\def\sq{\sqrt}
\def\tld{\tilde}
\def\we{\wedge}
\def\1{\mathds{1}}

\allowdisplaybreaks

\title{{\bf  A unified approach to gradient type formulas for BSDEs and some applications}
}

\author{
{\bf Xiliang Fan$^{a),b)}$,  Michael R\"{o}ckner$^{b),c)}$, Shao-Qin Zhang$^{d)}$}\\
\footnotesize{$^{a)}$School of Mathematics and Statistics, Anhui Normal University, Wuhu 241002, China}\\
\footnotesize{$^{b)}$Fakult\"{a}t f\"{u}r Mathematik, Universit\"{a}t Bielefeld, 33615 Bielefeld, Germany}\\
\footnotesize{$^{c)}$Academy of Mathematics and Systems Science, CAS, Beijing 100190, China}\\
\footnotesize{$^{d)}$School of Statistics and Mathematics, Central University of Finance and Economics, Beijing 100081, China}\\
\footnotesize{\sf{fanxiliang0515@163.com}(X. Fan),\ \sf{roeckner@math.uni-bielefeld.de}(M. R\"{o}ckner),\ \sf{zhangsq@cufe.edu.cn}(S. Zhang)}\\
}

\begin{document}

\maketitle

\begin{abstract}
In this paper we present a unified approach to establish gradient type formulas and Bismut type formulas for backward stochastic differential equations (BSDEs).
This approach relies on a mix of derivative formulas with respect to the conditional probability of forward SDEs and the expression of the solution of  BSDEs.
Some concrete examples are given to illustrate the results.
As applications, we provide representation formulas for the control solutions to McKean-Vlasov BSDEs and derive gradient estimates for related PDEs.
\end{abstract}\noindent

\noindent

\vskip 2cm

\section{Introduction}

Let $(\Omega,\sF,(\sF_t)_{t\in[0,T]},\P)$ be a filtered probability space
with $(\sF_t)_{t\in[0,T]}$ the natural completed and right continuous filtration generated by an $m$-dimensional Brownian motion $(W_t)_{0\leq t\leq T}$.
On $(\Omega,\sF,(\sF_t)_{t\in[0,T]},\P)$ we consider the following system of equations:
\beg{align}\label{Intro-1.1}
Y_t=g(X_T)+\int_t^Tf(r,X_r,Y_r,Z_r)\d r-\int_t^TZ_r\d W_r, \ \ \ t\in[0,T],
\end{align}
where $X$ is the solution to a stochastic differential equation (SDE):
\beg{align}\label{Intro-1.2}
X_t=x+\int_0^tb(r,X_r)\d r+\int_0^t\sigma(r,X_r)\d W_r, \ \ \ t\in[0,T].
\end{align}
Here, the coefficients $g:\R^d\rightarrow\R^l,f:[0,T]\times\R^d\times\R^l\times\R^l\otimes\R^m\rightarrow\R^l$
and $b:[0,T]\times\R^d\rightarrow\R^d, \sigma:[0,T]\times\R^d\rightarrow\R^d\otimes\R^m$ are deterministic functions.
The system \eqref{Intro-1.1}-\eqref{Intro-1.2} is called a (decoupled) forward-backward SDE (FBSDE),
in which the processes $X$ and $Y$ are called the forward component and the backward component, respectively.
The problem of existence and uniqueness for systems of this kind was first addressed by Pardoux and Peng \cite{Pardoux&Peng90a,Pardoux&Peng92a},
and since then there are a large number of papers dedicated to the study of FBSDEs due to their increasing importance in stochastic control and
mathematical finance (see, e.g., \cite{Karoui&Peng&Quenez97a,Ma&Protter&Yong94a,Ma&Yong99a,Yong&Zhou99a,JianfengZhang}).

In \cite{Pardoux&Peng92a}, Pardoux and Peng proved that the stochastic flow $(X^{t,x},Y^{t,x},Z^{t,x}),t\in[0,T],x\in\R^d$,
which is the solution to the system \eqref{Intro-1.1}-\eqref{Intro-1.2} restricted to the interval $[t,T]$ with $X_t^{t,x}=x$,
provides a nonlinear Feynman-Kac formula for the solution to a quasi-linear parabolic partial differential equation (PDE) of the form:
\beg{equation}\label{Intro-1.3}
\beg{cases}
\partial_t u(t,x)+\mathbb{L}u(t,x)+f(t,x,u(t,x),(\nn u\sigma)(t,x))=0,\ \ (t,x)\in[0,T]\times\R^d,\\
u(T,x)=g(x),\ \ x\in\R^d,
\end{cases}
\end{equation}
where $\mathbb{L}u=(\mathcal{L}u_1,\cdots,\mathcal{L}u_l)^*$ with
$\mathcal{L}u_i=\frac{1}{2}\mathrm{Tr}(\sigma\sigma^*\nn^2u_i)+\langle b,\nn u_i\rangle, i=1,\cdots,l$, and the notation $^*$ denotes the transpose.
More precisely, if the coefficients $g,f$ and $b,\sigma$ are sufficiently smooth, then
\beg{align}\label{Intro-1.4}
u(t,x)=Y_t^{t,x},\ \ \ (t,x)\in[0,T]\times\R^d
\end{align}
in the classical solution sense of PDE \eqref{Intro-1.3}.
Conversely, the solution $u$ to PDE \eqref{Intro-1.3} admits the following representation formula for the solution to
backward stochastic differential equation (BSDE) \eqref{Intro-1.1}:
\beg{align}\label{Intro-1.5}
Y_s^{t,x}=u(s,X_s^{t,x}),\ \ \ Z_s^{t,x}=\nn u(s,X_s^{t,x})\sigma(s,X_s^{t,x}), \ \ \ s\in[t,T],
\end{align}
which is also valid for fully coupled FBSDEs, i.e. $b$ and $\sigma$ may depend on $(Y,Z)$ (see \cite{Ma&Protter&Yong94a} for further details).
If the coefficients $g,f$ and $b,\sigma$ are Lipschitz continuous,
the relations \eqref{Intro-1.4} and  \eqref{Intro-1.5} between BSDE \eqref{Intro-1.1} and PDE \eqref{Intro-1.3}
remain true in the viscosity solution sense of PDE \eqref{Intro-1.3} (see, e.g., \cite{Bally&Pardoux&Stoic05a,Pardoux99a,Pardoux&Peng92a}).
We also point out that Crisan and Delarue \cite{Crisan&Delarue12a}, Zhu \cite{Zhu12a,Zhu12b} extended the above relations
to the case of generalized solution of PDE \eqref{Intro-1.3}.

When dealing with applications especially in the numerical analysis of BSDE,
one needs to investigate regularity properties of the solution $(Y,Z)$.
In \cite{MaZ}, Ma and Zhang first studied the path regularity for the $Z$ component when $g,f$ are Lipschitz continuous
and $\sigma$ is uniformly non-degenerate.
The key for their approach is to establish a probabilistic representation formula for $\nn u$
and then $Z$ by using Malliavin's integration by parts formula (see Remark \ref{Rem-Ex1} and Remark \ref{Rem1-Th3.1} (i) below).
That is, for any $(t,x)\in[0,T)\times\R^d$,
\beg{align}\label{Intro-1.6}
\nn u(t,x)=\E\left(g(X_T^{t,x})N_T^t+\int_t^T f(r,X_r^{t,x},Y_r^{t,x},Z_r^{t,x})N_r^t\d r\right),
\end{align}
where $N_r^t=\frac{1}{r-t}\int_t^r\sigma^{-1}(s,X_s)\nabla X_s\d W_s\cdot(\nabla X_t)^{-1}$.
Afterwards, Zhang \cite{ZJF} proved the path regularity of $Z$ in the degenerate case ($\sigma\neq0$)
where all processes concerned are one dimensional, $f$ is linear on $Z$ and $g$ maybe discontinuous.
The argument consists of using a representation formula for $\nn u$ similar to \eqref{Intro-1.6},
in which $N_\cdot^t\nabla X_t$ is a bounded variation process rather than a martingale
that makes the estimates of $\nn u$ more complicated.
We remark that although the conditions imposed on $\sigma$ in \cite{ZJF} don't require the invertibility for $\sigma\sigma^*$,
the invertibility of a linear functional of $\sigma\sigma^*$ is required instead.
As for another type of regularity, namely distributional regularity,
we refer to, e.g., \cite{Aboura&Bourguin13a,Antonelli&Kohatsu-Higa05a,Fan&Wu,Mastrolia16a,Mastrolia&Possama&Reveillac16a,Olivera&hamarova}
and the references therein.

The aim of this paper is to establish a gradient type formula for the solution $Y$
of the system \eqref{Intro-1.1}-\eqref{Intro-1.2}, which is stronger than \eqref{Intro-1.6}.
Our main result in that direction is Theorem \ref{thm-BSDE} below,
using a fundamentally different approach from \cite{MaZ,ZJF} and under more general assumptions.
We first establish a gradient type formula in the sense of conditional expectation for Markov processes
that is stronger than the Bismut type formula.
To the best of our knowledge the result is new, and its proof relies heavily on the Markov property.
Then, with the help of the relation \eqref{Intro-1.5} obtained under weaker conditions,
we are able to provide a gradient type formula for the associated FBSDE in a general setting within which the formula can be applied to a larger class of examples.
As a product of our main result, we establish a Bismut type formula for FBSDEs,
which generalize and improve the corresponding ones in the existing literature (see Corollary \ref{cor2} and Remark \ref{Rem-Ex1} below).
In addition, we discuss a series of examples,

$\bullet$ FBSDEs with non-degenerate forward SDEs, for which two different types of gradient type formulas are given,

$\bullet$ FBSDEs with forward Gruschin type processes,

$\bullet$ FBSDEs with forward stochastic Hamiltonian systems.\\
These examples should illustrate the power and flexibility of our unified method.
We believe that the method can also be used to handle other types of FBSDEs.

As an application of the main results, we wish to study McKean-Vlasov BSDEs and related nonlocal PDEs.
McKean-Vlasov SDEs (whose coefficients depend on the law of the solution), initiated by Kac \cite{Kac} and also known as mean-field equations, have been studied extensively in the past decade
as they naturally arise in the context of statistical physics,
and as they provide probabilistic representations for solutions of a class of nonlinear PDEs
which may involve the Lions derivative introduced by Lions in his lectures \cite{Cardaliaguet13}.
Recently, Buckdahn, Djehiche, Li and Peng \cite{Buckdahn&Djehiche&Li&Peng09a} first investigated nonlinear McKean-Vlasov BSDEs.
Since then, existence and uniqueness results of McKean-Vlasov FBSDEs and the theory of the associated nonlocal PDEs of mean-field type
have been studied in a variety of settings.
For examples, Chassagneux, Crisan and Delarue \cite{Chassagneux&Crisan&Delarue15} proved the existence and uniqueness of solutions to fully coupled McKean-Vlasov FBSDEs.
Carmona and Delarue \cite{Carmona&Delarue15} investigated McKean-Vlasov FBSDEs through the stochastic maximum principle;
Li \cite{Li18a} considered McKean-Vlasov FBSDEs driven by a Brownian motion and an independent Poisson random measure
and showed the existence and uniqueness of solutions of this kind and provided a
solution for the related nonlocal quasi-linear integral PDE of mean-field type via the value function.
The second objective of the present paper is to deepen the investigation of McKean-Vlasov FBSDEs.
With the main results above in hand,
we shall establish representation formulas for the control solution $Z^{t,\xi}_s$ of a McKean-Vlasov FBSDE which enable us to derive its path regularity,
and obtain gradient estimates for the solution to the associated nonlocal PDE of mean-field type.

The remainder of this paper is organized as follows.
In Section 2, we state and prove our main results concerning the gradient type formulas and Bismut type formulas of FBSDEs,
which are then applied to concrete FBSDEs associated with various types of forward SDEs such as non-degenerate SDEs,
Gruschin type processes and stochastic Hamiltonian systems.
In Section 3, using these formulas,
we establish representation formulas for the control solutions of McKean-Vlasov FBSDEs,
and provide gradient estimates for the corresponding PDEs.
Section 4 will be devoted to the proofs of some auxiliary lemmas.

$\mathbf{Notation.}$ The following notations are used in the sequel.

$\bullet$ We use $|\cdot|$ and $\<\cdot,\cdot\>$ for the Euclidean norm and the Euclidean inner product, respectively, and $\|\cdot\|$ for either the operator norm or the Euclidean norm if this will not lead to confusion.
For any $p\in[1,\infty)$, let $\|\cdot\|_p$ denote the $L^p(\P):=L^p(\Omega\ra\R^d,\sF,\P)$ norm,
and $\mu\text{-}\lim$ denotes the limit in measure $\mu$.
For any $0\leq t\leq r\leq T$,
$\sF_r^t$ denotes the completion of $\si\{W_s-W_t: t\leq s\leq r\}$ with the $\P$-null sets of $\sF$.

$\bullet\ \sB(\R^d;\R^l)$ denotes the space of all Borel measurable functions $\varphi:\R^d\rightarrow\R^l$
and $\sB_b(\R^d;\R^l)$ denotes the space of all those $\varphi\in\sB(\R^d;\R^l)$ that are bounded on $\R^d$
with the norm $\|\varphi\|_\infty:=\sup_{x\in\R^d}|\varphi(x)|$.
$C_b(\R^d;\R^l)$ is the set of all bounded continuous functions $\varphi:\R^d\rightarrow\R^l$.
$C^1(\R^d;\R^l)$ denotes the collection of all continuously differentiable functions $\varphi:\R^d\rightarrow\R^l$
and $C_b^1(\R^d;\R^l)$ denotes the collection of all those $\varphi\in C^1(\R^d;\R^l)$ with bounded derivatives.
$C^{0,1}([0,T]\times\R^d;\R^l)$ stands for the class of all continuous functions $\varphi:[0,T]\times\R^d\rightarrow\R^l$
such that they are continuously differentiable on the space variable.
For $f\in C^1(\R^d;\R^l)$, let $\nn f(x)=(\nn f_1(x),\cdots,\nn f_l(x))^*\in\R^l\otimes\R^d$ denote the gradient matrix.
When $l=1$, we often suppress $\R^l$ for simplicity.

$\bullet$ We let $C, C_{k_1,k_2,T}, C_{T,q_1,q_2}$, etc., denote generic constants, whose values may change from line to line.

\section{Gradient type formulas for FBSDEs}

The main objective of this section concerns the problem of gradient type formulas for FBSDEs.
We first construct gradient type formulas for general Markov processes.
In the second part of this section, we will show how to combine these results to establish gradient type formulas and Bismut type formulas for FBSDEs,
and their proofs are addressed in the third part.
Finally, we present some concrete examples to illustrate our main results.

\subsection{The case of Markov processes}

Let $\{X^t_r\}_{0\leq t\leq r\leq T}$ be a Markov process with respect to the filtration $\{\sF_r^t\}_{0\leq t\leq r\leq T}$.
We also write $X^{t,x}_r$ if  $X^{t}_t=x$.
Let $P_{t,r}\ps(x)=\E \ps(X^{t,x}_r)$,~$\ps\in\sB(\R^d)$, be the associated Markov semigroup.
The current part is dedicated to gradient type formulas for $X^{t,x}_r$ under the following assumptions on the $\{X^{t,x}_r\}_{0\leq t\leq r\leq T}$:
\beg{enumerate}
\item[(H1)] There exists a $\R^d\otimes\R^d$-valued adapted  process  $\{\nn  X_r^{t,x}\}_{0\leq t<r\leq T}$ such that
\beg{align}\label{L2-de}
\lim_{\ep\ra 0^+}\E\left| \ff {X_r^{t,x+\ep v}-X_r^{t,x}} {\ep}-\nn_v X_r^{t,x}\right|^2=0,\ \ v\in\R^d,
\end{align}
where $\nn_v X_r^{t,x}=(\nn X_r^{t,x})v$.
\item[(H2)]  For any $0\leq t<r\leq T$ and $\ps\in C^1_b(\R^d)$, $P_{t,r}\ps\in C^1(\R^d)$
and for any $ x\in\R^d$ there exists an integrable $\R^d$-valued, $\sF_r^t$-measurable random variable $ M^t_r(x) $ such that
\beg{align}\label{Bism-P}
\nn_v P_{t,r}\ps(x)=\E [ \ps(X^{t,x}_r)\<M_r^t(x), v\>],\ \ \ps\in C_b^1(\R^d),\ v\in\R^d.
\end{align}
\item[(H3)] For any $0\leq t\leq s<r\leq T$,
\beg{align}
&\lim_{y\ra  x}\E\left|M_r^s(X_s^{t,y})-M_{r}^s(X_s^{t,x})\right|^2=0,\label{L2-M}\\
&\lim_{\de\ra 0^+}\E\left|M_r^s(X_s^{t,x})-M_{r-\de}^s(X_s^{t,x})\right|^2=0.\label{add-Hy}
\end{align}
\end{enumerate}

With $M_r^t(\cdot)$ in hand, we can provide a representation formula for $\nn_v\E[\ph(X_r^{t,x})|\sF_s^t]$ defined as follows:
$$\nn_v\E[\ph(X_r^{t,x})|\sF_s^t]:=\lim\limits_{\epsilon\ra0^+}\ff {\E[\ph(X_r^{t,x+\epsilon v})|\sF_s^t]-\E[\ph(X_r^{t,x})|\sF_s^t]}{\epsilon},$$
if the limit exists in $L^p(\P)$ for some $p\geq 1$ (see Proposition \ref{lem-lpc} below,
which will play a crucial role in constructing gradient type formulas for FBSDEs in the next part).

\beg{prp}\label{lem-lpc}
Assume that (H1) and (H2) hold, and that for any $0\leq t\leq s<r\leq T$  there exist $k_0\geq 0$ and $C(t,s,r)>0$ such that
\beg{align}\label{in-K-pog}
K(t,s,r,x):=\|M_r^s(X^{t,x}_s)\|_2\leq C(t,s,r)(1+|x|^{k_0}),\ \ x\in\R^d,
\end{align}
and $\P$-a.s for any $\ps\in\sB_b(\R^d)$,
\beg{align}\label{equ-con}
\E\left[\ps(X_r^{t,x})G_r^s(X_s^{t,x})\big|\sF_s^t\right]=\left[\E\left(\ps(X_r^{s,y})G_r^s(y)\right)\right]\Big|_{y=X_s^{t,x}},
\ \  x,y\in\R^d,
\end{align}
where $G=M, |M|$.
Assume moreover that for any $k_1\geq 1$, there exist $k_2\geq 0$ and $C_{k_1,k_2,T}> 0$ such that
\beg{align}\label{mon-est-1}
\sup_{0\leq s\leq r\leq T}\|X_r^{s,y}\|_{k_1}\leq C_{k_1,k_2,T}(1+|y|^{k_2}),\ \ y\in\R^d.
\end{align}
\noindent (1) For any $\ph\in C(\R^d)$ with some $q\geq 0$ such that
\beg{align}\label{gg}
|\ph(x)|\leq C(1+|x|^q),\ \ x\in\R^d,
\end{align}
and for any $0\leq t\leq s<r\leq T$ and $x,v\in\R^d$, $\nn_v\E[\ph(X_r^{t,x})|\sF_s^t]$ exists in  $L^p(\P)$ for any $p\in [1,2)$ and moreover
\beg{align}\label{lpcov}
&\nn_v\E\left[\ph(X_r^{t,x})\big|\sF_s^t\right] =\E\left[ \ph(X_r^{t,x})\left\<M_r^s(X_s^{t,x}),\nn_v X_s^{t,x}\right\> \big| \sF_s^t\right].
\end{align}
\noindent (2) If (H3) holds, and if furthermore $C(t,s,\cdot)$ is bounded on any closed subinterval of $(s,T]$ and
\beg{align}\label{L2-con}
&\lim_{y\ra  x}\E\left\|\nn  X_s^{t,y} -\nn  X_s^{t,x} \right\|^2=0 ,
\end{align}
then for any  $\ph\in \sB(\R^d)$ satisfying  \eqref{gg},  \eqref{lpcov} also holds in $L^p(\P)$ for any $p\in [1,2)$.
\end{prp}

To prove this proposition, we first prove the following simple lemma.

\beg{lem}\label{lem-1}
Assume that $\{X^t_r\}_{0\leq t\leq r\leq T}$ satisfies (H2). Let $\ph\in \sB(\R^d)$  such that
\beg{align}\label{add-1(Lem2.2)}
\E \left[|\ph(X_r^{t,x})|(|M_r^t(x)|+1)\right] <\infty,\ \ 0\leq  t< r\leq T,\ x\in\R^d,
\end{align}
and suppose there exists a sequence $\{\ph_n\}_{n=1}^\infty\subset C_b^1(\R^d)$ such that for any $R>0$
\beg{align}\label{lim-sup}
\lim_{n\ra+\infty}\sup_{|x|\leq R}\E\left[|(\ph-\ph_n)(X^{t,x}_r)|\left(|M_r^t(x)|+1\right)\right] =0.
\end{align}
Then $P_{t,r}\ph\in C^1(\R^d)$ satisfying \eqref{Bism-P}. If in addition, \eqref{in-K-pog} and \eqref{mon-est-1} hold and $\ph$ satisfies \eqref{gg}, then  there exist $C>0$ and $q_0>0$ such that
\beg{align}\label{po-nn-ph}
\left|\nn P_{t,r}\ph(y)\right| \leq C(1+|y|^{q_0}),\ \ y\in\R^d.
\end{align}
\end{lem}

\beg{proof}
Since $\ph_n\in C_b^1(\R^d)$, (H2) yields that $P_{t,r}\ph_n\in C^1(\R^d)$. Moreover, for any $R>0$,  we have by (H2) and \eqref{lim-sup} that
\beg{align}\label{lim-phn}
&\lim_{n\ra+\infty}\sup_{|x|\leq  R}\left(\left|P_{t,r}\ph(x)-P_{t,r}\ph_n(x)\right|+\left|\E [\ph(X_r^{t,x})M_r^t(x)]-\nn P_{t,r}\ph_n(x)\right|\right)\nonumber\\
=&\lim_{n\ra+\infty}\sup_{|x|\leq  R}\left(\left|\E\ph(X_r^{t,x})-\E\ph_n(X_r^{t,x})\right|
+\left|\E [\ph(X_r^{t,x})M_r^t(x)]-\E [\ph_n(X_r^{t,x})M_r^t(x)] \right|\right)\nonumber\\
\leq& \lim_{n\ra+\infty}\sup_{|x|\leq R}\E \left[|(\ph-\ph_n)(X_r^{t,x})|\left(1+|M_r^t(x)|\right)\right]\nonumber\\
=&0.
\end{align}
Hence, $P_{t,r}\ph\in C^1(\R^d)$ and it satisfies \eqref{Bism-P}.

Furthermore, if \eqref{in-K-pog} and \eqref{mon-est-1} (with $k_1=2q$) hold and $\ph$ satisfies \eqref{gg},
then the H\"older inequality implies that
\beg{align}\label{in-po-nnP}
\left|\nn_vP_{t,r}\ph(y)\right|&= \left|\E\<\ph(X_r^{t,y})M_r^t(y),v\>\right|\nonumber\\
&\leq |v|\left(\E|\ph(X_r^{t,y})|^2\right)^{\ff 1 2}\left(\E|M_r^t(y)|^2\right)^{\ff 1 2}\nonumber\\
&\leq C|v|\left(\E(1+|X_r^{t,y}|^q)^2\right)^{\ff 1 2}K(t,t,r,y)\nonumber\\
&\leq C|v|(1+|y|^{qk_2})(1+|y|^{k_0})\nonumber\\
&\leq C(1+|y|^{q_0})|v|,
\end{align}
where $q_0=qk_2+k_0$.
\end{proof}

\noindent{\emph{\bf{Proof of Proposition \ref{lem-lpc}}}}

Let us first point out two facts for any $\ph\in\sB(\R^d)$ satisfying \eqref{gg}.

\textit{Fact (i):  The relation \eqref{equ-con} holds for all such $\ph$ $\P$-a.s}.
Indeed, similar to \eqref{in-po-nnP}, by \eqref{in-K-pog}, \eqref{mon-est-1} and \eqref{gg} we first have
\beg{align}\label{est-phM}
\E|\psi(X_r^{s,y})G_r^s(y)|\leq C(1+|y|^{q_0}),\ \ y\in\R^d
\end{align}
with $q_0=qk_2+k_0$, where $G=M, |M|$.
Then by the monotone class theorem we derive the desired result.

\textit{Fact (ii): There exists a bounded measurable sequence $\{\tld \ph_n\}_{n=1}^\infty$ such that \eqref{lim-sup} holds}.
Indeed, setting $\tld \ph_n:=(\ph\vee (-n))\we n$, we have
\beg{align}\label{in-b-unb}
&\E \left[|(\ph-\tld \ph_n)(X^{s,y}_r)|\left(|M_r^{s}(y)|+1\right)\right]\nonumber\\
&=\E \left[|\ph(X^{s,y}_r)|\1_{[|\ph(X^{s,y}_r)|\geq n]}\left(|M_r^{s}(y)|+1\right)\right]\nonumber\\
&\leq\left[\E\left(|\ph(X^{s,y}_r)|^{2}\1_{[|\ph(X_r^{s}(y))|\geq n]}\right)\right]^{1/2}\left[\E\left(|M_r^{s}(y)|+1\right)^2\right]^{\ff 1 2}.
\end{align}
On one hand, for any $\de>0$, we deduce by using \eqref{mon-est-1} (with $k_1=q(2+\de)$) and \eqref{gg} that
\beg{align}\label{add-1(Prop2.1)}
\E \left(|\ph(X^{s,y}_r)|^{2}\1_{[|\ph(X^{s,y}_r)|\geq n]}\right)& \leq \E\left[ |\ph(X^{s,y}_r)|^{2}\left(\ff {|\ph(X^{s,y}_r)|} n\right)^{\de}\1_{[|\ph(X^{s,y}_r)|\geq n]}\right]\cr
&\leq n^{-\de}\E |\ph(X^{s,y}_r)|^{2+\de}\cr
&\leq C n^{-\de}\E \left(1+|X^{s,y}_r|^{q(2+\de)}\right)\cr
&\leq C n^{-\de}\left(1+|y|^{q(2+\de)k_2}\right).
\end{align}
On the other hand, by the definition of the function $K(t,s,r,x)$, we get
\beg{align}\label{ineq-M}
\E |M_r^{s}(y)|^2 = \E  |M_r^{s}(X_s^{s,y})|^2 =K^2(s,s,r,y),
\end{align}
which implies that $\E (|M_r^{s}(y)|+1)^2$ is locally bounded in $y$.
Plugging \eqref{add-1(Prop2.1)} and \eqref{ineq-M} into \eqref{in-b-unb}, we obtain the second fact.

The rest of the proof is divided into three steps.

\textit{Step 1:} \textit{Claim: For $\ph\in\sB(\R^d)$ satisfying \eqref{gg} and such that
$P_{s,r}\ph$ belongs to $C^1(\R^d)$ and satisfies \eqref{Bism-P} and \eqref{po-nn-ph},
$\nn_v\E[\ph(X_r^{t,x})|\sF_s^t]$ exists in $L^p(\P)$ for any $p\in [1,2)$ and moreover  \eqref{lpcov} holds.}
Since $P_{s,r}\ph\in C^1(\R^d)$ and due to the Markov property of $X$, we first obtain for $\ep>0, v\in\R^d$ and $0\leq t\leq s<r\leq T$,
\beg{align*}
&\ff 1 {\ep} \left(\E[\ph(X_r^{t,x+\epsilon v})|\sF_s^t]-\E[\ph(X_r^{t,x})|\sF_s^t]\right)-\nn_{\nn_vX_s^{t,x}}P_{s,r}\ph(X_s^{t,x}) \\
&=\ff 1 {\ep} \left(P_{s,r}\ph(X_s^{t,x+\ep v})-P_{s,r}\ph(X_s^{t,x})\right)-\nn_{\nn_vX_s^{t,x}}P_{s,r}\ph(X_s^{t,x}) \\
&=\int_0^1 \nn_{\ff {X_s^{t,x+\ep v}-X_s^{t,x}} {\ep }}P_{s,r}\ph(X_s^{t,x}+\th (X_s^{t,x+\ep v}-X_s^{t,x}) )\d \th-\nn_{\nn_vX_s^{t,x}}P_{s,r}\ph(X_s^{t,x}) \\
&=\left\<\int_0^1\nn P_{s,r}\ph(X_s^{t,x}+\th (X_s^{t,x+\ep v}-X_s^{t,x}) )\d \th-\nn P_{s,r}\ph(X_s^{t,x}), \nn_v X_s^{t,x}\right\>\\
&\quad +\left\<\int_0^1 \nn P_{s,r}\ph(X_s^{t,x}+\th (X_s^{t,x+\ep v}-X_s^{t,x}) )\d \th, \ff {X_s^{t,x+\ep v}-X_s^{t,x}} {\ep}-\nn_vX_s^{t,x} \right\>.
\end{align*}
Then in order to prove that for any $p\in[1,2)$,
\beg{align}\label{de-ph}
\lim_{\ep\ra 0^+}\E\left|\ff 1 {\ep} \left(\E[\ph(X_r^{t,x+\epsilon v})|\sF_s^t]-\E[\ph(X_r^{t,x})|\sF_s^t]\right)
-\nn_{\nn_v X_s^{t,x}} P_{s,r}\ph(X_s^{t,x})\right|^p=0,
\end{align}
according to the H\"{o}lder inequality, (H1), \eqref{po-nn-ph} and \eqref{mon-est-1},
it is enough to show that,
\beg{align}\label{L2pp}
\lim_{\ep\ra 0^+}\int_0^1 \nn P_{s,r}\ph(X_s^{t,x}+\th (X_s^{t,x+\ep v}-X_s^{t,x}) )\d \th=\nn P_{s,r}\ph(X_s^{t,x})
\end{align}
in $L^{\ff {2p} {2-p}}(\P)$ for any $p\in [1,2)$.\\
Noting that $\nn P_{s,r}\ph\in C(\R^d)$, and that
$$\P\text{-}\lim\limits_{\ep\ra 0^+}X_s^{t,x+\ep v}=X_s^{t,x},$$
we deduce that
\beg{align}\label{cov-p-nnP}
\P\text{-}\lim_{\ep\ra 0^+}\nn P_{s,r}\ph(X_s^{t,x}+\th (X_s^{t,x+\ep v}-X_s^{t,x}) )=\nn P_{s,r}\ph(X_s^{t,x}),\ \ \th\in [0,1].
\end{align}
For any $\tld p\geq 1$, it follows by \eqref{po-nn-ph} and \eqref{mon-est-1} that
\beg{align*}
&\sup_{\ep\in (0,1],\th\in [0,1]}\E\left|\nn P_{s,r}\ph(X_s^{t,x}+\th (X_s^{t,x+\ep v}-X_s^{t,x}) )\right|^{\tld p}\\
&\qquad\leq C\sup_{\ep\in(0,1]}\E \left(1+(|X_s^{t,x}|+|X_s^{t,x+\ep v}|)^{q_0}\right)^{\tld p}\\
&\qquad <\infty.
\end{align*}
Then we conclude that
$$\left\{\left|\nn P_{s,r}\ph(X_s^{t,x}+\th (X_s^{t,x+\ep v}-X_s^{t,x}) )\right|^{\ff {2p} {2-p}}\right\}_{\ep\in (0,1]}$$
is uniformly integrable with respect to $\d \th\otimes \d \P$.
Combining this with \eqref{cov-p-nnP}, by the dominated convergence theorem we obtain
\beg{align*}
&\lim_{\ep\ra 0^+}\E\left|\int_0^1 \left(\nn P_{s,r}\ph(X_s^{t,x}+\th (X_s^{t,x+\ep v}-X_s^{t,x}) )-\nn P_{s,r}\ph(X_s^{t,x})\right)\d \th \right|^{\ff {2p} {2-p}}\\
&\leq\lim_{\ep\ra 0^+}\int_0^1\E\left| \nn P_{s,r}\ph(X_s^{t,x}+\th (X_s^{t,x+\ep v}-X_s^{t,x}) )-\nn P_{s,r}\ph(X_s^{t,x})\right|^{\ff {2p} {2-p}}\d \th \\
&=\int_0^1\E\left(\lim_{\ep\ra 0^+}\left| \nn P_{s,r}\ph(X_s^{t,x}+\th (X_s^{t,x+\ep v}-X_s^{t,x}) )-\nn P_{s,r}\ph(X_s^{t,x})\right|^{\ff {2p} {2-p}}\right)\d \th\\
&=0,
\end{align*}
which means that \eqref{L2pp} holds in $L^{\ff {2p} {2-p}}(\P)$.\\
Hence, we obtain that $\nn_v\E[\ph(X_r^{t,x})|\sF_s^t]$ exists in $L^p(\P)$ for any $p\in [1,2)$,
and moreover by \eqref{de-ph} and \eqref{Bism-P} the directional derivative satisfies
\beg{align}\label{equ-1}
\nn_v \E\left[\ph(X_r^{t,x})\Big|\sF_s^t\right]&=\nn_{\nn_v X_s^{t,x}} P_{s,r}\ph(X_s^{t,x})\nonumber\\
&=\E\left[\ph(X_r^{s,y})\<M_r^s(y),w\>\right]\Big|_{y=X_s^{t,x},w=\nn_v X_s^{t,x}}\nonumber\\
&=\left\<\E \left[\ph(X_r^{t,x}) M_r^s(X_s^{t,x})\Big|\sF_s^t\right],\nn_v X_s^{t,x}\right\>,
\end{align}
where we have used fact (i) in the last equality. \\
Note that again due to Fact (i), \eqref{est-phM} and \eqref{mon-est-1}, we obtain by the H\"{o}lder inequality that
\beg{align*}
\E \left|\ph(X_r^{t,x})\left\< M_r^s(X_s^{t,x}),\nn_v X_s^{t,x}\right\>\right|
&\leq \E \left[|\ph(X_r^{t,x})M_r^s(X_s^{t,x})|\cdot|\nn_v X_s^{t,x}|\right] \\
& =\E\left\{\E\left[ |\ph(X_r^{t,x})M_r^s(X_s^{t,x})|\Big|\sF_s^t\right]\cdot|\nn_v X_s^{t,x}|\right\}\\
& =\E \left\{\E\left[|\ph(X_r^{s,y})M_r^s(y)|\right]\Big|_{y=X_s^{t,x}}\cdot|\nn_vX_s^{t,x}|\right\}\\
& \leq C\left(\E\left(1+|X_s^{t,x}|^{q_0}\right)^2\right)^{\ff 1 2}\left(\E \left|\nn_v X_s^{t,x}\right|^2\right)^{\ff 1 2}\\
&<\infty.
\end{align*}
Then we have
$$\left\<\E \left[\ph(X_r^{t,x}) M_r^s(X_s^{t,x})\Big|\sF_s^t\right],\nn_v X_s^{t,x}\right\>=\E \left[\ph(X_r^{t,x})\left\< M_r^s(X_s^{t,x}),\nn_v X_s^{t,x}\right\>\Big|\sF_s^t\right].$$
This, together with \eqref{equ-1}, yields \eqref{lpcov}.

\textit{Step 2: Proof of assertion (1)}.
According to Step 1, it suffices to prove that
for $\ph\in C(\R^d)$ satisfying \eqref{gg}, $P_{s,r}\ph\in C^1(\R^d)$ and that it satisfies \eqref{Bism-P} and \eqref{po-nn-ph}.
To this end, we will invoke Lemma \ref{lem-1}.
Since $\ph\in C(\R^d)$, $\tld \ph_n$ defined as in Fact (ii) above belongs to $C_b(\R^d)$.
Then for each $n\geq 1$, there exists $\{\tld\ph_{n,m}\}_{m=1}^\infty\subset C_b^1(\R^d)$ such that
$\|\tld\ph_{n,m}\|_\infty\leq \|\tld\ph_n\|_\infty$ and for any $N>0$
\beg{align}\label{lim-M-nm}
\lim_{m\ra+\infty}\sup_{|x|\leq N}|\tld\ph_{n,m}(x)-\tld\ph_n(x)|=0.
\end{align}
Thus, a simple application of the H\"{o}lder inequality yields that for any $R>0$,
\beg{align*}
&\sup_{|y|\leq R}\E\left[ |(\tld\ph_{n,m}-\tld\ph_n)(X^{s,y}_r)|\left(|M_r^{s}(y)|+1\right)\right] \\
&\leq 2\|\tld\ph_n\|_\infty\sup_{|y|\leq R}\E\left[\left(|M_r^{s}(y)|+1\right) \1_{[ |X^{s,y}_r|\geq N]}\right]\\
&\quad+\left(\sup_{|x|\leq N}|\tld\ph_{n,m}(x)-\tld\ph_n(x)|\right)\sup_{|y|\leq R}\E \left(|M_r^{s}(y)|+1\right)\\
& \leq\ff 2 N\|\tld\ph_n\|_\infty\sup_{|y|\leq R}\left[\left(1+K(s,s,r,y)\right)\left(\E|X_r^{s,y}|^2\right)^{\ff 1 2}\right]\\
&\quad +\left(\sup_{|x|\leq N}|\tld\ph_{n,m}(x)-\tld\ph_n(x)|\right)\sup_{|y|\leq R}\left(1+K(s,s,r,y)\right),
\end{align*}
which, along with \eqref{lim-M-nm}, implies
$$\lim_{m\ra+\infty}\sup_{|y|\leq R}\E\left[|(\tld\ph_{n,m}-\tld\ph_n)(X^{s,y}_r)|\left(|M_r^{s}(y)|+1\right)\right]=0.$$
Then, combining this with Fact (ii), we derive that
there exists a sequence  $\{\ph_n\}_{n=1}^\infty\subset C_b^1(\R^d)$ such that \eqref{lim-sup} holds for $\ph$.
Thanks to \eqref{est-phM}, \eqref{gg} and \eqref{mon-est-1},
it is readily checked that \eqref{add-1(Lem2.2)} also holds for $\ph$.
Therefore, applying Lemma \ref{lem-1}, one sees that $P_{s,r}\ph$ belongs to $C^1(\R^d)$ and satisfies \eqref{Bism-P} and \eqref{po-nn-ph},
from which the assertion (1) follows.


\textit{Step 3: Proof of assertion (2)}.
Observe that assertion (1) holds for any $\ph\in C_b(\R^d)$ thanks to $ C_b(\R^d)\subset C(\R^d)$ and \eqref{gg} with $q=0$.
Then for any $\ph\in C_b(\R^d),\ep>0, v\in\R^d$ and $0\leq t\leq s<r\leq T$, we have
\beg{align}\label{ab-phF}
&\E\left[\ph(X_r^{t,x+\ep v})\Big|\sF_s^t\right] -  \E\left[\ph(X_r^{t,x})\Big|\sF_s^t\right] \nonumber\\
&=\int_0^{\ep} \nn_v\E\left[\ph(X_r^{t,x+\th v})\Big|\sF_s^t\right]\d \th\nonumber\\
&=\int_0^{\ep} \E \left[\left\<\ph(X_r^{t,x+\th v}) M_r^s(X_s^{t,x+\th v}),\nn_v X_s^{t,x+\th v}\right\>\Big|\sF_s^t\right]\d \th.
\end{align}
Let $\P_{X_r^{t,x+\ep v}}$ and $\P_{X_r^{t,x}}$ denote the law of  $X_r^{t,x+\ep v}$ and $X_r^{t,x}$, respectively,
and let
\beg{align*}
\mu_{\ep,r}^{t,s}(A)= \int_0^{\ep} \E\left( \1_A(X_r^{t,x+\th v})| M_r^s(X_s^{t,x+\th v})|\cdot|\nn_v X_s^{t,x+\th v}|\right)\d \th,\ \ A\in\sB(\R^d).
\end{align*}
From \eqref{in-K-pog} and \eqref{L2-con}, it is easy to see that $\mu_{\ep,r}^{t,s}$ is a finite measure on $\R^d$.
Then $C_b(\R^d)$ is dense in $L^1(\P_{X_r^{t,x+\ep v}}+\P_{X_r^{t,x}}+\mu_{\ep,r}^{t,s})$.
Hence, \eqref{ab-phF} holds for any $\ph\in\sB_b(\R^d)$.

Next, we intend to prove that for any $\ph\in\sB_b(\R^d)$,
$$\E \left[\left\<\ph(X_r^{t,x}) M_r^s(X_s^{t,x }),\nn_v X_s^{t,x }\right\>\Big|\sF_s^t\right],\ \ 0\leq t\leq s<r\leq T,$$
is continuous with respect to $x$ in $L^p(\P)$ for any $p\in [1,2)$.
Then it follows from \eqref{ab-phF} that \eqref{lpcov} holds in $L^p(\P)$ for $p\in[1,2)$ for any $\ph\in\sB_b(\R^d)$.\\
For any $\ph\in\sB_b(\R^d), x,y\in\R^d$ and $0<\de<r-s$, we obtain
\beg{align*}
&\qquad\E \left[\left\<\ph(X_r^{t,y}) M_r^s(X_s^{t,y }),\nn_v X_s^{t,y}\right\>\Big|\sF_s^t\right]-\E \left[\left\<\ph(X_r^{t,x}) M_r^s(X_s^{t,x }),\nn_v X_s^{t,x }\right\>\Big|\sF_s^t\right]\\
&\qquad =\E \left[\left\<\ph(X_r^{t,y}) M_r^s(X_s^{t,y }),\nn_v X_s^{t,y}-\nn_vX_s^{t,x}\right\>\Big|\sF_s^t\right] \\
&\qquad \quad +\E \left[\left\<\ph(X_r^{t,y})(M_r^s(X_s^{t,y })-M_r^s(X_s^{t,x })),\nn_v X_s^{t,x}\right\>\Big|\sF_s^t\right] \\
&\qquad\quad + \E \left[\left\<\left(\ph(X_r^{t,y})-\ph(X_r^{t,x}) \right) \left(M_{r }^s(X_s^{t,x })-M_{r-\de}^s(X_s^{t,x })\right),\nn_v X_s^{t,x}\right\>\Big|\sF_s^t\right]\\
&\qquad\quad +\E \left[\left\<\left(\ph(X_r^{t,y})-\ph(X_r^{t,x}) \right) M_{r-\de}^s(X_s^{t,x }),\nn_v X_s^{t,x}\right\>\Big|\sF_s^t\right]\\
&\qquad=: I_1+I_2+I_3+I_4.
\end{align*}
For any $p\in [1,2)$, by \eqref{equ-con}, \eqref{est-phM} (with $q=0$), \eqref{mon-est-1} and the H\"{o}lder inequality, we get
\beg{align}\label{in-nn-X0}
\E\left|I_1\right|^p &\leq \E\left|\E\left(|\ph(X_r^{t,y})|\cdot|M_r^s(X_s^{t,y})|\Big|\sF_s^t\right)|\nn_v X_s^{t,y}-\nn_vX_s^{t,x}|\right|^p\nonumber\\
&= \E\left|\left[\E\left(|\ph(X_r^{s,z})|\cdot|M_r^s(z)|\right)\right]\big|_{z=X_s^{t,y }}|\nn_v X_s^{t,y}-\nn_vX_s^{t,x}|\right|^p\nonumber\\
& \leq C\left(\E \left(1+|X_s^{t,y}|^{k_0}\right)^{\ff {2p} {2-p}}\right)^{\ff {2-p} {2}}\left(\E|\nn_v X_s^{t,y}-\nn_vX_s^{t,x}|^2\right) ^{\ff p 2}\nonumber\\
& \leq C\left(1+|y|^{k_0pk_2}\right)\left(\E|\nn_v X_s^{t,y}-\nn_vX_s^{t,x}|^2\right) ^{\ff p 2}.
\end{align}
Then \eqref{L2-con} yields that $\lim_{y\ra x}\E |I_1|^p=0$.\\
For the term $I_2$, from the H\"{o}lder inequality we have for any $p\in[1,2)$,
\beg{align*}
\E|I_2|^p\leq \E|\left\<\ph(X_r^{t,y})(M_r^s(X_s^{t,y })-M_r^s(X_s^{t,x })),\nn_v X_s^{t,x}\right\>|^p.
\end{align*}
By \eqref{L2-M}  of (H3) we arrive at
\beg{align*}
\P\text{-}\lim_{y\ra x}\left\<\ph(X_r^{t,y})(M_r^s(X_s^{t,y })-M_r^s(X_s^{t,x })),\nn_v X_s^{t,x}\right\>=0.
\end{align*}
Moreover, as in \eqref{in-nn-X0},  for any $p\in [1,2)$  and $R>0$, we have
\beg{align*}
\sup_{|y|\leq R}\E\left|\left\<\ph(X_r^{t,y}) M_r^s(X_s^{t,y }) ,\nn_v X_s^{t,x}\right\>\right|^p<\infty,
\end{align*}
which yields that
$$\left\{ | \<\ph(X_r^{t,y}) M_r^s(X_s^{t,y }) ,\nn_v X_s^{t,x} \> |^p\right\}_{|y|\leq R}$$
is uniformly integrable.
Consequently, we obtain that for $p\in [1,2)$ and $R>0$,
\beg{align*}
\left\{ | \<\ph(X_r^{t,y})(M_r^s(X_s^{t,y })-M_r^s(X_s^{t,x })),\nn_v X_s^{t,x} \> |^p\right\}_{|y|\leq R}
\end{align*}
is also uniformly integrable.
Then the dominated convergence theorem implies that $\lim_{y\ra x}\E |I_2|^p=0$.\\
For the term $I_3$, by the boundedness of $\ph$ and the H\"older inequality we get for $p\in [1,2)$,
\beg{align*}
\E \left| I_3\right|^p &\leq 2^p\|\ph\|_{\infty}^p\E\left|\E\left[ |M_{r }^s(X_s^{t,x })-M_{r-\de}^s(X_s^{t,x })|\Big| \sF_s^t\right] | \nn_vX_s^{t,x}|\right|^p\\
&\leq 2^p\|\ph\|_{\infty}^p \left(\E\left|\E\left[ |M_{r }^s(X_s^{t,x })-M_{r-\de}^s(X_s^{t,x })|\Big| \sF_s^t\right]\right|^{\ff {2p} {2-p}}\right)^{\ff {2-p} {2}}\left(\E|\nn_vX_s^{t,x}|^2\right)^{\ff {p} {2}}.
\end{align*}
Using \eqref{add-Hy} of (H3), we have
\beg{align*}
\P\text{-}\lim_{\de\ra 0^+}\E\left[ |M_{r }^s(X_s^{t,x })-M_{r-\de}^s(X_s^{t,x })|\Big| \sF_s^t\right]=0.
\end{align*}
By \eqref{equ-con} (with $\psi\equiv1$) and \eqref{in-K-pog}, we obtain that for any $\tilde{p}\geq1$,
\beg{align*}
&\E\left|\E\left[ |M_{r }^s(X_s^{t,x })-M_{r-\de}^s(X_s^{t,x })|\Big| \sF_s^t\right]\right|^{\tilde{p}}\\
& \leq  \E\left|\E\left[ |M_{r }^s(X_s^{t,x })|+|M_{r-\de}^s(X_s^{t,x })|\Big| \sF_s^t\right]\right|^{\tilde{p}}\\
&=\E\left[\left(\E|M_{r }^s(z)|+\E |M_{r-\de}^s(z)|\right)\Big|_{z=X_s^{t,x }} \right]^{\tilde{p}}\\
& \leq  \E\left[\left(K(s,s,r,z)+K(s,s,r-\de,z)\right)|_{z=X_s^{t,x }}\right]^{\tilde{p}}\\
& \leq  \left(C(s,s,r)+C(s,s,r-\de)\right)^{\tilde{p}}\E\left(1+|X_s^{t,x }|^{k_0}\right)^{\tilde{p}}.
\end{align*}
Combining this with \eqref{mon-est-1} and using the condition that $C(s,s, \cdot)$ is bounded on any closed subinterval of $(s,r]$,
we have
$$\sup_{\de <\de'}\E\left|\E\left[ |M_{r }^s(X_s^{t,x })-M_{r-\de}^s(X_s^{t,x })|\Big| \sF_s^t\right]\right|^{\tilde{p}}<\infty,
\ \ 0<\de'<r-s,\ \tilde{p}\geq1.$$
Then the dominated convergence theorem yields that $\lim_{\delta\ra0^+}\overline{\lim}_{y\ra x}\E |I_3|^p=0$.\\
For the term $I_4$, since $s<r-\de<r$ and $M_{r-\de}^s(X_s^{t,x })$ is measurable with respect to $F_{r-\de}^t$, we have
\beg{align*}
I_4 & =\E \left[\left\<\E\left[\left(\ph(X_r^{t,y})-\ph(X_r^{t,x})\right)\Big|\sF_{r-\de}^t\right] M_{r-\de}^s(X_s^{t,x }),\nn_v X_s^{t,x}\right\>\Big|\sF_s^t\right].
\end{align*}
Putting $\epsilon:=|x-y|$ and $v:=\ff {y-x} {|y-x|}$, by \eqref{ab-phF} we deduce that
\beg{align*}
&\qquad\E\left[\left(\ph(X_r^{t,y})-\ph(X_r^{t,x})\right)\Big|\sF_{r-\de}^t\right] \\
&\qquad = \int_0^{|x-y|}\E \left[\left\<\ph(X_r^{t,x+\th v}) M_r^{r-\de}(X_{r-\de}^{t,x+\th v}),\nn_v X_{r-\de}^{t,x+\th v}\right\>\Big|\sF_{r-\de}^t\right]\d \th.
\end{align*}
Then for any $p\in[1,2)$, it follows from \eqref{equ-con} (with $\psi\equiv1$), \eqref{in-K-pog} and \eqref{mon-est-1} that
\beg{align*}
&\qquad\E\left| \E\left[\left(\ph(X_r^{t,y})-\ph(X_r^{t,x})\right)\Big|\sF_{r-\de}^t\right] \right|^p\\
&\qquad\leq \|\ph\|_\infty|x-y|^{p-1}\int_0^{|x-y|}\E\left|\E\left[|M_r^{r-\de}(X_{r-\de}^{t,x+\th v})|\Big|\sF_{r-\de}^t\right]|\nn_v X_{r-\de}^{t,x+\th v}|\right|^p\d\th\\
&\qquad\leq C\|\ph\|_\infty|x-y|^{p-1}\int_0^{|x-y|}\E \left[\left(1+|X_{r-\de}^{t,x+\th v}|^{k_0}\right)^p |\nn_v X_{r-\de}^{t,x+\th v} |^p\right]\d\th\\
&\qquad\leq C\|\ph\|_\infty|x-y|^{p-1}\int_0^{|x-y|} \left(1+\|X_{r-\de}^{t,x+\th v}\|_{2k_0 p/(2-p)}^{\ff {2-p} 2} \right)  \|\nn_v X_{r-\de}^{t,x+\th v} \|_2^p\d\th\\
&\qquad\leq C\|\ph\|_\infty|x-y|^{p-1}\int_0^{|x-y|} \left(1+(1+|x+\th v|)^{\ff {(2-p)k_2} 2} \right)  \|\nn X_{r-\de}^{t,x+\th v} \|_2^p\d\th.
\end{align*}
Observe that by \eqref{L2-con}, $\|\nn X_{r-\de}^{t, x} \|_2^p$ is locally bounded with respect to $x$.
Consequently, we have
\beg{align*}
\lim_{y\ra x}\E\left| \E\left[\left(\ph(X_r^{t,y})-\ph(X_r^{t,x})\right)\Big|\sF_{r-\de}^t\right] \right|^p=0,
\end{align*}
which implies
\beg{align*}
\P\text{-}\lim_{y\ra x}\E\left[\left(\ph(X_r^{t,y})-\ph(X_r^{t,x})\right)\Big|\sF_{r-\de}^t\right]=0.
\end{align*}
Note that, again as in \eqref{in-nn-X0}, for any $p\in[1,2)$ we get,
\beg{align}\label{add-2(Prop2.1)}
\E\sup_{y\in\R^d} |I_4|^p & \leq2^p\|\ph\|_\infty^p \E\left(\E \left[ |M_{r-\de}^s(X_s^{t,x })|\Big|\sF_s^t\right]|\nn_v X_s^{t,x}|\right)^p\cr
& \leq C\left(1+|x|^{k_0pk_2}\right)\left(\E|\nn_vX_s^{t,x}|^2\right) ^{\ff p 2}
<\infty.
\end{align}
Then, by the dominated convergence theorem and \eqref{add-2(Prop2.1)} (with $p=1$) we obtain
\beg{align*}
\lim\limits_{y\ra x}\E|I_4| &\leq\lim\limits_{y\ra x}\E \left\{\left|\E\left[\left(\ph(X_r^{t,y})-\ph(X_r^{t,x})\right)\Big|\sF_{r-\de}^t\right]\right| \cdot|M_{r-\de}^s(X_s^{t,x })|\cdot|\nn_v X_s^{t,x}|\right\}=0,
\end{align*}
which leads to $\P\text{-}\lim_{y\ra x}|I_4|=0$.
Therefore, by the dominated convergence theorem and \eqref{add-2(Prop2.1)} again, we derive that $\lim_{y\ra x}\E|I_4|^p=0$ for any $p\in[1,2)$.

Let us now prove that for $\ph\in\sB(\R^d)$ satisfying \eqref{gg}, \eqref{lpcov} holds in $L^p(\P)$ for $p\in[1,2)$.
Notice that for such $\phi$, \eqref{add-1(Lem2.2)} holds true because of \eqref{gg} and  \eqref{est-phM},
and due to Fact (ii) there exists a sequence $\{\tld \ph_n\}_{n=1}^\infty\subset\sB_b(\R^d)$ such that \eqref{lim-sup} holds.
Since we have proved that \eqref{lpcov} holds for each $\tld \ph_n$,
by letting $s=t$ we get
\beg{align*}
(\nn_v P_{t,r}\tld \ph_n(x)=)\nn_v\E\left[\tld \ph_n(X_r^{t,x})\right]=\E \left[ \tld \ph_n(X^{t,x}_r)\<M_r^t(x), v\> \right],\ \ v\in\R^d.
\end{align*}
From the result proved in the above paragraphs,
we know that $\E[ \tld \ph_n(X^{t,x}_r)\<M_r^t(x), v\>]$ is continuous with respect to $x$,
that is, $\nn P_{s,r}\tld \ph_n\in C(\R^d,\R^d)$.
Then along the same lines as in \eqref{lim-phn} and \eqref{in-po-nnP},
we derive that $P_{s,r}\ph\in C^1(\R^d)$ and it satisfies \eqref{Bism-P} and \eqref{po-nn-ph},
which completes the proof of assertion (2), thanks to Step 1.
\qed


\subsection{Main results for FBSDEs}

We now consider the following decoupled FBSDE for $(t,x)\in[0,T]\times\R^d$:
\beg{align}\label{equ-X}
X^{t,x}_s & =x+\int_t^s b(r,X^{t,x}_r)\d r+\int_t^s\si(r,X^{t,x}_r)\d W_r,\ \ s\in[t,T],\\
Y_s^{t,x} & =g(X_T^{t,x})+\int_s^T f(r,X_r^{t,x},Y_r^{t,x},Z_r^{t,x})\d r-\int_s^T Z_r^{t,x}\d W_r,\ \ s\in[t,T].\label{BSDE1}
\end{align}
Here $b:[0,T]\times\R^d\ra\R^d,\si:[0,T]\times\R^d\ra\R^d\otimes\R^m$
and $g:\R^d\ra\R^l,f:[0,T]\times\R^d\times\R^l\times\R^l\otimes\R^m\ra\R^l$
are measurable mappings, $W$ is an $m$-dimensional Brownian motion.
Instead of  concrete conditions imposed on $b$ and $\si$, we assume that $X^{t,x}_s$ has the following properties:
\beg{enumerate}
\item [(C1)]  The system consisting of \eqref{equ-X} and the equation
\beg{align}\label{Jaco}
\nn  X_s^{t,x}=I_{d\times d}+\int_t^s \nn b(r,X_r^{t,x})\nn X_r^{t,x}\d r+\int_t^s\nn\si(r,X_r^{t,x})\nn X_r^{t,x}\d W_r,\ \ s\in[t,T]
\end{align}
has a unique strong solution, where $I_{d\times d}$ is the identity matrix in $\R^d\otimes \R^d$,
$\nn b$ and $\nn\si$ are the weak derivatives of $b$ and $\si$, respectively.
There exist $k_2\geq 1$, $k_3\geq 0$ such that for any $k_1\geq 1$, \eqref{mon-est-1}  holds for some $C_{k_1,k_2,T}>0$ and
\beg{align}\label{mon-est-2}
\sup_{0\leq s\leq r\leq T}\|\nn X_r^{s,y}\|_{k_1}&\leq C_{k_1,k_3,T}(1+|y|^{k_3}),\ \ y\in\R^d,
\end{align}
for some  $C_{k_1,k_3,T}>0$. Moreover, for any $t\leq s\leq T$ and $x,v\in\R^d$,
\beg{align}\label{add-der}
\lim_{\ep\ra 0^+}\E\sup_{s\in [t,T]}&\left|\ff {X_s^{t,x+\ep v}-X_s^{t,x}} {\ep}-\nn_v X_s^{t,x}\right|^2=0, \\
&\P\text{-}\lim_{y\ra x}\nn_v X_s^{t,y}=\nn_v X_s^{t,x}. \label{con-nn-P}
\end{align}

\item [(C2)]  There exist continuously differentiable $(b_n,\si_n)_{n\geq 1}$ with
\beg{align}\label{Revise-C1-b}
\sup_{r\in[0,T],\ x\in\R^d}\left( \| \nn b_n(r,x)\|+ \| \nn \si_n(r,x)\|\right)<\infty,\ \ \forall n\in\mathds{N},
\end{align}
such that for any sequence $x_n$ with $\lim_{n\ra +\infty}x_n=x$ and for the solution of  \eqref{equ-X} and \eqref{Jaco} with $(x,b,\si)$ replaced by $(x_n,b_n,\si_n)$, which we denote by $(X_s^{n,t,x_n},\nn X_s^{n,t,x_n})$, the following hold
\beg{align}\label{app-X1}
&\lim_{n\ra +\infty}\E \left(|X_T^{n,t,x_n}-X_T^{t,x}|^2+\int_t^T \left|X_r^{n,t,x_n}-X_r^{t,x}\right|^2\d r\right)=0,\\    \label{app-X2}
&\lim_{n\ra +\infty}\E \left(\|\nn X_T^{n,t,x_n}-\nn X_T^{t,x}\|^2+\int_t^T \left\|\nn X_r^{n,t,x_n}-\nn X_r^{t,x}\right\|^2\d r\right)=0,\\
&\lim_{n\ra+\infty} \|\si_n(s,x_n)-\si(s,x)\|=0, ~s\in [t,T].\label{app-si}
\end{align}
\end{enumerate}

\beg{rem}\label{Rem-C12}
It is clear that if $b(r,\cdot)$ and $\si(r,\cdot)$ are continuously differentiable and satisfy \eqref{C1-b} below, then (C1) and (C2) hold.
Besides, we may allow that $\si$ has polynomial growth which also ensures that (C1) and (C2) are satisfied (see Example \ref{ex.Gru} below).
\end{rem}

Now, we introduce the hypotheses on the coefficients $g$ and $f$ of the BSDE \eqref{BSDE1},
 under which we will be able to establish gradient type and Bismut type formulas for the solution $Y_s^{t,x}$.
\beg{enumerate}
\item[(A1)] $g$ is Lipschitz continuous, i.e. there exists a constant $K_1\geq 0$ such that
\beg{align*}
|g(x)-g(y)|\leq K_1|x-y|,\ \ x,y\in\R^d.
\end{align*}
\item[(A1$'$)] $g$ has $q$-th growth, i.e. there exist constants $\tld K_1\geq0$ and $q\geq0$ such that
\beg{align*}
|g(x)|\leq \tld K_1(1+|x|^q),\ \ x\in\R^d.
\end{align*}
\item[(A2)] There exist constants $K_2\geq0$ and $\tld K_2\geq0$ such that for all $(x_i,y_i,z_i)\in \R^{d}\times\R^{l}\times\R^{l}\otimes\R^m$,
\beg{align*}
\sup_{r\in [0,T]}|f(r,x_1,y_1,z_1)-f(r,x_2,y_2,z_2)|&\leq K_2(|x_1-x_2|+|y_1-y_2|+\|z_1-z_2\|),\\
\tld K_2:=\sup_{r\in [0,T]}&|f(r,0,0,0)|<\infty.
\end{align*}
\end{enumerate}

To simplify the notations, for any $t\leq s<r\leq T$ and $x,v\in R^d$, we set $\Th_r^{t,x}:=(X_r^{t,x},Y_r^{t,x},Z_r^{t,x})$
and $M_r^{t,s}(x,v):=\<M_r^s(X_s^{t,x}),\nn_v X_s^{t,x}\>$.
In particular, $M_r^{t,t}(x,v)=\<M_r^t(x),v\>$.
Now we formulate our main result, the proof of which we present in the next subsection.
\beg{thm}\label{thm-BSDE}
Let $0\leq t\leq s<T$ and $x_0,v\in\R^d$.
Assume that (C1), (C2), (H2) and (A2) are satisfied and that \eqref{equ-con} holds for any $\ps\in \sB_b(\R^d)$ $\P$-a.s. \\
(1) Let (A1) hold, and assume furthermore that  $\sigma(r,\cdot)$ is continuous and there exist $C>0$ and $k_4\geq 0$ such that
\beg{align}\label{Add-1(The2.3)}
\sup_{r\in[0,T]}\|\si(r,x)\|\leq C(1+|x|^{k_4}),\ \ x\in\R^d,
\end{align}
and that $K(t,s,r,\cdot)$ satisfies \eqref{in-K-pog} and for some $R>0$
\beg{align}\label{iint-K}
\int_s^T\sup_{|x-x_0|\leq R} K(t,s,r,x)\d r<\infty.
\end{align}
If in addition, either $g(\cdot)$ and $f(r,\cdot,\cdot,\cdot)$ are continuously differentiable, or \eqref{L2-M} holds,
then the following gradient type formula in $L^p(\P)$ holds for any $p\in[1,2)$:
\beg{align}\label{equ-deri}
\nn_v Y^{t,x_0}_s= \E\left[ g(X_T^{t,x_0})M_T^{t,s}(x_0,v)+ \int_s^T f(r,\Th_r^{t,x_0})M_r^{t,s}(x_0,v)\d r\Big| \sF_s^t \right].
\end{align}
\\
(2) Let (A1$'$) hold,
and assume furthermore that $k_2=1$,
that $\si$ is bounded on $[0,T]\times \R^d$,
and that $K(t,s,r,\cdot)$ satisfies \eqref{in-K-pog} with $k_0=0$ and there exists  $\be<1$ such that
\beg{align}\label{CC}
\sup_{0\leq t\leq s<r\leq T}\left[C(t,s,r)(r-s)^{\be}\right]<\infty.
\end{align}
If in addition, either $g\in C(\R^d;\R^l)$ and \eqref{L2-M} holds, or (H3) is satisfied,
then \eqref{equ-deri} also holds in $L^p(\P)$ for any $p\in[1,2)$.
\end{thm}

Concerning a Bismut type formula of $Y^{t,x}_s$, we have the following corollary.
The proof of this corollary follows from the same kind of arguments as Theorem \ref{thm-BSDE} with the conditional expectation $\E[\cdot|\sF_s^t] $ replaced by the expectation $\E[\cdot]$, and it is thus omitted here for the sake of conciseness.
\beg{cor}\label{cor2}
Let $0\leq t\leq s<T$ and $x_0, v\in\R^d$.
Assume that (C1), (C2), (H2) and (A2) are satisfied. \\
(1) Let (A1) hold, and assume furthermore that  $\sigma(r,\cdot)$ is continuous and there exist $C>0$ and $k_4\geq 0$ such that
\beg{align*}
\sup_{r\in[0,T]}\|\si(r,x)\|\leq C(1+|x|^{k_4}),\ \ x\in\R^d,
\end{align*}
and that
$K(t,t,r,\cdot)$ satisfies \eqref{in-K-pog} and for some $R>0$
\beg{align}\label{iint-K'}
\int_s^T\sup_{|x-x_0|\leq R} K(t,t,r,x)\d r<\infty.
\end{align}
If in addition, either $g(\cdot)$ and $f(r,\cdot,\cdot,\cdot)$ are continuously differentiable, or \eqref{L2-M} (with $s=t$) holds,
then the following Bismut type formula in $L^p(\P)$ holds for any $p\in[1,2)$:
\beg{align}\label{equ-Bis}
\nn_v\E Y^{t,x_0}_s= \E\left[ g(X_T^{t,x_0})\<M_T^t(x_0),v\>+ \int_s^T f(r,\Th_r^{t,x_0})\<M_r^t(x_0),v\>\d r\right].
\end{align}
\\
(2) Let (A1$'$) hold,
and assume furthermore that $k_2=1$,
that $\si$ is bounded on $[0,T]\times \R^d$,
and that $K(t,t,r,\cdot)$ satisfies \eqref{in-K-pog} with $k_0=0$ and for some $\be<1$
\beg{align*}
\sup_{0\leq t<r\leq T}\left[C(t,t,r)(r-t)^{\be}\right]<\infty.
\end{align*}
If in addition, either $g\in C(\R^d;\R^l)$ and \eqref{L2-M} (with $s=t$) holds, or (H3) (with $s=t$) is satisfied,
then \eqref{equ-Bis} also holds in $L^p(\P)$ for any $p\in[1,2)$.
\end{cor}


\subsection{Proof of the main result}

To prove Theorem \ref{thm-BSDE}, we need the following useful lemma.
Its proof is elementary but lengthy, and we defer it to the Appendix.
\beg{lem}\label{rep-YZ-X}
Assume that (C1), (C2), (A1) and (A2) hold,
and that $g(\cdot)$ and $f(r,\cdot,\cdot,\cdot)$ are continuously differentiable with any $r\in [0,T]$.
For each $(t,x)\in[0,T]\times\R^d$,
let $(X_\cdot^{t,x},Y_\cdot^{t,x},Z_\cdot^{t,x})$ be the solution of \eqref{equ-X}-\eqref{BSDE1}
and set $u(t,x):= Y_t^{t,x}$.
Then for any $t\in[0,T]$, $u(t,\cdot)$ is continuously differentiable with
\beg{align}\label{equ-nu0}
\nn u(t,x)=\E\left[ \nn g(X_T^{t,x})\nn X_T^{t,x}+\int_t^T \nn f(r,\Th_r^{t,x})\nn \Th_r^{t,x} \d r\right],\ \  x\in\R^d,
\end{align}
and there exist positive constants $C_{K_1,K_2,k_3,T}$ and $q_1$ such that
\beg{align}\label{ine-nu0}
\sup_{t\in [0,T]}\|\nn u(t,x)\|\leq C_{K_1,K_2,k_3,T}(1+|x|^{  q_1}),\ \ x\in\R^d.
\end{align}
Moreover, $Y_s^{t,x}=u(s,X_s^{t,x})$ and $Z_s^{t,x}=\nn u(s,X_s^{t,x})\si(s,X_s^{t,x}),\ \ \d s\otimes\d\P$-a.e.
\end{lem}

\noindent{\bf{\emph{Proof of Theorem \ref{thm-BSDE}}}}

First observe that, by \eqref{BSDE1} and the Fubini theorem we have
\beg{align}\label{Rep-Y-0}
Y_s^{t,x}&= \E \left[ g(X_T^{t,x})+ \int_s^T f(r,\Th_r^{t,x})\d r\Big|\sF_s^t\right]\cr
&= \E \left[ g(X_T^{t,x})\big|\sF_s^t\right]+ \int_s^T \E [f(r,\Th_r^{t,x})\big|\sF_s^t]\d r.
\end{align}
Below we shall apply Proposition \ref{lem-lpc} to $f$ and $g$ in two different cases (1) and (2), respectively.
Without loss of generality, we assume that $l=1$.
Otherwise, we may use $\<Y_s^{t,x},e\>,\<g,e\>$ and $\<f,e\>$ to replace $Y_s^{t,x},g$ and $f$ for a unit vector $e\in\R^l$.

We first prove assertion (1).
If (A1) holds, then $g\in C(\R^d)$ and it satisfies \eqref{gg} with $q=1$.
So, with the help of Proposition \ref{lem-lpc} (1) we have that
\beg{align}\label{equ-g-de}
\nn_v\E \left[ g(X_T^{t,x_0}) \Big| \sF_s^t\right]= \E\left[ g(X_T^{t,x_0})M_T^{t,s}(x_0,v)\Big|\sF_s^t\right]
\end{align}
in $L^p(\P)$ for any $p\in[1,2)$.

Now, we focus on dealing with the directional derivative of the second term of the right-hand side of \eqref{Rep-Y-0}.
For the coefficients $g(\cdot),f(r,\cdot,\cdot,\cdot)$ which fulfill (A1) and (A2) respectively,
by standard approximation arguments there exist continuously differentiable sequences $(g_n(\cdot))_{n\geq 1},(f_n(r,\cdot,\cdot,\cdot))_{n\geq 1}$,
in which for each $n\geq 1,g_n$ and $f_n$ satisfy (A1) and (A2) with the same Lipschitz constants $K_1,K_2$
and some positive constant $\bar{K}_2$ (
independent of $n$),
and moreover $g_n$ converges to $g$  uniformly on $\R^d$ and $f_n$ converges to $f$ uniformly on $[t,T]\times\R^d\times\R\times\R^m$.
We now let $(Y^{n,t,x},Z^{n,t,x})$ be the solution of \eqref{BSDE1} with coefficients $g,f$ replaced by $g_n,f_n$, respectively.
By the It\^{o} formula, one can show that
\beg{align}\label{con-Y-Zn}
&\lim_{n\ra+\infty}\left(\sup_{r\in [t,T]}\E\left|Y_r^{n,t,x}-Y_r^{t,x}\right|^2+\int_t^T\E\left|Z_r^{n,t,x}-Z_r^{t,x}\right|^2\d r\right)=0.
\end{align}
The above step is partially borrowed from \cite[Theorem 4.2, Page 1410-1411]{MaZ}.
By Lemma \ref{rep-YZ-X}, there exist $u_n:[0,T]\times\R^d\ra\R$ continuously differentiable in the space variable,
positive constants $C_{K_1,K_2,\bar{K}_2,T}$ and $q_1$ independent of $n$ such that
\beg{align}\label{Add-2(The2.3)}
\sup_{t\in[0,T]}\left(|u_n(t,x)|+|\nn u_n(t,x)|\right)\leq C_{K_1,K_2,\bar{K}_2,T}(1+|x|^{q_1}),\ \ x\in\R^d
\end{align}
and
$$\qquad\qquad Y_s^{n,t,x}=u_n(s,X^{t,x}_s),\ Z_s^{n,t,x}=\nn u_n(s,X_s^{t,x})\si(s,X^{t,x}_s),\ \ \d s\otimes\d\P\text{-a.e.}$$
Let $\Th_\cdot^{n,t,x}:=(X_\cdot^{ t,x },Y_\cdot^{n,t,x },Z_\cdot^{n,t,x })$.
Then, we have
$$ f(r,\Th_r^{n,t,x})=f(r,X_r^{t,x},u_n(r,X^{t,x}_r),\nn u_n(r,X_r^{t,x})\si(r,X^{t,x}_r)).$$
For any $r\in[t,T]$, let
$F_n(r,\cdot):=f(r,\cdot,u_n(r,\cdot),\nn u_n(r,\cdot)\si(r,\cdot))$.
Then, by (A2), the continuity of $\sigma(r,x)$ in $x$, \eqref{Add-1(The2.3)} and \eqref{Add-2(The2.3)},
we see that for a.e. $r\in[t,T]$, $F_n(r,\cdot)\in C(\R^d)$ and
\beg{align}\label{ineq-f-u}
|F_n(r,x)|&\leq \tld K_2+K_2\left(|x|+|u_n(r,x)|+ |\nn u_n(r,x)\si(r,x)|\right)\nonumber\\
& \leq  \tld K_2 +K_2\Big[|x|+C_{K_1,K_2,\bar{K}_2,T}(1+|x|^{q_1}) + \nonumber\\
&\qquad\qquad\ \ \ \ \  + C_{K_1,K_2,\bar{K}_2,T}C(1+|x|^{q_1})(1+|x|^{k_4})\Big]\nonumber\\
&\leq C_{K_1,K_2,\tld K_2,\bar{K}_2,T} \left( 1+|x|^{(q_1+k_4)\vee1}\right).
\end{align}
Consequently, using (C1), (H2) and \eqref{in-K-pog}-\eqref{equ-con}, we then apply Proposition \ref{lem-lpc} (1) to conclude that for a.e. $r\in[t,T]$,
\beg{align}\label{nn-v-Ef1}
\nn_v \E\left[ f(r,\Th_r^{n,t,x})  \Big| \sF_s^t \right]&=\nn_v \E\left[ F_n(r,X_r^{t,x})  \Big| \sF_s^t \right]\cr
&= \E\left[F_n(r,X_r^{t,x}) M_r^{t,s}(x,v) \Big| \sF_s^t \right]\nonumber \\
&= \E\left[ f(r,\Th_r^{n,t,x})  M_r^{t,s}(x,v) \Big| \sF_s^t \right],\ \ v\in\R^d
\end{align}
in $L^p(\P)$ for any $p\in [1,2)$.
In view of \eqref{ineq-f-u}, \eqref{mon-est-1}, \eqref{mon-est-2} and the fact that $C_{K_1,K_2,\tld K_2,\bar{K}_2,T}$ and $q_1$ above are independent of $n$, we obtain that for a.e. $r\in[t,T]$ and $p\in [1,2)$,
\beg{align}\label{supfM}
&\left\|\E\left[ \sup_{n\geq 1}|f(r,\Th_r^{n,t,x})  M_r^{t,s}(x,v) |\Big| \sF_s^t \right]\right\|_p \nonumber\\
& \leq C_{K_1,K_2,\tld K_2,\bar{K}_2,T} \left\|  (1+|X_r^{t,x}|^{(q_1+k_4)\vee1})|M_r^{s}(X_s^{t,x})|\cdot |\nn_vX_s^{t,x}| \right\|_p \nonumber\\
& \leq C_{K_1,K_2,\tld K_2,\bar{K}_2,T} \left\|  (1+|X_r^{t,x}|^{(q_1+k_4)\vee1})|\nn_vX_s^{t,x}|\right\|_{\ff {2p} {2-p}}\left\|M_r^{s}(X_s^{t,x})  \right\|_2 \nonumber\\
& \leq C_{K_1,K_2,\tld K_2,\bar{K}_2,k_2,k_3,k_4,q_1,p,T}\left(1+|x|^{((q_1+k_4)\vee1)k_2}\right)(1+|x|^{k_3})K(t,s,r,x)|v|.
\end{align}
This, together with \eqref{iint-K}, leads to
\beg{align}\label{insupfM}
\int_s^T\sup_{|x-x_0|\leq R}\left\|\E\left[ \sup_{n\geq 1}|f(r,\Th_r^{n,t,x})  M_r^{t,s}(x,v) |\Big| \sF_s^t \right]\right\|_p\d r<\infty.
\end{align}
Then, using \eqref{nn-v-Ef1}, \eqref{con-Y-Zn} and the dominated convergence theorem, we deduce that
the following relations hold in $L^p(\P)$ for $p\in [1,2)$ and {$x\in B_R(x_0):=\{z: |z-x_0|<R\}$:
\beg{align}\label{nvf}
\nn_v \int_s^T\E \left[ f(r,\Th_r^{n,t,x})\Big|\sF_s^t\right]\d r & = \int_s^T\nn_v\E \left[ f(r,\Th_r^{n,t,x})\Big|\sF_s^t\right]\d r\nonumber\\
& = \int_s^T\E\left[ f(r,\Th_r^{n,t,x})M_r^{t,s}(x,v) \Big| \sF_s^t \right]\d r
\end{align}
and
\beg{align}\label{lim-fM}
\lim_{n\ra+\infty}\int_s^T\E\left[ f(r,\Th_r^{n,t,x})M_r^{t,s}(x,v) \Big| \sF_s^t \right]\d r
 =\int_s^T\E\left[ f(r,\Th_r^{t,x})M_r^{t,s}(x,v) \Big| \sF_s^t \right]\d r.
\end{align}
Hence, \eqref{nvf}-\eqref{lim-fM} and the dominated convergence theorem imply that for any $\epsilon>0$ and $v\in\R^d$ with $x+\ep v\in B_R(x_0)$,
\beg{align}\label{LipifM}
&\int_s^T\E \left[ f(r,\Th_r^{t,x+\epsilon v})\Big|\sF_s^t\right]\d r-\int_s^T\E \left[ f(r,\Th_r^{ t,x})\Big|\sF_s^t\right]\d r\nonumber\\
& =\lim_{n\ra+\infty}\left(\int_s^T\E \left[ f(r,\Th_r^{n,t,x+\epsilon v})\Big|\sF_s^t\right]\d r-\int_s^T\E \left[ f(r,\Th_r^{n,t,x})\Big|\sF_s^t\right]\d r\right)\nonumber\\
&  = \lim_{n\ra+\infty}\int_0^\epsilon \left(\int_s^T\E\left[ f(r,\Th_r^{n,t,x+\th v})M_r^{t,s}(x+\th v,v) \Big| \sF_s^t \right]\d r\right) \d \th\nonumber\\
& = \int_0^\epsilon \left(\int_s^T\E\left[ f(r,\Th_r^{ t,x+\th v})M_r^{t,s}(x+\th v,v) \Big| \sF_s^t \right]\d r\right) \d \th,
\end{align}
where the first relation is due to (A2), \eqref{con-Y-Zn} and \eqref{ineq-f-u}.

Next, we will show that
$$\int_s^T\E\left[ f(r,\Th_r^{ t,x })M_r^{t,s}(x ,v) \Big| \sF_s^t \right]\d r$$
is continuous with respect to $x$ in a neighbourhood of $x_0$ in $L^p(\P)$ for $p\in[1,2)$.
Observe that this continuity, along with \eqref{LipifM}, easily yields that
\beg{align}\label{Add-5(The2.3)}
\nn_v \int_s^T\E \left[ f(r,\Th_r^{t,x_0})\Big|\sF_s^t\right]\d r=\int_s^T\E\left[ f(r,\Th_r^{t,x_0})M_r^{t,s}(x_0,v) \Big| \sF_s^t \right]\d r
\end{align}
in $L^p(\P)$ for $p\in[1,2)$,
which due to our previous equality \eqref{equ-g-de}, leads to \eqref{equ-deri}.
We now invoke the condition \eqref{L2-M} to get the continuity.
Indeed, due to (A1), (A2) and \eqref{add-der}, we can apply the It\^{o} formula and the B-D-G inequality to get
\beg{align*}
\lim_{y\ra x}\E \left[\sup_{r\in [t,T]}| Y_r^{t,y}-Y_r^{t,x}|^2 +\int_t^T \left|Z_r^{t,y}-Z_r^{t,x}\right|^2\d r\right]=0.
\end{align*}
Combining this with \eqref{L2-M} and \eqref{con-nn-P} yields
\beg{align}\label{Add-7(The2.3)}
\d r\otimes\d \P\text{-}\lim_{y\ra x}f(r,\Th_r^{ t, y})M_r^{t,s}(y ,v)=f(r,\Th_r^{ t,x })M_r^{t,s}(x ,v).
\end{align}
Observe that by \eqref{ineq-f-u}, we have
\beg{align*}
|f(r,\Th_r^{n,t,x})|\leq C_{K_1,K_2,\tld K_2,\bar{K}_2,T} \left( 1+|X_r^{t,x}|^{(q_1+k_4)\vee 1}\right),
\end{align*}
which, along with \eqref{con-Y-Zn}, implies
\beg{align*}
|f(r,\Th_r^{t,x})|\leq C_{K_1,K_2,\tld K_2,\bar{K}_2,T} \left( 1+|X_r^{t,x}|^{(q_1+k_4)\vee 1}\right).
\end{align*}
Then, as in \eqref{supfM} and \eqref{insupfM}, we derive that for any $p\in[1,2)$,
\beg{align}\label{Add-6(The2.3)}
\left\{\left|f(r,\Th_r^{ t,x })M_r^{t,s}(x ,v)\right|^p\right\}_{|x-x_0|\leq R}
\end{align}
is uniformly integrable and
\beg{align}\label{isupfM}
\int_s^T\sup_{|x-x_0|\leq R}\left\|\E\left[ |f(r,\Th_r^{ t,x})  M_r^{t,s}(x,v) |\Big| \sF_s^t \right]\right\|_p\d r<\infty.
\end{align}
Thus, by \eqref{Add-7(The2.3)} and \eqref{Add-6(The2.3)} the dominated convergence theorem implies that
$$\lim_{y\ra x}\E\left[ f(r,\Th_r^{ t, y})M_r^{t,s}(y ,v) \Big| \sF_s^t \right]=\E\left[ f(r,\Th_r^{ t,x })M_r^{t,s}(x ,v) \Big| \sF_s^t \right], \ \ x\in B_R(x_0),$$
in $L^p(\P)$ for any $p\in[1,2)$ and a.e. $r\in [t,T]$.
Consequently, again by the dominated convergence theorem and \eqref{isupfM} we get the desired continuity result.

Note that if $g(\cdot)$ and $f(r,\cdot,\cdot,\cdot)$ are continuously differentiable,
then we just set $g_n\equiv g$ and $f_n\equiv f$, and the assertion \eqref{equ-deri} follows from \eqref{nvf} and \eqref{equ-g-de}.

We now prove assertion (2).
Let $g\in C(\R^d)$ satisfying (A1$'$), then we apply directly Proposition \ref{lem-lpc} (1) to get \eqref{equ-g-de}.
While in the case of $g$ only satisfying (A1$'$),
we first observe that \eqref{L2-con} holds due to \eqref{mon-est-2} and \eqref{con-nn-P},
and then by (H3) and \eqref{CC} it is easily checked in our context that the hypotheses of Proposition \ref{lem-lpc} (2) are all satisfied,
which also yields \eqref{equ-g-de}.

Next, we investigate the directional derivative of the second term in the right-hand side of \eqref{Rep-Y-0} by adopting a similar approximation argument as for assertion (1).
For $f$ satisfying (A2), let $(f_n)_{n\geq 1}$ be as defined above.
For $g$ satisfying (A1$'$), we set $\tld g_n:=(g\vee (-n))\we n$ for each $n\geq 1$.
Then it follows from the Lusin theorem (see, e.g., \cite[Theorem 7.4.4]{Coh}), that there exist $\{\hat g_n\}_{n\geq 1}\subset C_b(\R^d)$ and compact sets $\{\sK_n\}_{n\geq 1}$ such that
\beg{align}\label{h-t-g}
\hat g_n\Big|_{\sK_n}=\tld g_n\Big|_{\sK_n},\qquad \|\hat g_n\|_\infty\leq n,\qquad\nu(\sK_n^c)\leq \ff 1 {n^3},
\end{align}
where $\nu=\P_{X_T^{t,x}}$ is the law of $X_T^{t,x}$.
By a standard approximation method, we may choose
$g_n\in C_b^1(\R^d)$ such that
\beg{align}\label{gn-hgn}
\sup_{x\in\sK_n}|g_n(x)-\hat g_n(x)|\leq \ff 1 n,\qquad \|g_n\|_\infty\leq n.
\end{align}

Now, let $(Y^{n,t,x}, Z^{n,t,x})$ be the solution of \eqref{BSDE1} with coefficients $(g_n, f_n)$ replacing $(g, f)$.
Applying the It\^{o} formula and the B-D-G inequality, we obtain that
\beg{align}\label{Add-3(The2.3)}
&\E\left(\sup_{r\in [t,T]}|Y_r^{n,t,x}-Y_r^{t,x}|^2+\int_t^T|Z_r^{n,t,x}-Z_r^{t,x}|^2\d r\right)\cr
&\leq C\left(\E| g_n(X_T^{t,x})-g(X_T^{t,x})|^2+\|f_n-f\|_\infty\right)
\end{align}
and
\beg{align}\label{Add-4(The2.3)}
\E\left(\sup_{r\in [t,T]}|Y_r^{n,t,x}|^2+\int_t^T|Z_r^{n,t,x}|^2\d r\right)
\leq C\left(\E| g_n(X_T^{t,x})|^2+1\right)
\end{align}
with some constant $C>0$ independent of $n$ and $t$.
Note that from \eqref{h-t-g} and \eqref{gn-hgn}, we have
\beg{align*}
|(g_n-g)(x)|&\leq |(g_n-\hat g_n)(x)|+|(\hat g_n-\tld g_n)(x)|+|(\tld g_n-g)(x)| \\
&=|(g_n-\hat g_n)(x)|\1_{\sK_n}(x)+|(g_n-\hat g_n)(x)|\1_{\sK_n^c}(x)\\
&\quad+|(\hat g_n-\tld g_n)(x)|\1_{\sK_n^c}(x)+|(\tld g_n-g)(x)|\1_{[|g(x)|\geq n]}\\
&\leq  \ff 1 {n}+4n\1_{\sK_n^c}(x) +|g(x)|\1_{[|g(x)|\geq n]},
\end{align*}
which yields
\beg{align*}
\E\left| g_n(X_T^{t,x})-g(X_T^{t,x})\right|^2&\leq \ff 3 {n^2}+48n^2\P\left(X_T^{t,x}\not\in\sK_n\right)+
3\E\left[|g(X_T^{t,x})|^2 \1_{[|g(X_T^{t,x})|\geq n]}\right].
\end{align*}
Then, again by \eqref{h-t-g}, (A1$'$) and \eqref{mon-est-1} we get
\beg{align*}
\lim_{n\ra +\infty}\E| g_n(X_T^{t,x})-g(X_T^{t,x})|^2=0.
\end{align*}
Consequently, substituting this into \eqref{Add-3(The2.3)} and using the uniform convergence of $f_n$ to $f$ we obtain
\beg{align}\label{cov-yzn}
\lim_{n\ra+\infty}\E\left(\sup_{r\in [t,T]}|Y_r^{n,t,x}-Y_r^{t,x}|^2+\int_t^T|Z_r^{n,t,x}-Z_r^{t,x}|^2\d r\right)=0.
\end{align}
We now intend to find an upper bound for \eqref{Add-4(The2.3)} independent of $n$ and $t$.
By \eqref{h-t-g}, \eqref{gn-hgn}, (A1$'$) and \eqref{mon-est-1} with $k_2=1$, we have
\beg{align}\label{Add-8(The2.3)}
\| g_n(X_T^{t,x})\|_2&\leq\|\hat g_n(X_T^{t,x})\1_{\sK_n}(X_T^{t,x})\|_2+\|\hat g_n(X_T^{t,x})\1_{\sK_n^c }(X_T^{t,x})\|_2
+\|g_n(X_T^{t,x})-\hat g_n(X_T^{t,x})\|_2\nonumber\\
& \leq \|g(X_T^{t,x})\|_2+\|\hat g_n(X_T^{t,x})\1_{\sK_n^c }(X_T^{t,x})\|_2\nonumber\\
&\quad+\ff 1 n+\|(g_n(X_T^{t,x})-\hat g_n(X_T^{t,x}))\1_{\sK_n^c }(X_T^{t,x})\|_2\nonumber\\
&\leq C(1+|x|^q)+n\sq{\nu(\sK_n^c)}+\ff 1 n+2n\sq{\nu(\sK_n^c)}\nonumber\\
&\leq C(1+|x|^q),
\end{align}
where $C$ is independent of $n$ and $t$.
Then, plugging this into \eqref{Add-4(The2.3)} yields
\beg{align}\label{sup-un1}
\sup_{n\geq 1,t\in [0,T] } \E\left(\sup_{r\in [t,T]}|Y_r^{n,t,x}|^2+\int_t^T|Z_r^{n,t,x}|^2\d r\right)\leq C(1+|x|^q)^2.
\end{align}

On the other hand, since $g_n,f_n$ are continuously differentiable, we can apply Lemma \ref{rep-YZ-X} to derive that
there exists $u_n:[0,T]\times\R^d\ra\R$ such that $\nn u_n(r,\cdot)$ is continuous for each $r\in[0,T]$ and
\beg{align}\label{Yn-un-1}
Y_r^{n,t,x}=u_n(r,X_r^{t,x}),\ Z_r^{n,t,x}=\nn u_n(r,X_r^{t,x})\si(r,X_r^{t,x}),\ \ \d r\otimes\d\P\text{-a.e.}
\end{align}
By \eqref{sup-un1} and the fact that $Y_t^{n,t,x}$ is deterministic, we have
\beg{align*}
\sup_{n\geq 1,x\in\R^d,t\in [0,T]}\ff {|u_n(t,x)|} {1+|x|^q}<\infty.
\end{align*}
Combining this with  \eqref{cov-yzn} and \eqref{Yn-un-1}, we get
\beg{align}\label{YnY}
\sup_{n\geq 1}|Y_r^{n,t,x}|+|Y_r^{t,x}|\leq C(1+|X_r^{t,x}|^q), \ \  \d r\otimes\d\P\text{-a.e.,}\ x\in\R^d.
\end{align}

Next, we want to establish an analogue of \eqref{YnY} for $Z_r^{n,t,x}$ and $Z_r^{t,x}$.
For this, define
$$A_s^{[n]}:=(T-s)^\be\sup_{x\in\R^d}\ff {|\nn u_n(s,x)|} {1+|x|^{q\vee 1}},\ s\in[t,T],\ \
B_t^{[n]}:=\sup_{s\in [t,T]}A_s^{[n]}.$$
Using an argument from step 1 in the proof of \cite[Theorem 3.2]{ZJF},
we shall prove that $B_t^{[n]}$ has an upper bound independent of $n$.
Note first that from assertion (1), there holds a gradient type formula in $L^p(\P)$ for $Y_s^{n,t,x}$ for any $s\in[t,T]$:
\beg{align*}
\nn_v Y^{n,t,x}_s= \E\left[g_n(X_T^{t,x})M_T^{t,s}(x,v)+ \int_s^T f_n(r,\Th_r^{n,t,x})M_r^{t,s}(x,v)\d r  \Big| \sF_s^t \right].
\end{align*}
Then letting $s=t$ and using \eqref{in-K-pog} with $k_0=0$, \eqref{CC} with $s=t$, \eqref{Add-8(The2.3)}, \eqref{Yn-un-1} and the H\"{o}lder inequality,
we obtain for any $\delta\in(0,T-t)$,
\beg{align*}
\|\nn Y_t^{n,t,x}\|&\leq \E\left[|g_n(X_T^{t,x})|\cdot|M_T^t(x)|\right]+\int_t^T\E |f_n(r,\Th_r^{n,t,x})M_r^t(x)|\d r\\
&\leq C(T-t)^{-\be}\|g_n(X_T^{t,x})\|_2\\
&\quad +C\int_{t}^T(r-t)^{-\be}\left(1+\|X_r^{t,x}\|_2 +\|Y_r^{n,t,x}\|_2 +\|Z_r^{n,t,x}\|_2 \right)\d r\\
&\leq \ff {C(1+|x|^q)} {(T-t)^{ \be}}+C (T-t)^{1-\be}\left(1+\sup_{r\in [t,T]}(\|X_r^{t,x}\|_2+\|Y_r^{n,t,x}\|_2)\right)\\
&\quad +C\left(\int_{t}^{t+\de} \ff { \|\nn u_n(r,X_r^{t,x})\si(r,X_r^{t,x})\|_2 } {(r-t)^{\be}} \d r+\de^{-\be}\int_{t+\de}^T \|Z_r^{n,t,x}\|_2 \d r\right).
\end{align*}
For the last two terms of the right-hand side of the above inequality,
according to the boundedness of $\sigma$, \eqref{mon-est-1} with $k_2=1$ and \eqref{sup-un1}, we deduce that
\beg{align*}
\int_{t}^{t+\de} \ff { \|\nn u_n(r,X_r^{t,x})\si(r,X_r^{t,x})\|_2 } {(r-t)^{\be}} \d r
& \leq\|\si\|_\infty\int_t^{t+\de} \ff {A_r^{[n]}[\E(1+|X_r^{t,x}|^{q\vee1})^2]^\ff 1 2} {(r-t)^{ \be}(T-r)^\be}\d r\\
& \leq \ff {C\de^{1-\be} B_t^{[n]} (1+|x|^{q\vee1}) } {(T-t-\de)^{\be}}
\end{align*}
and
\beg{align*}
\int_{t+\de}^T \|Z_r^{n,t,x}\|_2 \d r\leq (T-t-\de)^{\ff 1 2}\left(\int_{t }^T\E|Z_r^{n,t,x}|^2\d r\right)^{\ff 1 2}
\leq CT^{\ff 1 2} (1+|x|^q) .
\end{align*}
In addition, due to  \eqref{mon-est-1} with $k_2=1$ and again \eqref{sup-un1}, it is easily seen that
$$\sup_{r\in [t,T]}(\|X_r^{t,x}\|_2+\|Y_r^{n,t,x}\|_2)\leq C(1+|x|^{q\vee 1}).$$
Thus, there exists a positive constant $C$ independent of $n$ and $t$ such that
\beg{align*}
\|\nn Y_t^{n,t,x}\|\leq C\left(\ff 1 {(T-t)^{\be}}+1+\ff {\de^{1-\be}B_t^{[n]} } { (T-t-\de)^{\be}}+ \de^{-\be}\right)(1+|x|^{q\vee1}).
\end{align*}
Consequently, we get
\beg{align*}
A_t^{[n]} &= (T-t)^\be\sup_{x\in\R^d}\ff {\|\nn Y^{n,t,x}_t\|} {1+|x|^{q\vee1}}\\
&\leq C \left(1+(T-t)^{\be}+\ff {(T-t)^\be\de^{1-\be}B_t^{[n]} } { (T-t-\de)^{\be}} +\de^{-\be}(T-t)^{\be}\right)\\
&\leq C \left(1+\ff {(T-t)^\be\de^{1-\be}B_t^{[n]}} { (T-t-\de)^{\be}} +\de^{-\be}(T-t)^{\be}\right).
\end{align*}
Letting $\de=\ff {T-t} m$ with some $m\geq 2$ chosen later, we have
\beg{align*}
A_t^{[n]}&\leq C \left(1+m^{\be}+2^{\be}(T-t)^{1-\be} m^{-(1-\be)}B_t^{[n]}\right)\\
&\leq  \tld C(1+m^{\be}+ m^{-(1-\be)}B_t^{[n]})
\end{align*}
for some $\tld C>0$ independent of $n$ and $t$.
Here we have used $\ff{T-t}{T-t-\delta}\leq 2$.
Since $B_t^{[n]}$ is non-increasing with respect to $t$, it follows that
$$B_t^{[n]}\leq \tld C (1+m^{\be}+ m^{-(1-\be)}B_t^{[n]}).$$
Taking $m$ such that $\tld C m^{-(1-\be)}<1$, we arrive at
\beg{align}\label{B-sup}
\sup_{n\geq 1}B_t^{[n]}\leq \ff {\tld C (1+m^{\be})} {1-\tld C m^{-(1-\be)}}<\infty,
\end{align}
which means that $B_t^{[n]}$ has an upper bound independent of $n$.
Combining this with \eqref{Yn-un-1} and \eqref{B-sup},
we may find a constant $C>0$ independent of $n$  such that for any $x\in\R^d$,
\beg{align}\label{ZTr0}
|Z_r^{n,t,x}|&\leq |\nn u_n(r,X_r^{t,x})|\cdot\|\si(r,X_r^{t,x})\|\cr
&\leq \ff{\|\sigma\|_\infty B_t^{[n]}(1+|X_r^{t,x}|^{q\vee 1})}{(T-r)^{\be}}
\leq \ff {C(1+|X_r^{t,x}|^{q\vee 1})} {(T-r)^{\be}},\ \ \d r\otimes\d\P\text{-a.e.}
\end{align}
Since \eqref{cov-yzn} ensures that $Z_r^{n,t,x}$ converges  to $Z_r^{t,x}$ in the measure $\d r\otimes\d\P$, we have
\beg{align}\label{ZTr}
|Z_r^{t,x}|\leq \ff {C(1+|X_r^{t,x}|^{q\vee 1})} {(T-r)^{\be}},\ \ \d r\otimes\d\P\text{-a.e.},\ x\in\R^d,
\end{align}
which, along with \eqref{ZTr0}, is the desired analogue of \eqref{YnY} for $Z_r^{n,t,x}$ and $Z_r^{t,x}$.

By \eqref{CC}, \eqref{YnY} and \eqref{ZTr0}, we obtain as in \eqref{supfM} that
\beg{align*}
\left\|\E\left[ \sup_{n\geq 1}|f(r,\Th_r^{n,t,x})  M_r^{t,s}(x,v) |\Big| \sF_s^t \right]\right\|_p\leq \ff {C} {(T-r)^\be(r-s)^\be}(1+|x|^{
q\vee 1})(1+|x|^{k_3})|v|,
\end{align*}
which implies that \eqref{insupfM} holds in the present case due to $\be<1$.
Consequently, \eqref{nvf}, \eqref{lim-fM} and \eqref{LipifM} are also true.
Using \eqref{L2-M}, \eqref{YnY} and \eqref{ZTr} and following the same argument as in assertion (1),
one obtains that \eqref{Add-5(The2.3)} holds in $L^p(\P)$ for $p\in[1,2)$.
This completes our proof.
\qed

\subsection{Examples}

In this subsection, we present some examples to illustrate the above results.
As preparation we first state a lemma for verifying \eqref{equ-con} that holds for any $\psi\in \sB_b(\R^d)$ $\P$-a.s. and \eqref{iint-K} in the case of $M_r^t(x)$ being given by a stochastic integral,
whose proof is elementary and therefore postponed to the Appendix.

We set, for $T\geq r\geq t\geq 0$,
\beg{align*}
&\tld W(r-t) =W(r)-W(t),\\
W_{r,m}^{t,0}& =\{\ga\in C([t,r],\R^m)~|~\ga(t)=0\},
\end{align*}
and let  $\tld W([t,r])$ denote the path of $\tld W(\cdot-t)$ on $[t,r]$. Let
$$H_{t,r}: [t,r]\times \R^d\times W_{r,m}^{t,0}\times \R^d\times \R^d\otimes \R^d\ra \R^d\otimes\R^m$$
be measurable and let $\{H_{t,r}(\th, x, \tld W([t,r]), X_\th^{t,x},\nn X_\th^{t,x})\}_{\th\in[t,r]}$ be an adapted process.
We suppose that $\{M_r^t(x)\}_{0\leq t<r\leq T}$ in (H2) has the following form:
\beg{align}\label{equ-M}
M_r^t(x)=\int_t^r H_{t,r}(\th,x,\tld W([t,r]), X_\th^{t,x},\nn X_\th^{t,x})\d W_\th.
\end{align}

\beg{lem}\label{lem-SI}
Assume that the system \eqref{equ-X} and \eqref{Jaco} has a unique nonexplosive strong solution
and satisfies the moment conditions \eqref{mon-est-1} and \eqref{mon-est-2}.
Assume moreover that there exist positive constants $q_1, q_2$ and $C_{q_1,q_2,T}$ such that for all $ 0\leq t\leq s\leq \th <r\leq T$,
\beg{align}\label{ine-HH}
&\left\|H_{s,r}\left(\th,X_s^{t,x}, \tld W([s,r]), X_\th^{s,X_s^{t,x}},\nn X_\th^{s,y}\big|_{y=X_s^{t,x}}\right)\right\|\nonumber \\
&\leq C_{q_1,q_2,T}\left( K_\th(t,s,r, x)  +|X_\th^{t,x}|^{q_1}+\left\|\nn X_\th^{s,y}\Big|_{y=X_s^{t,x}}\right\|^{q_2}\right),
\end{align}
where  $\{K_\th(t,s,r,x)\}_{s\leq \th<r}$ is an adapted process depending on $t,s,r,x$ such that for any $t\leq s<r$,
$\E \int_s^r K^2_\th(t,s,r,\cdot) \d \th$ has polynomial growth.
Then for all $ 0\leq t\leq s <r\leq T$,
\eqref{equ-con} holds for any $\psi\in \sB_b(\R^d)$ $\P$-a.s. and $K(t,s,r,\cdot)$ has polynomial growth.
If in addition,
\beg{align}\label{int-K-1}
\int_s^T\sup_{|x-x_0|\leq R}\left(\E \int_s^r K^2_\th(t,s,r,x) \d \th\right)^{\ff 1 2}\d r<\infty,\ \ 0\leq t\leq s\leq T,
\end{align}
then \eqref{iint-K} holds.
\end{lem}

The first example is an FBSDE with non-degenerate forward SDE.

\beg{exa}\label{exa-non-deg-1}
Consider the FBSDE \eqref{equ-X}-\eqref{BSDE1}.
Assume that for any $r\in [0,T]$, $b(r,\cdot)$ and $\si(r,\cdot)$ are continuously differentiable with
\beg{align}\label{C1-b}
\sup_{r\in[0,T],\ x\in\R^d}\left( \| \nn b(r,x)\|+ \| \nn \si(r,x)\|\right)<\infty,
\end{align}
and that $\si$ is non-degenerate, i.e.  $\si\si^*$ is invertible and
$$\sup_{r\in [0,T],x\in\R^d}\left\|\left(\si^*(\si\si^*)^{-1}\right)(r,x)\right\|<\infty.$$
Assume moreover that (A2) holds  and one of the following assumptions is satisfied:\\
(i) (A1) holds.\\
(ii) (A1$'$) holds and  $\si$ is bounded on $[0,T]\times \R^d$.\\
Then the gradient type formula \eqref{equ-deri} in Theorem \ref{thm-BSDE} holds for
\beg{align*}
M_r^{t,s}(x_0,v)=\int_s^r \ff 1 {r-s}\left\<\left(\si^*(\si\si^*)^{-1}\right)(\th,X_\th^{ t,x_0})\nn_v X_\th^{t,x_0},\d W_\th\right\>.
\end{align*}
\end{exa}


\beg{proof}
By \eqref{C1-b}, we know that the system  \eqref{equ-X} and \eqref{Jaco} has a unique strong solution and the moment conditions \eqref{mon-est-1} and \eqref{mon-est-2} hold with $k_2=1$ and $k_3=0$, respectively.
It is readily checked that (C1) and (C2) hold.

Owing to the condition that $\si$ is non-degenerate and \eqref{C1-b} again,
it is well known that the derivative formula \eqref{Bism-P} holds with $\<M_r^{t}(x),v\>$ given as follows (see, e.g., \cite[Theorem 2.1]{Elworthy&Li94}):
\beg{align}\label{Add0-Ex1}
\<M_r^{t}(x),v\>=\int_t^r\ff {1} {r-t}\left\<\left(\si^*(\si\si^*)^{-1}\right)(\th,X_\th^{t,x})\nn_v X_\th^{t,x},\d W_\th\right\>,\ \ v,x\in\R^d,\ 0\leq t<r\leq T.
\end{align}
Moreover, since $\nn b(t,\cdot) $ and  $\nn \si(t,\cdot)$ are bounded and continuous,
we can derive from \eqref{Bism-P} that $\nn P_{t,r}\ps\in C(\R^d)$ for any $\ps\in C_b^1(\R^d)$.  Hence (H2) holds.

Next, we shall apply Lemma \ref{lem-SI} to verify \eqref{equ-con} and \eqref{iint-K}.
First note that for any $x\in\R^d$ and $0\leq t<r\leq T$, by \eqref{Add0-Ex1} we can write $M_r^{t}(x)$ as follows:
\beg{align}\label{Add1-Ex1}
M_r^{t}(x)=\int_t^r H_{t,r}(\th,x,\tld W([t,r]), X_\th^{t,x},\nn X_\th^{t,x})\d W_\th,
\end{align}
where
\beg{align*}
H_{t,r}(\th,x,\tld W([t,r]), X_\th^{t,x},\nn X_\th^{t,x})=\ff {1} {r-t}\left(\left(\si^*(\si\si^*)^{-1}\right)(\th,X_\th^{t,x})\nn X_\th^{t,x}\right)^*.
\end{align*}
Then we deduce that for any $ 0\leq t\leq s\leq \th <r\leq T$,
\beg{align*}
&\left\|H_{s,r}\left(\th,X_s^{t,x}, \tld W([s,r]), X_\th^{s,X_s^{t,x}},\nn X_\th^{s,y}\big|_{y=X_s^{t,x}}\right)\right\|\cr
&=\ff {1} {r-s}\left\|\left(\left(\si^*(\si\si^*)^{-1}\right)(\th,X_\th^{s,X_s^{t,x}})\nn X_\th^{s,y}\big|_{y=X_s^{t,x}}\right)^*\right\|\cr
&\leq \ff 1 {r-s}\sup_{\th\in[s,r],x\in\R^d}\left\|\left(\si^*(\si\si^*)^{-1}\right)(\th,x)\right\|\cdot\left\|\nn X_\th^{s,y}\Big|_{y=X_s^{t,x}}\right\|\cr
&=:K_\th(t,s,r,x).
\end{align*}
Combining this with \eqref{Add1-Ex1} yields that for any $x\in\R^d$ and $0\leq t\leq s<r\leq T$,
\beg{align}\label{Add2-Ex1}
&K^2(t,s,r,x)=\E \left|M_r^{s}(X_s^{t,x})\right|^2 \cr
&=\E\left|\int_s^r H_{s,r}\left(\th, X_s^{t,x}, \tld W([s,r]), X_\th^{s,X_s^{t,x}},\nn X_\th^{s,y}\Big|_{y=X_s^{t,x}}\right)\d W_\th\right|^2 \cr
& =\E\int_s^r \left\|H_{s,r}\left(\th, X_s^{t,x}, \tld W([s,r]), X_\th^{t,x},\nn X_\th^{s,y}\Big|_{y=X_s^{t,x}}\right)\right\|^2\d \th\cr
&\leq\E\int_s^r K_\th^2(t,s,r,x)\d\th\leq \ff {C_{\si,T} }  {r-s}\sup_{\th\in [s,r],y\in\R^d}\E \|\nn X_\th^{s,y} \|^2\leq\ff {C_{b,\si,T} } {r-s},
\end{align}
where the last inequality is due to \eqref{mon-est-2} with $k_1=2$ and $k_3=0$.
Hence, by Lemma \ref{lem-SI} we get \eqref{equ-con} and \eqref{iint-K}.
In addition, $K(t,s,r,\cdot)$ satisfies \eqref{in-K-pog} with $k_0=0$ and \eqref{CC} holds with $\be=1/2$.

Note that by the continuity and the non-degeneracy of $\si$ and the continuity of $\nn b(r,\cdot)$ and $\nn \si(r,\cdot)$,
it is routine to show that (H3) holds.
Moreover, we have
\beg{align}\label{Add3-Ex1}
M_r^{t,s}(x,v)& =\<M_r^s(X_s^{t,x}),\nn_v X_s^{t,x}\>\cr
&=  \left\<\int_s^r \ff 1 {r-s}\left(\left(\si^*(\si\si^*)^{-1}\right)(\th,X_\th^{s,X_s^{t,x}})\nn X_\th^{s,y}\Big|_{y=X_s^{t,x}}\right)^*\d W_\th,{\nn_v X_s^{t,x}}\right\> \cr
&=  \int_s^r \ff 1 {r-s}\left\<\left(\si^*(\si\si^*)^{-1}\right)(\th,X_\th^{s,X_s^{t,x}})\nn_{\nn_v X_s^{t,x}} X_\th^{s,y}\Big|_{y=X_s^{t,x}},\d W_\th\right\> \cr
&= \ff 1 {r-s} \int_s^r \left\<\left(\si^*(\si\si^*)^{-1}\right)(\th,X_\th^{ t,x })\nn_v X_\th^{t,x},\d W_\th\right\>.
\end{align}

Therefore, taking into account the assumptions (i) and (ii), the assertions follow from Theorem \ref{thm-BSDE} (1) and (2), respectively.
\end{proof}

\beg{rem}\label{Rem-Ex1}
For the model studied in Example \ref{exa-non-deg-1},
two other works \cite{MaZ,Fuhrman&Tessitore02} give some partial papers related to our gradient type formulas.
In \cite[Theorem 4.2]{MaZ} dealing with the case of $l=1$ and $m=d$,
the authors aim to provide a representation formula for $Z$.
Their proof, however, actually leads to the same gradient type formula as that of Example \ref{exa-non-deg-1},
when $g\in C^1_b(\R^d)$ and $f\in C^{0,1}_b([0,T]\times\R^d\times\R\times\R^d)$.
In \cite[Theorem 3.10]{Fuhrman&Tessitore02}, under some regularity conditions on the coefficients the authors proved the following ``Bismut type" formula:
\beg{align*}
\E\nn_v Y^{t,x_0}_s= \E\left[ g(X_T^{t,x_0})\widetilde{M}_T^{t,s}(x_0,v)+ \int_s^T f(r,\Th_r^{t,x_0})\widetilde{M}_r^{t,s}(x_0,v)\d r\right],
\end{align*}
where $\widetilde{M}_r^{t,s}(x_0,v)= \ff 1 {r-t}\int_t^r\<\si^{-1}(\th,X_\th^{ t,x_0})\nn_v X_\th^{t,x_0},\d W_\th\>$,
which is slightly different from $M_r^{t,s}(x_0,v)$ in Example \ref{exa-non-deg-1}.
Both papers use the same arguments, which consist of utilizing a relationship between the Malliavin derivatives $(DX,DY,DZ)$ and partial derivatives in the initial value $x$ of $(X,Y,Z)$ (see, e.g., \cite[Lemma 2.4]{Pardoux&Peng92a} and \cite[Lemma 2.4]{MaZ}), together with the Malliavin integration by parts formula.
In comparison with \cite{MaZ,Fuhrman&Tessitore02},
the gradient type formula in Example \ref{exa-non-deg-1} is obtained under more general conditions,
and moreover our method is more powerful and flexible, which is illustrated by Examples \ref{Add-exa-nondeg2}, \ref{ex.Gru} and  \ref{ex.Ham} below.
We also point out that under non-degenerate assumptions which additionally require $\sigma$ to be independent of $x$,
\cite{Masiero15} extended the above ``Bismut type" formula stated in \cite[Theorem 3.10]{Fuhrman&Tessitore02} to the case of $f$ with quadratic growth with respect to $Z$.
\end{rem}

For the FBSDE stated in Example \ref{exa-non-deg-1}, we can provide an alternative version of the gradient type formula.

\beg{exa}\label{Add-exa-nondeg2}
Let the same assumptions as for Example \ref{exa-non-deg-1} hold.
For each $0\leq t<s\leq T$, define
$$\varrho_{t,s}(\th): =c^{-1}\left(1-e^{-c(s-\th)}\right),\ \ \th\in[t,s),$$
with some positive constant $c$ depending only on $b$ and $\sigma$, and let
$\{G_s^{t,x}(r)\}_{r\in [t,s)}$ satisfy the linear equation on $\R^d\otimes\R^d$:
\beg{align}\label{Add4-Ex2}
G_{s }^{t,x}(r) =I_{d\times d}+\int_t^{r}\left(\nn b(\th,X_\th^{t,x})G_s^{t,x}(\th)-\ff {G_s^{t,x}(\th)} {\varrho_{t,s}(\th)} \right)\d \th+ \int_t^{r}\nn \si(\th,X_\th^{t,x})G_s^{t,x}(\th)\d W_\th,
\end{align}
where  $(\int_t^{r}\nn \si(\th,X_\th^{t,x})G_s^{t,x}(\th)\d W_\th)v:=\int_t^{r}\nn_{G_s^{t,x}(\th)v}\si(\th,X_\th^{t,x})\d W_\th, v\in\R^d.$
Then the gradient type formula \eqref{equ-deri} in Theorem \ref{thm-BSDE} holds for
\beg{align}\label{M-exa-2}
M_r^{t,s}(x_0,v)=\int_s^r\ff 1 {\varrho_{s,r}(\th)}\left\<\left(\si^*(\si\si^*)^{-1}\right)(\th,X_\th^{t,x_0})G_r^{s,X_s^{t,x_0}}(\th)\nn_v X_s^{t,x_0},\d W_\th\right\>.
\end{align}
\end{exa}


\beg{proof}
For any $0\leq t<s\leq T$, we first introduce the following linear equation on $\R^d$:
\beg{align*}
\Ga_{s }^{t,x}(r)&=v+\int_t^r\left(\nn_{\Ga_{s }^{t,x}(\th)}b(\th,X_\th^{t,x})-\ff {\Ga_{s }^{t,x}(\th)} {\varrho_{t,s}(\th)}\right)\d \th\cr
&\qquad+\int_t^r\nn_{\Ga_{s}^{t,x}( \th)}\si(\th,X_\th^{t,x})\d W_\th,\ \ r\in[t,s).
\end{align*}
Obviously, our assumptions imply that the equation has a unique solution $\{\Ga_s^{t,x}(r)\}_{r\in[t,s)}$,
and moreover, $\Ga_{s }^{t,x}(r)= G_{s }^{t,x}(r) v$ for any $r\in[t,s)$.
According to \cite[Theorem 3.5.1 or Theorem 4.3.7]{WB},
the derivative formula \eqref{Bism-P} holds with $\<M_r^{t}(x),v\>$ given as follows:
\beg{align*}
\<M_r^{t}(x),v\>=\int_t^r\ff 1 {\varrho_{t,r}(\th)}\left\< \left(\si^*(\si\si^*)^{-1}\right)(\th,X_\th^{t,x})\Ga_r^{t,x}(\th),\d W_\th\right\>,
\ \ v,x\in\R^d,\ 0\leq t<r\leq T.
\end{align*}
This, along with $\Ga_r^{t,x}(\th)= G_r^{t,x}(\th) v$ for $\th\in[t,r)$, implies that for any $x\in\R^d$ and $0\leq t< r\leq T$,
\beg{align*}
M_r^{t}(x)=\int_t^r H_{t,r}(\th,x,\tld W([t,r]), X_\th^{t,x},\nn X_\th^{t,x})\d W_\th,
\end{align*}
where
\beg{align*}
H_{t,r}(\th,x,\tld W([t,r]), X_\th^{t,x},\nn X_\th^{t,x})=
\ff 1 {\varrho_{t,r}(\th)}\left(\left(\si^*(\si\si^*)^{-1}\right)(\th,X_\th^{t,x})G_r^{t,x}(\th)\right)^*.
\end{align*}
Then we get that for any $ 0\leq t\leq s\leq \th <r\leq T$,
\beg{align*}
&\left\|H_{s,r}\left(\th, X_s^{t,x},\tld W([s,r]),X_\th^{s,X_s^{t,x}},\nn X_\th^{s,y}\Big|_{y=X_s^{t,x}}\right)\right\|\\
&=\ff 1 {\varrho_{s,r}(\th)}\left\|\left(\left(\si^*(\si\si^*)^{-1}\right)(\th,X_\th^{s,X_s^{t,x}})G_r^{s,X_s^{t,x}}(\th)\right)^*\right\|\\
&\leq \ff 1 {\varrho_{s,r}(\th)}\sup_{\th\in[s,r],x\in\R^d}\left\|\left(\si^*(\si\si^*)^{-1}\right)(\th,x)\right\|\cdot\left\|G_r^{s,X_s^{t,x}}(\th)\right\| \cr
&=:K_\th(t,s,r,x).
\end{align*}
Similar to \eqref{Add2-Ex1} in Example  \ref{exa-non-deg-1}, we have
\beg{align}\label{Add5-Ex2}
K^2(t,s,r,x)=\E \left|M_r^{s}(X_s^{t,x})\right|^2&\leq \E\int_s^r K^2_\th(t,s,r,x) \d \th\cr
&\leq C_{\si,T}\E\int_s^r \ff {\|G_r^{s,X_s^{t,x}}(\th)\|^2}{\varrho^2_{s,r}(\th)}\d \th\cr
&\leq C_{b,\si,T}\left(1+\ff  1 {1-e^{-c(r-s)}}\right),
\end{align}
where the last inequality is due to \cite[Theorem 3.5.1 and (3.65), or Theorem 4.3.7 and (4.51)]{WB}.
Then, \eqref{equ-con} and \eqref{iint-K} follow from Lemma \ref{lem-SI} and $K(t,s,r,\cdot)$ satisfies \eqref{in-K-pog} with $k_0=0$.
Besides, a direct and easy computation shows that \eqref{CC} also holds with $\be=1/2$ as in Example \ref{exa-non-deg-1}.


Finally, we are to verify that (H3) holds and $M_r^{t,s}(x,v)$ is given by \eqref{M-exa-2}.
As in \eqref{Add3-Ex1} of Example \ref{exa-non-deg-1}, we obtain
\beg{align*}
M_r^{t,s}(x,v)& =\<M_r^s(X_s^{t,x}),\nn_v X_s^{t,x}\>\\
&=\left\<\int_s^r \ff 1 {\varrho_{s,r}(\th)}\left(\left(\si^*(\si\si^*)^{-1}\right)(\th,X_\th^{s,X_s^{t,x}})G_r^{s,X_s^{t,x}}(\th)\right)^*
\d W_\th,{\nn_v X_s^{t,x}}\right\> \\
&=\int_s^r\ff 1 {\varrho_{s,r}(\th)}\left\<\left(\si^*(\si\si^*)^{-1}\right)(\th,X_\th^{t,x})G_r^{s,X_s^{t,x}}(\th)\nn_v X_s^{t,x},\d W_\th\right\>,
\end{align*}
which is consistent with \eqref{M-exa-2}.
As for (H3), it is standard to prove that its first part \eqref{L2-M} holds due to the continuity and non-degeneracy of $\si$ and the continuity of $\nn b(r,\cdot)$ and  $\nn \si(r,\cdot)$.
Let us now deal with its second part \eqref{add-Hy}.
Following the same arguments as in \cite[Lemma 4.3.8]{WB} and \cite[Theorem 4.3.7 and (4.51)]{WB}, we derive that for any $p>1$,
\beg{align}\label{Add2-Ex2}
\E\sup_{\th\in[s,r)}\|G_r^{s,X_s^{t,x}}(\th)\|^p\leq C_{b,\sigma,T}
\end{align}
and
\beg{align}\label{Add3-Ex2}
\sup_{0<\de\leq \ff {r-s} 2}\int_s^r\ff {\E\|G_{r-\de}^{s,X_s^{t,x}}(\th)\|^p}
{\varrho_{s,r-\de}^p(\th)}\1_{[s,r-\de)}(\th)\d\th
\leq \sup_{0<\de\leq \ff {r-s} 2}C_{b,\sigma,T}\left(1+\ff  1 {1-e^{-c(r-\delta-s)}} \right)<\infty,
\end{align}
respectively.
By \eqref{Add4-Ex2}, we have for $\th\in[s,r-\delta)$,
\beg{align*}
G_{r-\de}^{s,X_s^{t,x}}(\th)-G_r^{s,X_s^{t,x}}(\th)&=\int_s^\th\nn b(\tilde{\th},X_{\tilde{\th}}^{s,X_s^{t,x}})\left(G_{r-\de}^{s,X_s^{t,x}}(\tilde{\th})-G_r^{s,X_s^{t,x}}(\tilde{\th})\right)\d\tilde{\th}\cr
&\quad-\int_s^\th\left(\ff{G_{r-\de}^{s,X_s^{t,x}}(\tilde{\th})}{\varrho_{s,r-\delta}(\tilde{\th})}
-\ff{G_r^{s,X_s^{t,x}}(\tilde{\th})}{\varrho_{s,r}(\tilde{\th})}\right)\d\tilde{\th}\cr
&\quad+\int_s^\th\nn\si(\tilde{\th},X_{\tilde{\th}}^{s,X_s^{t,x}})\left(G_{r-\de}^{s,X_s^{t,x}}(\tilde{\th})-G_r^{s,X_s^{t,x}}(\tilde{\th})\right)\d W_{\tilde{\th}}.
\end{align*}
Applying the It\^{o} formula and using the boundedness of $\nn b$ and $\nn\si$ lead to
\beg{align*}
&\E\|G_{r-\de}^{s,X_s^{t,x}}(\th)-G_r^{s,X_s^{t,x}}(\th)\|^2\cr
&\leq C_{b,\si}\E\int_s^\th\|G_{r-\de}^{s,X_s^{t,x}}(\tilde{\th})-G_r^{s,X_s^{t,x}}(\tilde{\th})\|^2\d\tilde{\th}
-2\E\int_s^\th\ff{\|G_{r-\de}^{s,X_s^{t,x}}(\tilde{\th})-G_r^{s,X_s^{t,x}}(\tilde{\th})\|^2}{\varrho_{s,r-\delta}(\tilde{\th})}
\d\tilde{\th}\cr
&\quad+2\E\int_s^\th\left(\ff 1{\varrho_{s,r}(\tilde{\th})}-\ff 1 {\varrho_{s,r-\delta}(\tilde{\th})}\right)
\left\<G_{r-\de}^{s,X_s^{t,x}}(\tilde{\th})-G_r^{s,X_s^{t,x}}(\tilde{\th}),G_r^{s,X_s^{t,x}}(\tilde{\th})\right\>\d\tilde{\th}\cr
&\leq C_{b,\si}\int_s^\th\E\|G_{r-\de}^{s,X_s^{t,x}}(\tilde{\th})-G_r^{s,X_s^{t,x}}(\tilde{\th})\|^2\d\tilde{\th}\cr
&\quad+\int_s^\th\left(\ff 1{\varrho_{s,r}(\tilde{\th})}-\ff 1 {\varrho_{s,r-\delta}(\tilde{\th})}\right)^2\E\|G_r^{s,X_s^{t,x}}(\tilde{\th})\|^2\d\tilde{\th},\ \ \th\in[s,r-\delta),
\end{align*}
where the second inequality is due to the positivity of $\varrho_{s,r-\delta}(\cdot)$ and the Young inequality.
The Gronwall inequality implies that for any $\th\in[s,r-\delta)$,
\beg{align*}
\E\|G_{r-\de}^{s,X_s^{t,x}}(\th)-G_r^{s,X_s^{t,x}}(\th)\|^2
\leq C_{b,\si,T} \int_s^\th\left(\ff 1{\varrho_{s,r}(\tilde{\th})}-\ff 1 {\varrho_{s,r-\delta}(\tilde{\th})}\right)^2\E\|G_r^{s,X_s^{t,x}}(\tilde{\th})\|^2\d\tilde{\th}.
\end{align*}
Then using \eqref{Add2-Ex2} and taking into account the fact that $\lim_{\de\ra 0^+}\varrho_{s,r-\de}(\tilde{\th})=\varrho_{s,r}(\tilde{\th})$,
by the dominated convergence theorem we obtain that for every $r'\in[s,r)$,
\beg{align*}
\lim_{\de\ra 0^+,\de<r-r'}\sup_{s\leq \th \leq r'}\E \|G_{r-\de}^{s,X_s^{t,x}}(\th)-G_r^{s,X_s^{t,x}}(\th)\|^2=0,
\end{align*}
Consequently, we obtain that for every $\th\in[s,r)$,
$$\P\text{-}\lim_{\de\ra 0^+}\left(\ff{G_{r-\de}^{s,X_s^{t,x}}(\th)}{\varrho_{s,r-\de}(\th)}\1_{[s,r-\de)}(\th)-\ff{G_r^{s,X_s^{t,x}}(\th)}{\varrho_{s,r}(\th)}\right)=0.$$
Note that the term involving $\delta$ above is well defined  as long as $\delta$ is small enough.
Thus, combining this with the boundedness of $\si^*(\si\si^*)^{-1}$ and \eqref{Add3-Ex2}, we apply the dominated convergence theorem to get
\beg{align*}
&\quad\lim_{\de\ra 0^+}\E\left|M_{r-\de}^s(X_s^{t,x})-M_r^s(X_s^{t,x})\right|^2\cr
&=\lim_{\de\ra 0^+}\E\bigg|\int_s^r\ff 1 {\varrho_{s,r-\delta}(\th)}\left(\left(\si^*(\si\si^*)^{-1}\right)(\th,X_\th^{s,X_s^{t,x}})G_{r-\delta}^{s,X_s^{t,x}}(\th)\right)^*\1_{[s,r-\de)}(\th)
\d W_\th\cr
&\qquad\qquad\quad-\int_s^r \ff 1 {\varrho_{s,r}(\th)}\left(\left(\si^*(\si\si^*)^{-1}\right)(\th,X_\th^{s,X_s^{t,x}})G_r^{s,X_s^{t,x}}(\th)\right)^*
\d W_\th\bigg|^2\cr
&=\lim_{\de\ra 0^+}\E\int_s^r\bigg\|\left(\si^*(\si\si^*)^{-1}\right)(\th,X_\th^{s,X_s^{t,x}})
\bigg(\ff{G_{r-\de}^{s,X_s^{t,x}}(\th)}{\varrho_{s,r-\de}(\th)}\1_{[s,r-\de)}(\th)-\ff{G_r^{s,X_s^{t,x}}(\th)}{\varrho_{s,r}(\th)}\bigg)\bigg\|^2\d\th\cr
&=0,
\end{align*}
which means that \eqref{add-Hy} holds.

Therefore, similar to Example \ref{exa-non-deg-1}, the assertions follow from Theorem \ref{thm-BSDE} (1) and (2), respectively.
\end{proof}

We may also consider the following FBSDE with forward Gruschin type process.

\beg{exa}\label{ex.Gru}
Let $d=d_1+d_2,x=(x^{(1)},x^{(2)})$ and $W=(W^{(1)},W^{(2)})$ be a Brownian motion on $\R^{d_1+d_2}$,
and consider the FBSDE \eqref{equ-X}-\eqref{BSDE1},
where the forward SDE \eqref{equ-X} is of the following form:
\beg{equation}\label{Gru}
\beg{cases}
X_s^{(1),t,x}=x^{(1)}+\int_t^s \d W^{(1)}_r,\ \ s\in[t,T],\\
X_s^{(2),t,x}=x^{(2)}+\int_t^s\si(X_r^{(1),t,x})\d W^{(2)}_r,\ \ s\in[t,T],
\end{cases}
\end{equation}
with $\si\in C^1(\R^{d_1};\R^{d_2}\otimes\R^{d_2})$ which might be degenerate.
Assume that there exist $\al\in [1,\ff {d_1} 2+1)$ and constants $a_1,a_2>0$ such that
\beg{align}\label{non-deg-0}
\|\si(x)\|\geq a_1|x|^\al,\ \ \ \|\si(x)\|+\|\nn\si(x)\|\cdot|x|\leq a_2|x|^\al,\ \  x\in\R^{d_1},
\end{align}
and that (A1) and (A2) are satisfied.
Then the gradient type formula \eqref{equ-deri} in Theorem \ref{thm-BSDE} holds for
\beg{align*}
&M_r^{t,s}(x_0,v) =\ff {\<v^{(1)},W^{(1)}_{r}-W^{(1)}_s\>} {r-s}\\
&-\mathrm{Tr}\left(Q_{s,r}^{-1} (x^{(1)}_0+W^{(1)}_s-W^{(1)}_t)\int_s^r\ff {r-\th} {r-s}\left((\nn_{v^{(1)}}\si)\si^*\right)(x^{(1)}_0+W^{(1)}_\th-W^{(1)}_t)\d \th\right)\\
&+\left\<Q_{s,r}^{-1}(x^{(1)}_0+W^{(1)}_s-W^{(1)}_t)\left[
\nn_{v}X_s^{(2),t,x_0}+\int_s^r\ff {r-\th} {r-s}\nn_{v^{(1)}}\si(x^{(1)}_0+W^{(1)}_\th-W^{(1)}_t)\d W^{(2)}_\th\right],\right.\\
&\left.\qquad\int_s^r\si(x^{(1)}_0+W^{(1)}_\th-W^{(1)}_t)\d W^{(2)}_\th\right\>
\end{align*}
provided that $x^{(1)}_0\neq0$ when $s=t$.
Here $v=(v^{(1)},v^{(2)})\in\R^{d_1+d_2}$ and
\beg{align}\label{Add2-Ex3}
Q_{s,r}(y):=\int_s^r (\si\si^*)(y+W^{(1)}_\th-W^{(1)}_s)\d \th,\ \ y\in\R^{d_1}.
\end{align}
\end{exa}


\beg{proof}
By \eqref{Gru} and \eqref{non-deg-0}, one can show that for any $k\geq 1$,
\beg{align*}
\sup_{0\leq t\leq s\leq T}\E |X^{t,x}_s|^k&\leq C_{k,\al,a_2,T}\left(1+|x|^k+|x^{(1)}|^{k\al}\right)
\end{align*}
and
\beg{align*}
\sup_{0\leq t\leq s\leq T}\E |\nn X^{t,x}_s|^k&\leq C_{k,\al,a_2,T}\left(1+|x^{(1)}|^{k(\al-1)}\right),
\end{align*}
which imply that the moment conditions \eqref{mon-est-1} and  \eqref{mon-est-2} hold with $k_2=\al$ and $k_3=\al-1$, respectively.
Moreover, \eqref{add-der} and \eqref{con-nn-P} in (C1) follow easily from \eqref{Gru} and \eqref{non-deg-0}.
Hence (C1) holds.

Now, we are to verify that (C2) holds.
Let
$$\ph_n(x):=\left[1+\left((|x|-n)^+\right)^2\right]^{-\ff 1 2},\ \ x\in\R^{d_1},\ n\geq 1,$$
and $\si_n(x):=\si(\ph_n(x)x)$.
Then $\ph_n(x)=1$ for any $|x|\leq n$, and
\beg{align}\label{si-n-si}
\lim_{n\ra+\infty}\sup_{|x|\leq R}\|\si( \ph_n(x)x)-\si( x)\|=0,\ \ R>0.
\end{align}
Moreover, we have $\nn\si_n( \cdot)\in C(\R^{d_1};\R^{d_2}\otimes\R^{d_2}\otimes\R^{d_1})$ with
\beg{align}\label{Add1-Ex3}
\nn\si_n( x)
&=\nn\si(\ph_n(x)x)\left(\nn \ph_n(x)\otimes x+\ph_n(x)I_{d_1\times\d_1}\right)\cr
&=\nn\si(\ph_n(x)x)\left[-\ff { (|x|-n)^+(x \otimes x)} {\left(1+\left((|x|-n)^+\right)^2\right)^{\ff 3 2}|x|}+\ph_n(x)I_{d_1\times\d_1}\right].
\end{align}
Consequently, by \eqref{non-deg-0} we get
\beg{align*}
\|\nn\si_n\|_\infty&\leq a_2\sup_{x\in\R^{d_1}}\left[\left(\ff {(|x|-n)^+|x|} {\left(1+\left((|x|-n)^+\right)^2\right)^{\ff 3 2}} +1\right)|\ph_n(x)x|^{\al-1}\right]\leq a_2(n+2)^\al,
\end{align*}
i.e. \eqref{Revise-C1-b} holds.
For any $(x_n)_{n\geq 1}$ with $\lim_{n\ra+\infty}x_n=x$,
let $(X_s^{n,t,x_n})_{s\in[t,T]}$ be the solution of the system \eqref{Gru} with $\si$ replaced by $\si_n$.
Since $X_s^{(1),n,t,x_n}-X_s^{(1),t,x}=x^{(1)}_n-x^{(1)}$ and $\nn X_s^{(1),n,t,x_n}=\nn X_s^{(1),t,x}$, we only need to
check \eqref{app-X1} and \eqref{app-X2} for $X_s^{(2),n,t,x}$.
By the It\^o formula and the B-D-G inequality, we have
\beg{align}\label{Xn-X-G}
&\E\sup_{s\in [t,T]} \left|X_s^{(2),n,t,x_n}-X_s^{(2),t,x}\right|^2\nonumber\\
&\leq C\left( \left|x_n^{(2)}-x^{(2)}\right|^2+\int_t^T\E \left\|\si_n(X_r^{(1),n,t,x_n})-\si(X_r^{(1),t,x})\right\|^2\d r\right).
\end{align}
By \eqref{si-n-si} and the continuity of $\si$, we obtain that for each $r\in[t,T]$,
\beg{align*}
&\quad\lim_{n\ra+\infty}\left\|\si_n(X_r^{(1),n,t,x_n})-\si(X_r^{(1),t,x})\right\|\\
&=\lim_{n\ra+\infty}\left\|\si_n(x^{(1)}_n+W_r^{(1)}-W_t^{(1)})-\si(x^{(1)}+W_r^{(1)}-W_t^{(1)}) \right\|\\
&\leq\lim_{n\ra+\infty}\bigg[\left\|\si_n(x^{(1)}_n+W_r^{(1)}-W_t^{(1)})-\si(x^{(1)}_n+W_r^{(1)}-W_t^{(1)}) \right\|\\
&\qquad\qquad\quad+\left\|\si(x^{(1)}_n+W_r^{(1)}-W_t^{(1)})-\si(x^{(1)}+W_r^{(1)}-W_t^{(1)}) \right\|\bigg]\\
&=0.
\end{align*}
Combining this with \eqref{Xn-X-G} and applying the dominated convergence theorem, we conclude that
\beg{align*}
\lim_{n\ra+\infty}\E\sup_{s\in [t,T]} \left|X_s^{(2),n,t,x_n}-X_s^{(2),t,x}\right|^2=0,
\end{align*}
i.e. \eqref{app-X1} holds.
As for \eqref{app-X2}, owing to \eqref{Add1-Ex3} and the continuity of $\nn\si$ , we derive that for each $r\in[t,T]$,
\beg{align*}
&\quad\lim_{n\ra+\infty}\left\|\nn\si_n(X_r^{(1),n,t,x_n})-\nn \si(X_r^{(1),t,x})\right\|\cr
&=\lim_{n\ra+\infty}\left\|\nn\si(\ph_n(y)y)\left(\nn \ph_n(y)\otimes y+\ph_n(y)I_{d_1\times d_1}\right)\Big|_{y=X_r^{(1),n,t,x_n}}-\nn \si(X_r^{(1),t,x})\right\|\cr
&=0.
\end{align*}
Then using the It\^o formula and the B-D-G inequality and applying the dominated convergence theorem again, we deduce that
\beg{align*}
&\quad\lim_{n\ra+\infty} \E\sup_{s\in [t,T]} \left\|\nn X_s^{(2),n,t,x_n}-\nn X_s^{(2),t,x}\right\|^2\\
&\leq \lim_{n\ra+\infty}\int_t^T\E\left\|\nn\si_n(X_r^{(1),n,t,x_n})-\nn\si(X_r^{(1), t,x })\right\|^2\d r\\
&=0,
\end{align*}
where the choice of the dominating function can be determined by \eqref{non-deg-0} and \eqref{Add1-Ex3}.
That is, \eqref{app-X2} holds.
Besides, it is clear that \eqref{app-si} holds.
Hence (C2) holds.

Note that, according to \cite[Lemma 3.1]{W14} and \eqref{non-deg-0},
we obtain that for any $0\leq t<r\leq T, z\in\R^{d_1}$ and $p>0$, $Q_{t,r}(z)$ defined in \eqref{Add2-Ex3} is invertible and
\beg{align}\label{mon-Q-1}
\E \|Q_{t,r}^{-1}(z)\|^p\leq \ff C {(r-t)^p [|z|^2+(r-t)]^{p\al}}.
\end{align}
By \cite[Theorem 1.1 and Corollary 1.2]{W14},
it follows that for any $0\leq t<r\leq T$ and $x=(x^{(1)},x^{(2)})$, $v=(v^{(1)},v^{(2)})\in\R^{d_1+d_2}$,
\beg{align*}
\nn_v P_{t,r}\psi(x^{(1)},x^{(2)})=\E\left[\psi(X^{(1),t,x}_r,X^{(2),t,x}_r)\<M_r^t(x),v\>\right],\ \ \psi\in C_b^1(\R^{d_1+d_2})
\end{align*}
with
\beg{align*}
&\<M_r^t(x),v\>=\ff {\<v^{(1)},W^{(1)}_{r}-W^{(1)}_t\>} {r-t}\cr
&-\mathrm{Tr}\left(Q_{t,r}^{-1}(x^{(1)})\int_t^r\ff {r-\th} {r-t}\ \left((\nn_{v^{(1)}}\si)\si^*\right)(x^{(1)}+W^{(1)}_\th-W^{(1)}_t)\d\th\right)\\
&+\left\<Q_{t,r}^{-1}(x^{(1)})\left[v^{(2)}+\int_t^r\ff {r-\th} {r-t}\nn_{v^{(1)}}\si(x^{(1)}+W^{(1)}_\th-W^{(1)}_t)\d W^{(2)}_\th\right],\right.\\
&\qquad\left.\int_t^r\si(x^{(1)}+W^{(1)}_\th-W^{(1)}_t)\d W^{(2)}_\th\right\>.
\end{align*}
Moreover, since
$$Q_{s,r}^{-1}(y)-Q_{s,r}^{-1}(x)=Q^{-1}_{s,r}(x)\left(Q_{s,r}(x)-Q_{s,r}(y)\right)Q^{-1}_{s,r}(y),$$
by the dominated convergence theorem and \eqref{mon-Q-1} we get
\beg{align}\label{Add9-Ex3}
\lim_{y\ra x}\E\|Q_{s,r}^{-1}(y)-Q_{s,r}^{-1}(x)\|^2=0.
\end{align}
Then we conclude that  $\nn P_{t,r}\ps\in C(\R^d)$ for any $\ps\in C_b^1(\R^d)$.
Hence (H2) holds.

Now, let $\{e_i\}_{i=1}^{d_1+d_2}$  be the canonical ONB of $\R^{d_1+d_2}$.
Since $Q_{s,r}$ is independent of $\sF_s^t$ for any $0\leq t\leq s<r\leq T$,
it is obvious that for every $e_i$ and $y\in\R^{d_1+d_2}$, $\langle M_r^s(y),e_i\rangle$ is also independent of $\sF_s^t$,
Consequently, it is readily checked that \eqref{equ-con} holds due to the fact that
$ M_r^s(y)=\sum_{i=1}^{d_1+d_2}\langle M_r^s(y),e_i\rangle e_i$.

Next, we shall verify that $K(t,s,r,\cdot)$ satisfies \eqref{in-K-pog} and \eqref{iint-K}.
First note that,
\beg{align}\label{mon-M-2'}
 K^2(t,s,r,x)=\E |M_r^{ s}( X_s^{t,x})|^2 &= \sum_{i=1}^{d_1+d_2}\E\left[\left(\E|\<M_r^s(X_s^{t,x}),e_i\>|^2\Big|\sF_s^t\right)\right]\nonumber\\
&=\sum_{i=1}^{d_1+d_2}\E\left[\left(\E|\<M_r^{s}(y),e_i\>|^2\right)\Big|_{y=X_s^{t,x} }\right].
\end{align}
Along the same lines as in the proof of \cite[Corollary 1.2]{W14}, we have
\beg{align*}
\E |\<M_r^{s}(y),e_i\>|^2\leq C\left(\ff {|e_i^{(1)}|^2} {r-s}+\ff {|e_i^{(2)}|^2} {(r-s)(|y^{(1)}|^2+(r-s))^\al}\right).
\end{align*}
Plugging this into \eqref{mon-M-2'}, we obtain
\beg{align}\label{Add3-Ex3}
 K^2(t,s,r,x)&\leq \ff {C} {r-s}\left( 1 + \E \ff 1{\left(|x^{(1)} +W^{(1)}_s-W^{(1)}_t|^2+(r-s)\right)^{\al}} \right).
\end{align}
Consequently, it is easy to see that $K(t,s,r,\cdot)$ satisfies \eqref{in-K-pog} with $k_0=0$ and
$C(t,s,r)=\ff {C} {\sqrt{r-s}}(1+\ff {1} {(r-s)^{\al/2}})$.
Since $\al\in [1,\ff {d_1} 2+1)$, we can take $\be\in (\al-1,\ff {d_1}2)$, which implies that
$$d_1-1-2\be>-1 \ \ \text{and}  \ \ (1+\al-\be)/2<1.$$
Then by \eqref{Add3-Ex3}, we deduce that there exists $R>0$ such that
\beg{align*}
&\int_s^T\sup_{|x-x_0|\leq R}K(t,s,r,x)\d r\\
& \leq C\int_s^T\ff 1{\sqrt{r-s}}\d r+C\int_s^T\sup_{|x-x_0|\leq R}\left(\E \ff {(r-s)^{-1}} {(|x^{(1)} +W^{(1)}_s-W^{(1)}_t|^2+(r-s))^{\al}}\right)^{\ff 1 2}\d r\\
& \leq C+C\int_s^T\sup_{|x-x_0|\leq R}\left(\E \ff {(r-s)^{-1}|x^{(1)} +W^{(1)}_s-W^{(1)}_t|^{-2\be}} {(|x^{(1)}+W^{(1)}_s-W^{(1)}_t|^2+(r-s))^{\al- \be}}\right)^{\ff 1 2}\d r\\
& \leq C+C\left(\sup_{|x-x_0|\leq R}\E|x^{(1)} +W^{(1)}_s-W^{(1)}_t|^{-2\be}\right)^{\ff 1 2}\int_s^T(r-s)^{-(1+\al-\be)/2}\d r\\
&<\infty,
\end{align*}
where we have used the condition that $x^{(1)}_0\neq0$ when $s=t$.
We emphasize that here we only consider $\al\geq \be$ in the second to last inequality since the integral in the second inequality is finite when $\al<\be$.
Hence \eqref{iint-K} holds.

Finally, we are to verify that \eqref{L2-M} holds.
Owing to \eqref{non-deg-0}, \eqref{mon-Q-1}, \eqref{Add9-Ex3} and since $\nn \si\in C(\R^{d_1},\R^{d_2}\otimes\R^{d_2}\otimes\R^{d_1})$,
we can apply the dominated convergence theorem to get
\beg{align}\label{Add5-Ex3}
\lim_{y\ra x}\E |M^{s}_r(y)-M^{s}_r(x)|^2=0.
\end{align}
Consequently, using the independence of $\{M_r^s(x)\}_{x\in\R^{d_1+d_2}}$ and $\sF_s^t$ and applying the dominated convergence theorem again,
we obtain
\beg{align*}
\lim_{y\ra x}\E|M^{s}_r(X_s^{t,y})-M^{s}_r(X_s^{t,x})|^2&=\lim_{y\ra x}\E\left[\left(\E|M^{s}_r(z_1)-M^{s}_r(z_2)|^2\right)\Big|_{z_1=X_s^{t,y},z_2=X_s^{t,x}}\right]\\
&=0.
\end{align*}
Indeed, for $\varphi(z_1,z_2):=\E|M^{s}_r(z_1)-M^{s}_r(z_2)|^2$, it suffices to prove that
for every sequence $(y_n)_{n\geq1}$ with $\lim_{n\ra+\infty}y_n=x$,
there exists a subsequence $(y_{n_k})_{k\geq1}$ such that
$$\lim_{k\ra+\infty}\E\varphi(X_s^{t,y_{n_k}},X_s^{t,x})=0.$$
By \eqref{Gru}-\eqref{non-deg-0} and the It\^{o} formula, it is easy to see that $\lim_{n\ra+\infty}\E|X_s^{t,y_n}-X_s^{t,x}|^2=0$.
Then there exists a subsequence $(y_{n_k})_{k\geq1}$ such that $\lim_{k\ra+\infty}X_s^{t,y_{n_k}}=X_s^{t,x},\ \P$-a.s.
Due to the definition of $M^{s}_r(\cdot)$ and \eqref{mon-Q-1}, one can show that for each $p>1$,
$$\sup_{k\geq1}\E\varphi^p(X_s^{t,y_{n_k}},X_s^{t,x})<\infty$$
So, by the dominated convergence theorem and \eqref{Add5-Ex3} we get the desired result.
Hence \eqref{L2-M} holds.

Therefore, the assertion  follows from Theorem \ref{thm-BSDE} (1).
\end{proof}

\beg{rem}\label{Rem-Ex3}
It is not hard to extend the above result to the FBSDE with the following general Gruschin type process:
\beg{equation*}
\beg{cases}
X_s^{(1),t,x}=x^{(1)}+\int_t^sb_1(X_r^{(1),t,x})\d r+\int_t^s\si_1(X_r^{(1),t,x}) \d W^{(1)}_r,\ \ s\in[t,T],\\
X_s^{(2),t,x}=x^{(2)}+\int_t^sb_2(X_r^{(1),t,x})\d r+\int_t^s\si_2(X_r^{(1),t,x})\d W^{(2)}_r,\ \ s\in[t,T],
\end{cases}
\end{equation*}
where $b_1\in C_b^1(\R^{d_1};\R^{d_1})$ and $b_2\in C_b^1(\R^{d_1};\R^{d_2})$,
$\si_1\in C_b^1(\R^{d_1};\R^{d_1}\otimes\R^{d_1})$ is invertible and $\si^{-1}_1$ is bounded,
 $\si_2\in C^1(\R^{d_1};\R^{d_2}\otimes\R^{d_2})$  might be degenerate.
\end{rem}

Now we apply Corollary \ref{cor2} to the FBSDE with forward stochastic Hamiltonian system.
Let $d=d_1+d_2$ and $W$ be a $d_2$-dimensional Brownian motion.
Consider the FBSDE \eqref{equ-X}-\eqref{BSDE1}, where the forward SDE \eqref{equ-X} is of the following form:
\beg{align*}
X_s^{t,x}&=x+\int_t^sb(r,X_r^{t,x})\d r+\int_t^s(0,\sigma(r)\d W_r)\\
&=x+\int_t^s(BX^{(2),t,x}_r,\tilde{b}(r,X^{(1),t,x}_r,X^{(2),t,x}_r))\d r+\int_t^s(0,\sigma(r)\d W_r),\ s\in[t,T],
\end{align*}
where $B$ is a $d_1\times d_2$-matrix with rank $d_1$, $\tilde{b}:[0,T]\times\R^{d_1}\times\R^{d_2}\ra\R^{d_2}$ and
$\si(r)$ is an invertible $d_2\times d_2$-matrix which is continuous in $r\in[0,T]$.
Since $B$ has rank $d_1$, we know that $d_2\geq d_1$ and for every $y^{(1)}\in\R^{d_1}$,
$$B^{-1}y^{(1)}:=\left\{y^{(2)}\in\R^{d_2}: By^{(2)}=y^{(1)}\right\}\neq\emptyset.$$
We set
$$|B^{-1}y^{(1)}|:=\inf\left\{|y^{(2)}|: y^{(2)}\in B^{-1}y^{(1)}\right\}.$$
Then it is easy to show that
$$\|B^{-1}\|:=\sup\left\{|B^{-1}y^{(1)}|: y^{(1)}\in\R^{d_1}\ \text{and}\  |y^{(1)}|\leq 1\right\}<\infty.$$
In addition, let $L^\infty([0,T],C^1_b(\R^{d_1}\otimes\R^{d_2};\R^{d_2}))$ be the space of all a.e. bounded measurable functions
$\varphi:[0,T]\ra C^1_b(\R^{d_1}\otimes\R^{d_2};\R^{d_2})$ such that
$$\mathop{\mathrm{ess\,sup}}\limits_{s\in[0,T],y\in\R^{d_1}\otimes\R^{d_2}}(|\varphi(s,y)|+\|\nn\varphi(s,y)\|)<\infty.$$
\beg{exa}\label{ex.Ham}
Assume that $\tilde{b}\in L^\infty([0,T],C^1_b(\R^{d_1}\otimes\R^{d_2};\R^{d_2}))$ and that (A1) and (A2) are satisfied.
Then for $0\leq t<s\leq T$, the Bismut type formula \eqref{equ-Bis} in Corollary \ref{cor2} holds for
\beg{align*}
\<M^t_r(x_0),v\>=\int_t^r \left\<\si^{-1}(\th)\left(\chi''_{t,r}(\th)\tilde{v}^{(2)}
-\kappa''_{t,r}(\th)v^{(2)}+\nn_{\Xi_{t,r}(v,\tilde{v}^{(2)},\th)}\tilde{b}(\th,X_\th^{t,x_0})\right),\d W_\th\right\>,
\end{align*}
where $v=(v^{(1)},v^{(2)})\in\R^{d_1+d_2}, \tilde{v}^{(2)}\in B^{-1} v^{(1)}$ and
\beg{align}\label{Add1-Ex4}
\chi_{t,r}(\th)&=\ff {(\th-t)^2(3r-t-2\th)} {(r-t)^3},\ \ \kappa_{t,r}(\th)=\ff {(\th-t)(r-\th)^2} {(r-t)^2},\ \th\in [t,r],\\
\Xi_{t,r}(v,\tilde{v}^{(2)},\th)&=\left(\left(1-\chi_{t,r}(\th)\right)v^{(1)}+\kappa_{t,r}(\th)Bv^{(2)},
\kappa'_{t,r}(\th)v^{(2)}-\chi'_{t,r}(\th)\tilde{v}^{(2)}\right).\label{Add2-Ex4}
\end{align}
\end{exa}

\beg{proof}
By our assumptions, one can show that the moment conditions \eqref{mon-est-1} and \eqref{mon-est-2} hold with $k_2=1$ and $k_3=0$, respectively.
Furthermore, it is easy to check that (C1) and (C2) hold.

According to \cite[Theorem 2.2]{Guillin&Wang12a}, the derivative formula \eqref{Bism-P} holds with $\<M_r^{t}(x),v\>$ given as follows:
for any $x,v\in\R^{d_1+d_2}, \ 0\leq t< r\leq T$,
\beg{align*}
\<M_r^{t}(x),v\>&=\int_t^r \left\<\si^{-1}(\th)\left(\chi''_{t,r}(\th)\tilde{v}^{(2)}
-\kappa''_{t,r}(\th)v^{(2)}+\nn_{\Xi_{t,r}(v,\tilde{v}^{(2)},\th)}\tilde{b}(\th,X_\th^{t,x})\right),\d W_\th\right\>\cr
&=:\int_t^r \left\<\left(H_{t,r}(\th,x,\tld W([t,r]), X_\th^{t,x},\nn X_\th^{t,x})\right)^*v,\d W_\th\right\>,
\end{align*}
where $\chi_{t,r}(\th),\kappa_{t,r}(\th)$ and $\Xi_{t,r}(v,\tilde{v}^{(2)},\th)$ are defined in \eqref{Add1-Ex4} and \eqref{Add2-Ex4}.
Moreover, similarly to Example \ref{exa-non-deg-1}, we obtain that $\nn P_{t,r}\ps\in C(\R^d)$ for any $\ps\in C_b^1(\R^d)$.  Hence (H2) holds.

Notice that by a direct calculation, we have
\beg{align*}
&\left|\left(H_{t,r}(\th,x,\tld W([t,r]), X_\th^{t,x},\nn X_\th^{t,x})\right)^*v\right|^2\\
&\leq C_T\sup_{\th\in [0,T]}\|\si^{-1}(\th)\|^2 \left[\ff {\|B^{-1}\|^2|v^{(1)}|^2} {(r-t)^4} +\ff {|v^{(2)}|^2} {(r-t)^2}\right.\\
&\qquad\qquad \left.+\|\nn \tilde{b}(\th,\cdot)\|_\infty^2 \left(|v^{(1)}|^2+(r-t)^2|Bv^{(2)}|^2+|v^{(2)}|^2+\ff {\|B^{-1}\|^2|v^{(1)}|^2} {(r-t)^2}\right)\right]\\
&\leq C_{B,\tilde{b},\si,T}\left(1+\ff 1{(r-t)^4}\right)|v|^2.
\end{align*}
Consequently, we get for any $x,v\in\R^{d_1+d_2}$,
\beg{align}\label{Add3-Ex4}
\E|\<M_r^{t}(x),v\>|^2&=\E\int_t^r\left|\left(H_{t,r}(\th,x,\tld W([t,r]), X_\th^{t,x},\nn X_\th^{t,x})\right)^*v\right|^2\d\th\cr
&\leq C_{B,\tilde{b},\si,T}\left(r-t+\ff 1{(r-t)^3}\right)|v|^2.
\end{align}
Then we obtain that for any $0\leq t<s<r\leq T$ and $R>0$,
\beg{align*}
\int_s^T\sup_{|x-x_0|\leq R} K(t,t,r,x)\d r&=\int_s^T\sup_{|x-x_0|\leq R}\left(\E|M_r^t(x)|^2\right)^\ff 1 2\d r\\
&=\int_s^T\sup_{|x-x_0|\leq R}\left(\sum_{i=1}^{d_1+d_2}\E|\<M_r^t(x),e_i\>|^2\right)^\ff 1 2\d r\\
&\leq C_{T,B,\si}\int_s^T\left(\sqrt{r-t}+\ff 1{(r-t)^{\ff 3 2}}\right)\d r\\
&<\infty,
\end{align*}
where $\{e_i\}_{i=1}^{d_1+d_2}$ is the canonical ONB of $\R^{d_1+d_2}$.
Hence, $K(t,t,r,\cdot)$ satisfies \eqref{in-K-pog} with $k_0=0$ and \eqref{iint-K'} holds.

Finally, since $\chi_{t,r}$ and $\kappa_{t,r}$ are deterministic and $\nn\tilde{b}$ is bounded,
it is readily verified that \eqref{L2-M} (with $s=t$) holds.

Therefore, the assertion follows from Corollary \ref{cor2} (1).
\end{proof}

\beg{rem}\label{Rem-Ex4}
(i) If $d_1=d_2, B=\si(r)=I_{d_1\times d_1}$ and
$$\tilde{b}(r,x)=-\nn\varphi(x^{(1)})-Cx^{(2)},$$
the forward SDE is knows as ``stochastic damping Hamiltonian system" in probability (see, e.g., \cite{Bakry etal08a} and \cite{Wu01a}).

(ii) We mention that $1/(r-t)^3$ in \eqref{Add3-Ex4} seems to be optimal,
meaning that the term $\E\int_s^T f(r,\Th_r^{t,x_0})\<M_r^t(x_0),v\>\d r$ in \eqref{equ-Bis} might be ill-defined when $t=s$.
For example, let $t=s$ and suppose that $f$ is bounded, then by \eqref{Add3-Ex4}
\beg{align*}
\left|\E\int_s^T f(r,\Th_r^{t,x_0})\<M_r^t(x_0),v\>\d r\right|\leq C_{B,\tilde{b},\si,T,f}|v|\int_s^T\left(\sqrt{r-s}+\ff 1{(r-s)^{\ff 3 2}}\right)\d r=\infty.
\end{align*}
So, the Bismut type formula is thus not available in the case of $t=s$.
\end{rem}

\section{Applications to McKean-Vlasov FBSDEs and related PDEs}

In this section, we will apply the gradient type and Bismut type formulas to the study of McKean-Vlasov FBSDEs.
Our goal is twofold.
Firstly, we want to establish the representation formulae for the control solution $Z^{t,\xi}_s$,
which enables us to derive its path regularity.
Secondly, we are concerned with the gradient estimates for the solution of a nonlocal PDE of mean-field type,
which is associated with a McKean-Vlasov FBSDE.

Let $W$ be a $d$-dimensional Brownian motion,
$(\sF_t)_{t\in[0,T]}$ and $\{\sF_r^t\}_{0\leq t\leq r\leq T}$ the corresponding filtrations as before,
and denote by $\sG\subset\sF$ a sub-$\si$-algebra independent of $W$.
For an initial datum $(t,\xi)\in[0,T]\times L^2(\Omega\rightarrow\R^d, \sG,\P)$ we consider the following McKean-Vlasov FBSDE:
for $s\in[t,T]$,
\begin{equation}\label{3.1}
\left\{
\begin{array}{ll}
X_s^{t,\xi}=\xi+\int_t^sb(r,X_r^{t,\xi},\P_{X_r^{t,\xi}})\d r+\int_t^s\sigma(r,X_r^{t,\xi},\P_{X_r^{t,\xi}})\d W_r,\\
Y_s^{t,\xi}=g(X_T^{t,\xi},\P_{X_T^{t,\xi}})+\int_s^Tf(r,X_r^{t,\xi},Y_r^{t,\xi},Z_r^{t,\xi},\P_{(X_r^{t,\xi},Y_r^{t,\xi},Z_r^{t,\xi})})\d r+\int_s^TZ_r^{t,\xi}\d W_r,
\end{array}\right.
\end{equation}
where $\P_{X_r^{t,\xi}}$ and $\P_{(X_r^{t,\xi},Y_r^{t,\xi},Z_r^{t,\xi})}$ denote the distributions of $X_r^{t,\xi}$ and $(X_r^{t,\xi},Y_r^{t,\xi},Z_r^{t,\xi})$ under $\P$, respectively.
Note that the coefficients $b,\sigma$ and $g,f$ depend on the law of the solution,
and the above forward-backward system is decoupled, as the first equation does not contain the solution $(Y^{t,\xi},Z^{t,\xi})$ of the second one.

Let $\mathcal{P}_2(\R^d)$ be the collection of all probability measures on $\R^d$ with finite second moment.
Define the $2$-Wasserstein distance on $\mathcal{P}_2(\R^d)$ by
\beg{align*}
\mathbb{W}_2(\mu,\nu):=\inf\limits_{\pi\in\mathfrak{C}(\mu,\nu)}\left(\int_{\R^d\times\R^d}|x-y|^2\pi(\d x,\d y)\right)^\frac{1}{2},
\end{align*}
where $\mathfrak{C}(\mu,\nu)$ is the set of all probability measures on $\R^d\times\R^d$ with marginal laws $\mu$ and $\nu$.
Then $(\mathcal{P}_2(\R^d),\mathbb{W}_2)$ is a Polish space.
In addition, let $\delta_0$ be the Dirac measure at $0$ and denote by
$\mathcal{S}^2([t,T];\R^k)$ (respectively, $\mathcal{H}^2([t,T];\R^k)$)  the space of all $\R^k$-valued and $(\mathcal{F}_t\vee\sG)_{0\leq t\leq T}$-adapted continuous processes (respectively, predictable processes) $\phi=(\phi_s)_{s\in[t,T]}$
with $\E(\sup_{s\in[t,T]}|\phi_s|^2)<\infty$ (respectively, $\E(\int_t^T|\phi_s|^2\d s)<\infty$).
Throughout this section, we will make the following assumptions on the coefficients $b,\sigma$ and $g,f$.
\beg{enumerate}
\item[(\~{H}1)] $b:[0,T]\times\R^d\times\mathcal{P}_2(\R^d)\rightarrow\R^d$ and $\sigma:[0,T]\times\R^d\times\mathcal{P}_2(\R^d)\rightarrow\R^d\otimes\R^d$
are continuous in time and there exists a constant $L_1>0$ such that
\beg{align*}
&|b(r,x_1,\mu_1)-b(r,x_2,\mu_2)|+\|\sigma(r,x_1,\mu_1)-\sigma(r,x_2,\mu_2)\|\cr
&\leq L_1(|x_1-x_2|+\mathbb{W}_2(\mu_1,\mu_2)),\ \ r\in[0,T],\ x_1,x_2\in\R^d,\ \mu_1,\mu_2\in\mathcal{P}_2(\R^d).
\end{align*}
Besides,  $b(r,0,\delta_0)$ and $\sigma(r,0,\delta_0)$ are bounded functions of $r\in[0,T]$.
\item[(\~{H}2)] $g:\R^d\times\mathcal{P}_2(\R^d)\rightarrow\R$ and $f:[0,T]\times\R^d\times\R\times\R^d\times\mathcal{P}_2(\R^d\times\R\times\R^d)\rightarrow\R$ are
Lipschitz continuous, that is, there exists a constant $L_2>0$ such that
for all $r\in[0,T],x_1,x_2\in\R^d,y_1,y_2\in\R,z_1,z_2\in\R^d,\mu_1,\mu_2\in\mathcal{P}_2(\R^d\times\R\times\R^d),\mu,\tilde{\mu}\in\mathcal{P}_2(\R^d)$,
$$|g(x_1,\mu)-g(x_2,\tilde{\mu})|\leq L_2(|x_1-x_2|+\mathbb{W}_2(\mu,\tilde{\mu}))$$
and
$$|f(r,x_1,y_1,z_1,\mu_1)-f(r,x_2,y_2,z_2,\mu_2)|\leq L_2(|x_1-x_2|+|y_1-y_2|+|z_1-z_2|+\mathbb{W}_2(\mu_1,\mu_2)).$$
Besides, assume that $|g(0,\delta_0)|+\sup_{r\in [0,T]}|f(r,0,0,0,\delta_0)|<\infty.$
\end{enumerate}

We observe that under (\~{H}1) and (\~{H}2), the system \eqref{3.1} has a unique solution $(X^{t,\xi},Y^{t,\xi},Z^{t,\xi})\in
\mathcal{S}^2([t,T];\R^d)\times\mathcal{S}^2([t,T];\R)\times\mathcal{H}^2([t,T];\R^d)$, see, for instance \cite[Theorem A.1]{Li18a},
where $b,\sigma$ and $f$ are independent of the time variable and the driven noises are a Brownian motion and an independent Poisson random measure.
Since the proof of \cite[Theorem A.1]{Li18a} also applies to the present case, we omit the proof here.
In order to solve a class of PDEs of mean-field type (see \eqref{3.6} below),
the authors of \cite{Buckdahn&L&Peng&ainer17a} and \cite{Li18a}
introduce SDEs accompanying the system \eqref{3.1} to deal with the homogeneous and non-homogeneous case, respectively,
which then allows to provide probabilistic representations for the solution of \eqref{3.6} (see also \cite{Crisan&McMurray18a} for the homogeneous case).
In this part, we will adopt this idea to establish the representation formulas for $Z^{t,\xi}$
and the gradient estimates for \eqref{3.6} via the gradient type and Bismut type formulas for BSDEs obtained in the last section.
To formulate the accompanying SDEs, let $(X^{t,\xi},Y^{t,\xi},Z^{t,\xi})$ be the solution of \eqref{3.1}.
Now, for any $x\in\R^d$, consider the following SDEs: for $s\in[t,T]$,
\beg{align}\label{3.2}
X_s^{t,x,\P_\xi}=x+\int_t^sb(r,X_r^{t,x,\P_\xi},\P_{X_r^{t,\xi}})\d r+\int_t^s\sigma(r,X_r^{t,x,\P_\xi},\P_{X_r^{t,\xi}})\d W_r
\end{align}
and
\beg{align}\label{3.3}
Y_s^{t,x,\P_\xi}&=g(X_T^{t,x,\P_\xi},\P_{X_T^{t,\xi}})
+\int_s^Tf(r,X_r^{t,x,\P_\xi},Y_r^{t,x,\P_\xi},Z_r^{t,x,\P_\xi},\P_{(X_r^{t,\xi},Y_r^{t,\xi},Z_r^{t,\xi})})\d r\nonumber\\
&\quad+\int_s^TZ_r^{t,x,\P_\xi}\d W_r.
\end{align}
Since the distribution dependence in the coefficients of \eqref{3.2} is $\P_{X_r^{t,\xi}}$,
rather than that of the solution to itself, it becomes a classical SDE.
Then, it is easy to see that under (\~{H}1) there exists a unique solution $X^{t,x,\P_\xi}$ to \eqref{3.2}.
Similarly, under (\~{H}1) and (\~{H}2) \eqref{3.3} also has a unique solution $(Y^{t,x,\P_\xi},Z^{t,x,\P_\xi})$.
Moreover,
\beg{align}\label{3.5}
X^{t,\xi}=X^{t,x,\P_\xi}|_{x=\xi}, \ \ Y^{t,\xi}=Y^{t,x,\P_\xi}|_{x=\xi}, \ \ Z^{t,\xi}=Z^{t,x,\P_\xi}|_{x=\xi}.
\end{align}
For further details, we, for instance, refer to \cite{Buckdahn&L&Peng&ainer17a,Crisan&McMurray18a,Li18a}.

Below is our first result which asserts that the control solution $Z_s^{t,\xi}$ of the system \eqref{3.1} admits a representation formula.
\beg{thm}\label{Theor3.1}
Let assumptions (\~{H}1) and (\~{H}2) hold.
Assume moreover that for each $r\in [0,T]$ and $\mu\in\mathcal{P}_2(\R^d)$, $b(r,\cdot,\mu)$ and $\si(r,\cdot,\mu)$ are continuously differentiable with
\beg{align*}
\sup_{(r,x,\mu)\in[0,T]\times\R^d\times\mathcal{P}_2(\R^d)}\left( \| \nn b(r,x,\mu)\|+ \| \nn\si(r,x,\mu)\|\right)<\infty,
\end{align*}
and that $\sigma$ is non-degenerate.
Then for each $s\in[t,T)$,
\beg{align}\label{Theor3.1-1}
Z_s^{t,\xi}&
=\E\left[g(X_T^{t,\xi},\P_{X_T^{t,\xi}})N_T^{t,s}(\xi)+ \int_s^T f(r,X_r^{t,\xi},Y_r^{t,\xi},Z_r^{t,\xi},\P_{(X_r^{t,\xi},Y_r^{t,\xi},Z_r^{t,\xi})})N_r^{t,s}(\xi)\d r \Big| \sF_s^t\vee\sG \right]\nonumber\\
&\qquad\times(\nabla X_s^{t,\xi,\P_\xi})^{-1}\sigma(s,X_s^{t,\xi},\P_{X_s^{t,\xi}}).
\end{align}
In particular,
\beg{align}\label{Theor3.1-2}
Z_t^{t,\xi}
&=\E\bigg[g(X_T^{t,\xi},\P_{X_T^{t,\xi}})N_T^{t,t}(\xi)
+\int_t^Tf(r,X_r^{t,\xi},Y_r^{t,\xi},Z_r^{t,\xi},\P_{(X_r^{t,\xi},Y_r^{t,\xi},Z_r^{t,\xi})})N_r^{t,t}(\xi)\d r\Big|\sG\bigg]\nonumber\\
&\qquad\times\sigma(t,\xi,\P_\xi)
\end{align}
Here for $t\leq s<r\leq T, \nabla X_s^{t,\xi,\P_\xi}=\nabla X_s^{t,x,\P_\xi}\big|_{x=\xi}$ and
\beg{align*}
N_r^{t,s}(\xi)=
\int_s^r  \ff 1 {r-s}\left(\left(\si^*(\si\si^*)^{-1}\right)(\th,X_\th^{t,\xi},\P_{X_\th^{t,\xi}})\nn X_\th^{t,\xi,\P_\xi}\right)^*\d W_\th.
\end{align*}
\end{thm}


\beg{proof}
We first assume that $g(\cdot,\mu)$ and $f(r,\cdot,\cdot,\cdot,\mu)$ are also continuously differentiable for any $r\in [0,T]$ and $\mu\in\mathcal{P}_2(\R^d)$.
For fixed $(t,\xi)$, considering the coefficients $\tilde{b}(s,x):=b(s,x,\P_{X_s^{t,\xi}}), \tilde{\sigma}(s,x):=\sigma(s,x,\P_{X_s^{t,\xi}})$
and $\tilde{g}(x):=g(x,\P_{X_T^{t,\xi}}),\tilde{f}(s,x,y,z):=f(s,x,y,z,\P_{(X_s^{t,\xi},Y_s^{t,\xi},Z_s^{t,\xi})})$
and then applying Lemma \ref{rep-YZ-X} or \cite[Theorem 3.1]{MaZ} to the system \eqref{3.2}-\eqref{3.3},
we see that there exists $\tilde{u}:[0,T]\times\R^d\ra\R$ such that $\nn\tilde{u}(r,\cdot)$ is continuous for each $r\in[0,T]$ and
$$Y_s^{t,x,\P_\xi}=\tilde{u}(s,X_s^{t,x,\P_\xi}),\ \ Z_s^{t,x,\P_\xi}=\nn\tilde{u}(s,X_s^{t,x,\P_\xi})\tilde{\sigma}(s,X_s^{t,x,\P_\xi}),\ \ s\in[t,T],\ \P\text{-}a.s.$$
Then in view of Example \ref{exa-non-deg-1}, we obtain that for any $s\in[t,T)$,
\beg{align}\label{Pf-Theor3.1-1}
Z_s^{t,x,\P_\xi}&=\nn Y_s^{t,x,\P_\xi}(\nabla X_s^{t,x,\P_\xi})^{-1}\tilde{\sigma}(s,X_s^{t,x,\P_\xi})\cr
&=\E\left[\tilde{g}(X_T^{t,x,\P_\xi})N_T^{t,s}(x)+ \int_s^T \tilde{f}(r,X_r^{t,x,\P_\xi},Y_r^{t,x,\P_\xi},Z_r^{t,x,\P_\xi})N_r^{t,s}(x)\d r  \Big| \sF_s^t \right]\cr
&\qquad\times(\nabla X_s^{t,x,\P_\xi})^{-1}\tilde{\sigma}(s,X_s^{t,x,\P_\xi}),
\end{align}
where for any $t\leq s<r\leq T$,
\beg{align}\label{Pf-Theor3.1-2}
N_r^{t,s}(x)=
\int_s^r  \ff 1 {r-s}\left(\left(\tilde{\si}^*(\tilde{\si}\tilde{\si}^*)^{-1}\right)(\th,X_\th^{ t,x,\P_\xi})\nn X_\th^{t,x,\P_\xi}\right)^*\d W_\th.
\end{align}
By the same approximation argument as in the proof of Theorem \ref{thm-BSDE} (1) or \cite[Theorem 4.2, Page 1410-1411]{MaZ},
we may extend the formula \eqref{Pf-Theor3.1-1} to all $f$ and $g$ satisfying only (\~{H}2).
Therefore, letting $x=\xi$ in \eqref{Pf-Theor3.1-1}-\eqref{Pf-Theor3.1-2} and taking into account \eqref{3.5} and the definitions of $\tilde{\si}, \tilde{g}$ and $\tilde{f}$, we derive \eqref{Theor3.1-1}.
We stress here that based on the regular conditional probability and the monotone class theorem,
the fact that $\sF_s^t$ is independent of $\sG$ converts $\E[\cdot|\sF_s^t]$ in \eqref{Pf-Theor3.1-1} into
$\E[\cdot|\sF_s^t\vee\sG]$ in \eqref{Theor3.1-1}.
The proof is standard, and hence, omitted here.

Finally, it is easy to see that \eqref{Theor3.1-2} follows from \eqref{Theor3.1-1} since the process $\nn X_s^{t,\xi,\P_\xi}:=\nn X_s^{t,x,\P_\xi}|_{x=\xi},s\in[t,T]$ satisfies the following SDE:
\beg{align*}
\nn X_s^{t,\xi,\P_\xi}=I_{d\times d}+\int_t^s\nn b(r,X_r^{t,\xi},\P_{X_r^{t,\xi}})\nn X_r^{t,\xi,\P_\xi}\d r
+\int_t^s\nn\sigma(r,X_r^{t,\xi},\P_{X_r^{t,\xi}})\nn X_r^{t,\xi,\P_\xi}\d W_r,
\end{align*}
which is a slight variation of \cite[Remark 4.1]{Buckdahn&L&Peng&ainer17a}.
\end{proof}

\beg{rem}\label{Rem1-Th3.1}
(i) If the coefficients $b, \si$ and $g, f$ of FBSDE \eqref{3.1} are assumed to be functions with no dependence on the measure and $\xi=x$ (namely, $\xi$ is deterministic),
then it is readily checked that our representation formula \eqref{Theor3.1-1} coincides with that of \cite[Theorem 4.2]{MaZ}.
So, our formula \eqref{Theor3.1-1} is a generalized version of \cite[Theorem 4.2]{MaZ}.

(ii) Combining the above proof and Example \ref{Add-exa-nondeg2}, we can alternatively derive that
\beg{align*}
N_r^{t,s}(\xi)=
\int_s^r  \ff 1 {\varrho_{s,r}(\th)}\left(\left(\si^*(\si\si^*)^{-1}\right)(\th,X_\th^{t,\xi},\P_{X_\th^{t,\xi}})
G_r^{s,y}(\th)\big|_{y=X_s^{t,\xi}}\nn X_s^{t,\xi,\P_\xi}\right)^*\d W_\th,
\end{align*}
where $\{G_r^{s,y}(\th)\}_{\th\in[s,r)}$ solves \eqref{Add4-Ex2} for $\nn\tilde{b}$ and $\nn\tilde{\si}$ replacing $\nn b$ and $\nn \si$ there, respectively.
Besides, under conditions on $\si$ similar to Example \ref{ex.Gru} and the assumptions (\~{H}1) and (\~{H}2),
the representation formula for $Z^{t,\xi}$ can be extended to McKean-Vlasov FBSDE with forward Gruschin type process.
\end{rem}

Now let $V(t,x,\P_\xi):=Y_t^{t,x,\P_\xi}$.
Notice that $V(t,x,\P_\xi)$ is deterministic since $Y_t^{t,x,\P_\xi}$ is measurable both with respect to $\mathcal{F}_t$ and $\mathcal{F}_T^t$.
By the It\^{o} formula for distribution dependent functionals (see, e.g., \cite[Theorem 7.1]{Buckdahn&L&Peng&ainer17a} and \cite[Theorem 2.1]{Li18a}),
it is shown in \cite[Theorem 9.2]{Li18a} that under appropriate regularity conditions on the coefficients of $b,\sigma,g$ and $f$,
$V$ is the solution to the following nonlocal PDE of mean-field type:
 \begin{equation}\label{3.6}
\left\{
\begin{array}{ll}
(\partial_t+\mathbb{L}_{b,\sigma})V(t,x,\P_\xi)=f\left(x,V(t,x,\P_\xi),(\nn V\sigma)(t,x,\P_\xi),\P_{\left(\xi,V(t,\xi,\P_\xi),(\nn V\sigma)(t,\xi,\P_\xi)\right)}\right),\\
 \  \ \ \ \ \  \  \ \ \ \ \  \  \ \ \ \ \ \  \ \ \ \ \ \  \ \ \ \ \  \  \ \ \ \ \ \ \ \ \ \ \ \  \ \ \ \ \ \  \ \
 (t,x,\xi)\in[0,T]\times\R^d\times L^2(\Omega\rightarrow\R^d, \mathcal{F}_t,\P),\\
V(T,x,\P_\xi)=g(x,\P_\xi),\ \ (x,\P_\xi)\in\R^d\times\mathcal{P}_2(\R^d),
\end{array} \right.
\end{equation}
where the operator $\mathbb{L}_{b,\sigma}$, introduced in \cite{Buckdahn&L&Peng&ainer17a}, is defined as follows: for any $(t,x,\mu)\in[0,T]\times\R^d\times\mathcal{P}_2(\R^d)$,
\beg{align*}
\mathbb{L}_{b,\sigma}V(t,x,\mu)&=\langle b,\nn V\rangle(t,x,\mu)
+\frac{1}{2}\mathrm{Tr}\left(\sigma\sigma^*\nn^2 V\right)(t,x,\mu)\cr
&\quad+\int_{\R^d}\big\langle b(t,y,\mu),\partial_\mu V(t,x,\mu)(y)\big\rangle\mu(\d y)\cr
&\quad+\frac{1}{2}\int_{\R^d}\mathrm{Tr}\Big((\sigma\sigma^*)(t,y,\mu)\nn(\partial_\mu V(t,x,\mu)(\cdot))(y)\Big)\mu(\d y).
\end{align*}
Here $\nn V, \nn^2 V$ and $\partial_\mu V$ denote the gradient, the Hessian and the Lions derivative of $V$, respectively.
We remark that, in contrast to the classical case,
the last two terms in the description of $\mathbb{L}_{b,\sigma}V(t,x,\mu)$ involve the derivatives with respect to the measure variable.
Recall that the notion of this type of derivative was introduced
by P.-L. Lions in his lectures \cite{Cardaliaguet13} at the Coll\`{e}ge de France,
see also \cite{Albeverio&Kondratiev&Rockner96,Buckdahn&L&Peng&ainer17a,Ren&Wang19a}.
Below we will provide the uniform gradient estimates for PDE \eqref{3.6}.

\beg{thm}\label{Theor3.2}
Let the assumptions in Theorem \ref{Theor3.1} hold and $V(t,x,\P_\xi)=Y_t^{t,x,\P_\xi}$ be the solution to the PDE \eqref{3.6} (which is assumed to exist and to be unique).
Then there exists a constant $C_{L_1,L_2,T}>0$ such that
\beg{align}\label{Theor3.2-1}
|\nn V(t,x,\P_\xi)|\leq C_{L_1,L_2,T}\left(\ff 1 {\sqrt{T-t}}+\sqrt{T-t}\right)(1+|x|+\|\xi\|_2).
\end{align}
\end{thm}

\beg{proof}
We consider the functions $\tilde{b},\tilde{\si}$ and $\tilde{g},\tilde{f}$ introduced in the proof of Theorem \ref{Theor3.1}.
With these notations in hand, we can apply Corollary \ref{cor2} and \eqref{Add1-Ex1} in Example \ref{exa-non-deg-1} to the system \eqref{3.2}-\eqref{3.3}, and we obtain that
for any $0\leq t\leq s<T$ and $v\in\R^d$,
\beg{align}\label{Pf-Theor3.2-1}
\nn_v\E Y^{t,x,\P_\xi}_s&=\E\bigg[\tilde{g}(X_T^{t,x,\P_\xi})
\<M_T^t(x),v\>+\int_s^T \tilde{f}(r,X_r^{t,x,\P_\xi},Y_r^{t,x,\P_\xi},Z_r^{t,x,\P_\xi})\<M_r^t(x),v\>\d r\bigg],
\end{align}
where for each $t<r\leq T$,
\beg{align*}
M_r^t(x)&=\int_t^r\ff {1} {r-t}\left(\left(\tilde{\si}^*(\tilde{\si}\tilde{\si}^*)^{-1}\right)(\th,X_\th^{t,x,\P_\xi})\nn X_\th^{t,x,\P_\xi}\right)^*\d W_\th\\
&=\int_t^r\ff {1} {r-t}\left(\left(\si^*(\si\si^*)^{-1}\right)(\th,X_\th^{t,x,\P_\xi},\P_{X_\th^{t,\xi}})\nn X_\th^{t,x,\P_\xi}\right)^*\d W_\th.
\end{align*}
Since $Y^{t,x,\P_\xi}_t$ is deterministic, setting $s=t$ in \eqref{Pf-Theor3.2-1} implies
\beg{align}\label{Pf-Theor3.2-2}
&\quad\nn V(t,x,\P_\xi)=\nn Y^{t,x,\P_\xi}_t\nonumber\\
&=\E\bigg[g(X_T^{t,x,\P_\xi},\P_{X_T^{t,\xi}})M_T^t(x)
+\int_t^T  f(r,X_r^{t,x,\P_\xi},Y_r^{t,x,\P_\xi},Z_r^{t,x,\P_\xi},\P_{({X_r^{t,\xi},Y_r^{t,\xi},Z_r^{t,\xi}})})M_r^t(x)\d r\bigg].
\end{align}

In order to obtain an upper bound for $\nn V(t,x,\P_\xi)$,
we start with the following moment estimates of the solutions to the equations \eqref{3.1}-\eqref{3.3}:
\beg{align*}
\E\left[\sup\limits_{s\in[t,T]}|X_s^{t,x,\P_\xi}|^2+\sup\limits_{s\in[t,T]}|Y_s^{t,x,\P_\xi}|^2+\mathop{\mathrm{ess\,sup}}\limits_{s\in[t,T]}|Z_s^{t,x,\P_\xi}|^2\right]
\leq C_{L_1,L_2,T}(1+|x|^2)
\end{align*}
and
\beg{align*}
\E\left[\sup\limits_{s\in[t,T]}|X_s^{t,\xi}|^2+\sup\limits_{s\in[t,T]}|Y_s^{t,\xi}|^2+\mathop{\mathrm{ess\,sup}}\limits_{s\in[t,T]}|Z_s^{t,\xi}|^2\right]
\leq C_{L_1,L_2,T}(1+\|\xi\|_2^2),
\end{align*}
which can be proved in the spirit of the proof of \cite[Theorem 3.3]{MaZ}.
Then by the H\"{o}lder inequality and the Lipschitz continuity of $g$ and $f$, we have
\beg{align}\label{Pf-Theor3.2-5}
\left|\E\left[g(X_T^{t,x,\P_\xi},\P_{X_T^{t,\xi}})M_T^t(x)\right]\right|
&\leq\left[\E\left|g(X_T^{t,x,\P_\xi},\P_{X_T^{t,\xi}})\right|^2\right]^\frac{1}{2}\|M_T^t(x)\|_2\nonumber\\
&\leq C_{L_1,L_2,T}\left[\E\left(1+|X_T^{t,x,\P_\xi}|+\|X_T^{t,\xi}\|_2\right)^2\right]^\frac{1}{2}\|M_T^t(x)\|_2\nonumber\\
&\leq C_{L_1,L_2,T}(1+|x|+\|\xi\|_2)\|M_T^t(x)\|_2
\end{align}
and
\beg{align}\label{Pf-Theor3.2-6}
&\left|\E\int_t^T f(r,X_r^{t,x,\P_\xi},Y_r^{t,x,\P_\xi},Z_r^{t,x,\P_\xi},\P_{({X_r^{t,\xi},Y_r^{t,\xi},Z_r^{t,\xi}})})M_r^t(x)\d r\right|\cr
&\leq C_{L_1,L_2,T}\E\int_t^T \Big[\left(1+|X_r^{t,x,\P_\xi}|+|Y_r^{t,x,\P_\xi}|+|Z_r^{t,x,\P_\xi}|
+\|X_r^{t,\xi}\|_2+\|Y_r^{t,\xi}\|_2+\|Z_r^{t,\xi}\|_2\right)\cr
&\qquad\qquad\qquad\qquad\times|M_r^t(x)|\Big]\d r\cr
&\leq C_{L_1,L_2,T}(1+|x|+\|\xi\|_2)\int_t^T\|M_r^t(x)\|_2\d r.
\end{align}
Note that due to \eqref{Add2-Ex1} (with $s=t$) in Example \ref{exa-non-deg-1}, we can get for $t<r\leq T$,
\beg{align*}
\|M_r^t(x)\|_2\leq \ff {C_{L_1,T}}{\sqrt{r-t}}.
\end{align*}
Therefore, substituting this into \eqref{Pf-Theor3.2-5}-\eqref{Pf-Theor3.2-6} and going back to relation \eqref{Pf-Theor3.2-2},
we get the desired assertion.
\end{proof}

\beg{rem}\label{Rem1-Th3.2}
(i) In \cite[Remark 5.6 and Theorem 5.8]{Crisan&McMurray18a},
by developing integration by parts formulas for McKean-Vlasov SDEs with uniformly elliptic coefficients,
the authors proved the existence of a classical solution to the homogeneous PDE \eqref{3.6} (namely, $f=0$ in \eqref{3.6})
with a class of non-differentiable terminal conditions $g$
(including, for instance, $g(x,\mu)=\int\varphi(x,y)\mu(\d y)$, where $\varphi:\R^d\times\R^d\ra\R$ is continuous with polynomial growth,
more interesting examples can be found in \cite[Example 5.4]{Crisan&McMurray18a}),
and obtained a gradient bound of the form $C(T-t)^{-1/2}(1+|x|+\|\xi\|_2)^q$ with positive constants $C$ and $q$, which essentially coincides with our estimate \eqref{Theor3.2-1}.
Therefore, the result stated in Theorem \ref{Theor3.2} can be regarded as a generalization of that in \cite{Crisan&McMurray18a}.

(ii) One can provide an alternative version of the previous estimate,
whose proof follows the same lines as the one of Theorem \ref{Theor3.2}.
Owing to \eqref{Add5-Ex2} (with $s=t$) in Example \ref{Add-exa-nondeg2}, we have for $t<r\leq T$,
\beg{align*}
\|M_r^{t}(x)\|_2\leq C_{L_1,T}\left(1+\ff  1 {\sqrt{1-e^{-c(r-t)}}}\right).
\end{align*}
Then, we obtain
\beg{align*}
|\nn V(t,x,\P_\xi)|\leq C_{L_1,L_2,T}\left(1+\ff  1 {\sqrt{1-e^{-c(T-t)}}}+\int_t^T\ff  1 {\sqrt{1-e^{-c(r-t)}}}\d r\right)(1+|x|+\|\xi\|_2).
\end{align*}
In particular, as $T$ tends to $t$, we have
\beg{align*}
1+\ff  1 {\sqrt{1-e^{-c(T-t)}}}+\int_t^T\ff  1 {\sqrt{1-e^{-c(r-t)}}}\d r=O\left(\ff 1 {\sqrt{T-t}}\right),
\end{align*}
which means that this type of gradient estimate has the same order as the right-hand side of \eqref{Theor3.2-1} when $T$ goes to $t$.

(iii) For PDEs \eqref{3.6} associated with a McKean-Vlasov FBSDE with forward Gruschin type process, it follows from \eqref{Add3-Ex3} (with $s=t$) in Example \ref{ex.Gru} that
\beg{align*}
\|M_r^{t}(x)\|_2\leq \ff {C_{L_1,T}} {\sqrt{r-t}}\left(1+\ff 1{\left(|x^{(1)}|^2+(r-t)\right)^{\ff \al 2}} \right).
\end{align*}
Then, we have
\beg{align*}
|\nn V(t,x,\P_\xi)|&\leq C_{L_1,L_2,T}\Bigg[\ff 1{\sqrt{T-t}}\left(1+\ff 1{\left(|x^{(1)}|^2+(T-t)\right)^{\ff \al 2}} \right)\cr
&\qquad\qquad\quad+\int_t^T\ff 1 {\sqrt{r-t}}\left(1+\ff 1{\left(|x^{(1)}|^2+(r-t)\right)^{\ff \al 2}} \right)\d r\Bigg](1+|x|+\|\xi\|_2),
\end{align*}
provided that $x^{(1)}\neq0$.

(iv) For PDEs \eqref{3.6} associated with a McKean-Vlasov FBSDE with forward stochastic Hamiltonian system,
in light of \eqref{Add3-Ex4} in Example \ref{ex.Ham} there will appear
\beg{align*}
\int_t^T\|M_r^{t}(x)\|_2\d r\leq C_{B,\tld{b},\si,T}\int_t^T\left(\sqrt{r-t}+\ff 1{(r-t)^{\ff 3 2}}\right)\d r=\infty.
\end{align*}
So, we cannot be sure that the left-hand side is finite and thus cannot directly adopt the above arguments to obtain an efficient estimate in this setting.
\end{rem}

\section{Appendix: proofs of auxiliary lemmas}

\subsection{Proof of Lemma \ref{rep-YZ-X}}

\beg{proof} We split the proof into three steps.

\textit{Step 1: Claim: $u(t,\cdot)$ is continuously differentiable for every $t\in[0,T]$} and satisfies \eqref{equ-nu0}.
Let $(t,x)\in[0,T]\times\R^d$ be fixed.
Owing to (A1), (A2) and the moment conditions \eqref{mon-est-1} and \eqref{mon-est-2} in (C1),
the following equations have a unique solution $(Y_s^{t,x},Z^{t,x}_s)_{s\in[t,T]}$ and $(\nn Y_s^{t,x},\nn Z^{t,x}_s)_{s\in[t,T]}$:
\beg{align}\label{add5-(Lem2.5)}
Y_s^{ t,x}&=g (X_T^{ t,x})+\int_s^T f(r,\Th_r^{ t,x })\d r-\int_s^T Z_r^{ t,x }\d W_r,\\
\nn Y_s^{ t,x}&=\nn g (X_T^{ t,x })\nn X_T^{ t,x }+\int_s^T \nn f (r,\Th_r^{ t,x })\nn\Th_r^{ t,x }\d r-\int_s^T \nn Z_r^{ t,x }\d W_r.\label{new-Y}
\end{align}
Moreover, as in the proof of \cite[Theorem 3.1, Page 1398-1399]{MaZ}, one can show that for any $v\in\R^d$,
\beg{align*}
\lim_{\ep\ra 0^+}\E\left(\sup_{s\in [t,T]}\left|\ff {Y_s^{t,x+\ep v}-Y_s^{t,x}} {\ep}-\nn_v Y_s^{t,x}\right|^2+\int_t^T\left\|\ff {Z_r^{t,x+\ep v}-Z_r^{t,x}} {\ep}-\nn_v Z_r^{t,x}\right\|^2\d r\right)=0.
\end{align*}
That is, $u(t,\cdot)$ is differentiable along the direction $v$ and $\nn_v u(t,x)=\nn_v Y_t^{t,x}$.

Next, we focus on proving that $\nn Y_t^{t,\cdot}$ is continuous.
To the end, applying the It\^o formula to \eqref{add5-(Lem2.5)} and using (A1)-(A2), we first obtain
\beg{align*}
&|Y_s^{t,y}-Y_s^{t,x}|^2+\int_s^T \|Z_r^{t,y}-Z_r^{t,x}\|^2\d r\\
&=|g(X_T^{t,y})-g(X_T^{t,x})|^2+2\int_s^T\<f(r,\Th_r^{t,y})-f(r,\Th_r^{t,x}),Y_r^{t,y}-Y_r^{t,x}\>\d r\\
&\quad-2\int_s^T\<Y_r^{t,y}-Y_r^{t,x},(Z_r^{t,y}-Z_r^{t,x})\d W_r\>\\
&\leq K_1^2|X_T^{t,y}-X_T^{t,x}|^2+2K_2\int_s^T|\Th_r^{t,y}-\Th_r^{t,x}|\cdot|Y_r^{t,y}-Y_r^{t,x}|\d r\\
&\quad+2\left|\int_s^T\<Y_r^{t,y}-Y_r^{t,x},(Z_r^{t,y}-Z_r^{t,x})\d W_r\>\right|.
\end{align*}
By the H\"{o}lder inequality, the Gronwall inequality and the B-D-G inequality, we deduce that there exists a constant $C>0$ such that
\beg{align*}
&\qquad\E \left(\sup_{s\in [t,T]}|Y_s^{t,y}-Y_s^{t,x}|^2+\int_t^T \|Z_r^{t,y}-Z_r^{t,x}\|^2\d r\right)\\
&\qquad \leq C\left(\E|X_T^{t,y}-X_T^{t,x}|^2+\E\int_t^T\left|X_r^{t,y}-X_r^{t,x}\right|^2\d r\right).
\end{align*}
Observe that from \eqref{mon-est-2} (with $k_1=2$) and \eqref{add-der} of (C1), it follows that for each $r\in[t,T]$,
\beg{align}\label{add1-(Lem2.5)}
\E|X_r^{t,y}-X_r^{t,x}|^2=\E\left|\int_0^{|y-x|}\nn_v X_r^{t,x+\th v}\d\th\right|^2
&\leq |y-x|\int_0^{|y-x|}\E|\nn_v X_r^{t,x+\th v}|^2\d\th\cr
&\leq C\left[1+(|x|+|y-x|)^{k_3}\right]^2|y-x|^2,\cr
\end{align}
where $v:=\ff {y-x}{|y-x|}$.
Then we have
\beg{align}\label{add2-(Lem2.5)}
\lim_{y\ra x}\E \left(\sup_{s\in [t,T]}|Y_s^{t,y}-Y_s^{t,x}|^2+\int_t^T\|Z_r^{t,y}-Z_r^{t,x}\|^2\d r\right)=0.
\end{align}
Similarly, applying the It\^o formula to \eqref{new-Y} and using (A1)-(A2) again, we get
\beg{align*}
&\|\nn Y_s^{t,y}-\nn Y_s^{t,x}\|^2+\int_s^T\|\nn Z_r^{t,y}-\nn Z_r^{t,x}\|^2\d r\\
&=\|\nn g(X_T^{t,y})\nn X_T^{t,y}-\nn g(X_T^{t,x})\nn X_T^{t,x}\|^2\\
& \quad +2\int_s^T \left\<\nn f(r,\Th_r^{t,y})\nn \Th_r^{t,y}-\nn f(r,\Th_r^{ t,x })\nn \Th_r^{ t,x },\nn Y_r^{t,y}-\nn Y_r^{t,x}\right\>\d r\\
& \quad -2\int_s^T\<\nn Y_r^{t,y}-\nn Y_r^{t,x}, \left(\nn Z_r^{t,y}-\nn Z_r^{t,x}\right)\d W_r\>\\
& \leq 2K_1^2\|\nn X_T^{t,y}-\nn X_T^{t,x}\|^2+2\|\nn g(X_T^{t,y})-\nn g(X_T^{t,x})\|^2\cdot\|\nn X_T^{t,x}\|^2\\
& \quad +2K_2\int_s^T\|\nn \Th_r^{t,y}-\nn \Th_r^{t,x}\|\cdot\|\nn Y_r^{t,y}-\nn Y_r^{t,x}\|\d r\\
& \quad +2\int_s^T\|\nn f(r,\Th^{t,y}_r)-\nn f(r,\Th_r^{t,x})\|\cdot\|\nn\Th_r^{t,x}\|\cdot\|\nn Y_r^{t,y}-\nn Y_r^{t,x}\|\d r\\
& \quad +2\left|\int_s^T\<\nn Y_r^{t,y}-\nn Y_r^{t,x}, \left(\nn Z_r^{t,y}-\nn Z_r^{t,x}\right)\d W_r\>\right|.
\end{align*}
Then, with the help of the H\"{o}lder, Gronwall and B-D-G inequalities again,
we derive that there exists a constant $C>0$ such that
\beg{align}\label{add0-(Lem2.5)}
&\E\left(\sup_{s\in[t,T]}\|\nn Y_s^{t,y}-\nn Y_s^{t,x}\|^2+\int_t^T\|\nn Z_r^{t,y}-\nn Z_r^{t,x}\|^2\d r\right)\cr
&\leq C\bigg[\E\|\nn X_T^{t,y}-\nn X_T^{t,x}\|^2+\E\left(\|\nn g(X_T^{t,y})-\nn g(X_T^{t,x})\|^2\cdot\|\nn X_T^{t,x}\|^2\right)\cr
&\quad\quad\quad +\E\int_t^T\|\nn X_r^{t,y}-\nn X_r^{t,x}\|^2\d r\cr
&\quad\quad\quad +\E\int_t^T\|\nn f(r,\Th_r^{t,y})-\nn f(r,\Th_r^{t,x}) \|^2\cdot\|\nn\Th_r^{t,x}\|^2\d r\bigg].
\end{align}
By \eqref{mon-est-2} and \eqref{con-nn-P} of (C1), the dominated convergence theorem implies
\beg{align}\label{add3-(Lem2.5)}
\lim_{y\ra x}\left[\E\|\nn X_T^{t,y}-\nn X_T^{t,x}\|^2+\E\int_t^T\|\nn X_r^{t,y}-\nn X_r^{t,x}\|^2\d r\right]=0.
\end{align}
In view of \eqref{add1-(Lem2.5)} and \eqref{add2-(Lem2.5)}, we obtain that as $y$ goes to $x$,
$X_T^{t,y}$ and $\Th_r^{t,y}$ converge to $X_T^{t,x}$ and $\Th_r^{t,x}$ in measures $\d\P$ and $\d r\otimes\d\P$, respectively.
Then, since $\nn g(\cdot)$ and $\nn f(r,\cdot,\cdot,\cdot)$ are continuous, we have
$$\P\text{-}\lim_{y\ra x}\|\nn g(X_T^{t,y})-\nn g(X_T^{t,x}) \|^2=0$$
and
$$\d r\otimes\d \P\text{-}\lim_{y\ra x}\|\nn f(r,\Th_r^{t,y})-\nn f(r,\Th_r^{t,x}) \|^2=0.$$
Consequently, due to the fact that $\nn g$ and $\nn f$ are bounded,
we can apply the dominated convergence theorem to get
\beg{align*}
&\lim_{y\ra x}\bigg[\E\left(\|\nn g(X_T^{t,y})-\nn g(X_T^{t,x})\|^2\cdot\|\nn X_T^{t,x}\|^2\right)\cr
&\qquad\quad+\E\int_t^T\|\nn f(r,\Th_r^{t,y})-\nn f(r,\Th_r^{t,x}) \|^2\cdot\|\nn\Th_r^{t,x}\|^2\d r\bigg]=0.
\end{align*}
Plugging this and \eqref{add3-(Lem2.5)} into \eqref{add0-(Lem2.5)}, we obtain
\beg{align}\label{add4-(Lem2.5)}
\lim_{y\ra x}\E\left(\sup_{s\in[t,T]}\|\nn Y_s^{t,y}-\nn Y_s^{t,x}\|^2+\int_t^T\|\nn Z_r^{t,y}-\nn Z_r^{t,x}\|^2\d r\right)=0,
\end{align}
which implies the continuity of $\nn Y_t^{t,\cdot}$.
We therefore conclude that $\nn u(t,\cdot)$ is continuous and $\nn u(t,x)=\nn Y_t^{t,x}$.

Note that $\nn Y_t^{t,x}$ is deterministic.
Then letting $s=t$ and taking expectation in \eqref{new-Y} yield \eqref{equ-nu0}.

\textit{Step 2: Claim: \eqref{ine-nu0} holds}.
Applying the It\^{o} formula to \eqref{new-Y}, we have
\beg{align*}
&\|\nn Y_s^{t,x}\|^2+\int_s^T\|\nn Z_r^{t,x}\|^2\d r+2\int_s^T\<\nn Y_r^{t,x}, \nn Z_r^{ t,x }\d W_r\>\\
& = \|\nn g (X_T^{ t,x })\nn X_T^{ t,x }\|^2+2\int_s^T\<\nn f (r,\Th_r^{ t,x })\nn\Th_r^{ t,x },\nn Y_r^{t,x}\>\d r\\
& \leq K_1^2\|\nn X_T^{ t,x }\|^2+2K_2\int_s^T \|\nn \Th_r^{t,x}\|\cdot\|\nn Y_r^{t,x}\|\d r,
\end{align*}
where we have used the assumptions (A1) and (A2).
Then it follows from the H\"older, Gronwall and B-D-G inequalities that
\beg{align*}
&\E\left(\sup_{s\in [t,T]}\|\nn Y_s^{t,x}\|^2+\int_t^T\|\nn Z_r^{t,x}\|^2\d r\right)\\
&\leq C_{K_1,K_2,T}\left(\E\|\nn X_T^{ t,x }\|^2+\int_t^T\E\|\nn X_r^{t,x}\|^2\d r\right)
\end{align*}
for some constant $C_{K_1,K_2,T}>0$.
Combining this with $\nn u(t,x)=\nn Y_t^{t,x}$ and recalling that $\nn Y_t^{t,x}$ is deterministic,
we obtain that there exists some constant $q_1>0$ such that
\beg{align*}
\|\nn u(t,x)\|&\leq   C_{K_1,K_2,T}\left[\left(\E\|\nn X_T^{ t,x }\|^2\right)^{\ff 1 2}+\left(\int_t^T\E \|\nn X_r^{t,x}\|^2\d r\right)^{\ff 1 2}\right]\\
&\leq C_{K_1,K_2,k_3,T}(1+|x|^{q_1}),\ \ x\in\R^d,
\end{align*}
where the last inequality is due to \eqref{mon-est-2} (with $k_1=2$ and $q_1=k_3$) of (C1).
This yields the desired assertion.


\textit{Step 3: Claim: $Y_s^{t,x}=u(s,X_s^{t,x})$ and $Z_s^{t,x}=\nn u(s,X_s^{t,x})\si(s,X_s^{t,x}), \d s\otimes\d\P$-a.e.}
For $x_n$ and $(X_s^{n,t,x_n},\nn X_s^{n,t,x_n})_{s\in[t,T]}$ as in (C2),
let $(Y_s^{n,t,x_n},Z_s^{n,t,x_n})_{s\in[t,T]}$ be the solution of the following equation:
\beg{align*}
Y_s^{n,t,x_n}&=g (X_T^{n,t,x_n})+\int_s^T f (r,\Th_r^{n,t,x_n})\d r-\int_s^T Z_r^{n,t,x_n}\d W_r,
\end{align*}
where $\Th_\cdot^{n,t,x_n}:=(X_\cdot^{n,t,x_n},Y_\cdot^{n,t,x_n},Z_\cdot^{n,t,x_n})$.
Since $g(\cdot)$ and $f(r,\cdot,\cdot,\cdot)$ are continuously differentiable, we arrive at
\beg{align*}
\nn Y_s^{n,t,x_n}&=\nn g (X_T^{n,t,x_n})\nn X_T^{n,t,x_n}+\int_s^T\nn f (r,\Th_r^{n,t,x_n})\nn\Th_r^{n,t,x_n}\d r- \int_s^T \nn Z_r^{n,t,x_n}\d W_r.
\end{align*}
Using \eqref{app-X1} and \eqref{app-X2} of (C2) and following the same arguments as for \eqref{add2-(Lem2.5)} and \eqref{add4-(Lem2.5)},
we have
\beg{align}\label{app-Y1}
\lim_{n\ra+\infty}\E \left(\sup_{s\in [t,T]}|Y_s^{n,t,x_n}-Y_s^{t,x}|^2+\int_t^T\|Z_r^{n,t,x_n}-Z_r^{t,x}\|^2\d r\right)=0
\end{align}
and
\beg{align}\label{app-nnY}
\lim_{n\ra+\infty}\E\left(\sup_{s\in[t,T]}\|\nn Y_s^{n,t,x_n}-\nn Y_s^{t,x}\|^2+ \int_t^T\|\nn Z_r^{n,t,x_n}-\nn Z_r^{t,x}\|^2\d r\right)=0.
\end{align}
We now set $u_n(t,x_n):=Y_t^{n,t,x_n}$.
According to \cite[Theorem 3.1 and Remeark 3.3]{ZJF}, we obtain that $u_n\in C^{0,1}([0,T]\times\R^d;\R^l)$ and
\beg{align}\label{Add6-(Lem2.5)}
Y_s^{n,t,x_n}=u_n(s,X_s^{n,t,x_n}),\ \ Z_s^{n,t,x_n}=\nn u_n(s,X_s^{n,t,x_n})\si_n(s,X_s^{n,t,x_n}),\ \ s\in[t,T].
\end{align}
Then by \eqref{app-Y1}, \eqref{app-nnY} and Step 1, we have that for any $t\in [0,T]$,
\beg{align}\label{Add7-(Lem2.5)}
\lim_{n\ra+\infty}u_n(t,x_n)&=\lim_{n\ra+\infty}Y_t^{n,t,x_n}=Y_t^{t,x}=u(t,x)
\end{align}
and
\beg{align}\label{Add8-(Lem2.5)}
\lim_{n\ra+\infty}\nn u_n(t,x_n)&=\lim_{n\ra+\infty}\nn Y_t^{n,t,x_n}=\nn Y_t^{t,x}=\nn u(t,x).
\end{align}
Observe that from \eqref{app-X1} of (C2), there exists a subsequence of $X_s^{n,t,x_n}$, which we still denote by $ X_s^{n,t,x_n}$,  such that
$$\lim_{n\ra +\infty}X_s^{n,t,x_n}=X_s^{t,x},\ \ \d s\otimes \d\P\text{-a.e.,}$$
and moreover from \eqref{app-si} of (C2), it follows that
\beg{align*}
\lim_{n\ra+\infty}\si_n(s,X_s^{n,t,x_n})=\si(s,X_s^{t,x}),\ \ \d s\otimes\d\P\text{-a.e}.
\end{align*}
Consequently, by \eqref{Add7-(Lem2.5)} and \eqref{Add8-(Lem2.5)} we derive that $\d s\otimes \d\P$-a.e.,
\beg{align*}
\lim_{n\ra+\infty}u_n(s,X_s^{n,t,x_n})=u(s,X_s^{t,x})
\end{align*}
and
\beg{align*}
\lim_{n\ra+\infty}\nn u_n(s,X_s^{n,t,x_n})\si_n(s,X_s^{n,t,x_n})=\nn u(s,X_s^{t,x})\si(s,X_s^{t,x}).
\end{align*}
Note that by \eqref{app-Y1} and selecting a subsequence if necessary, we obtain that for each $s\in [t,T]$,
$\lim_{n\ra+\infty}Y_s^{n,t,x_n}=Y_s^{t,x}$, $\P$-a.s. and $\lim_{n\ra+\infty}Z_s^{n,t,x_n}=Z_s^{t,x}$, $\d s\otimes \d\P$-a.e.
Hence, because of \eqref{Add6-(Lem2.5)}, it is easy to see that the desired relations hold.
This completes the proof.
\end{proof}

\subsection{Proof of Lemma \ref{lem-SI}}

\beg{proof}
Note that for each $0\leq t\leq s<r\leq T$,
$\tld W([s,r])$ is independent of $\sF_s^t$ and $X_s^{t,x}\in \sF_s^t$.
Then by \eqref{equ-M}, we have for any $\psi\in \sB_b(\R^d)$,
\beg{align*}
&\E\left[\ps(X_r^{t,x})M_r^s(X_s^{t,x})\big|\sF_s^t\right]\cr
&=\E\left[ \ps(X_r^{t,x})\int_s^r H_{s,r}\left(\th, X_s^{t,x}, \tld W([s,r]),X_\th^{s,X_s^{t,x}},\nn X_\th^{s,y}\Big|_{y=X_s^{t,x}}\right)\d W_\th \Big| \sF_s^t\right]\cr
&= \E\left[ \ps(X_r^{s,X_s^{t,x}}) \int_s^r H_{s,r}\left(\th, X_s^{t,x}, \tld W([s,r]), X_\theta^{s,X_s^{t,x}},\nn X_\th^{s,y}\Big|_{y=X_s^{t,x}}\right)\d W_\th \Big| \sF_s^t\right]\cr
&= \left[E\left(\ps(X_r^{s,y}) \int_s^r H_{s,r}(\th, y, \tld  W([s,r]),X_\th^{s,y},\nn X_\th^{s,y})\d W_\th\right) \right]\Big|_{y=X_s^{t,x}}\cr
&= \left[\E\left(\ps(X_r^{s,y})M_r^s(y)\right) \right]\big|_{y=X_s^{t,x}}.
\end{align*}
Observe that all equalities above still hold true for  $|M_r^s(X_s^{t,x})|$ replacing $M_r^s(X_s^{t,x})$.
Hence \eqref{equ-con} holds for each $\ps\in \sB_b(\R^d)$ $\P$-a.s.

Let us now investigate the function $K(t,s,r,\cdot)$.
According to \eqref{equ-M} and \eqref{ine-HH}, we obtain
\beg{align}\label{Add0-Lem2.6}
&K^2(t,s,r,x)=\E \left|M_r^{s}(X_s^{t,x})\right|^2 \cr
&=\E\left|\int_s^r H_{s,r}\left(\th, X_s^{t,x}, \tld W([s,r]), X_\th^{s,X_s^{t,x}},\nn X_\th^{s,y}\Big|_{y=X_s^{t,x}}\right)\d W_\th\right|^2 \cr
& =\E\int_s^r \left\|H_{s,r}\left(\th, X_s^{t,x}, \tld W([s,r]), X_\th^{t,x},\nn X_\th^{s,y}\Big|_{y=X_s^{t,x}}\right)\right\|^2\d \th\cr
& \leq 3C^2_{q_1,q_2,T}\E\int_s^r \left( K^2_\th(t,s,r,x)+|X_\th^{t,x}|^{2q_1}+\left\|\nn X_\th^{s,y}\Big|_{y=X_s^{t,x}}\right\|^{2q_2}\right) \d \th\cr
&\leq C_{q_1,q_2,k_2,k_3,T} \left( \E\int_s^r K^2_\th(t,s,r,x)\d \th+1+|x|^{2q_1k_2}+\E\left|X_s^{t,x}\right|^{2q_2k_3}\right)\cr
& \leq C_{q_1,q_2,k_2,k_3,T}\left(\E\int_s^r K^2_\th(t,s,r,x)\d \th+1+|x|^{2q_1k_2}+|x|^{2q_2k_3k_2}\right).
\end{align}
Here, we have used the moment conditions \eqref{mon-est-1} and \eqref{mon-est-2} in the last two inequalities.
Due to the assumption that the term $\E\int_s^r K^2_\th(t,s,r,\cdot)\d \th$ is of polynomial growth,
so is $K(t,s,r,\cdot)$.
Finally, it is easily seen that \eqref{iint-K} follows from \eqref{Add0-Lem2.6} and \eqref{int-K-1}, which completes the proof.
\end{proof}

\textbf{Acknowledgement}

X. Fan and M. R\"{o}ckner are grateful to the financial support by the DFG
through the CRC 1283 Taming uncertainty and profiting from randomness and low
regularity in analysis, stochastics and their applications.
X. Fan is partially supported by the Natural Science Foundation of Anhui Province (No. 2008085MA10).
S.-Q. Zhang is supported in part by the National Natural Science Foundation of China (No. 11771326, 11901604).

\end{document}